\documentclass[12pt]{article}
\usepackage{amsmath,amssymb,amsfonts,amsthm}
\usepackage[all]{xy}

\newcounter{item}[section]
\newcounter{kirshr}
\newcounter{kirsha}
\newcounter{kirshb}
\newenvironment{enumroman}{\setcounter{kirshr}{1}
\begin{list}{(\roman{kirshr})}{\usecounter{kirshr}} }{\end{list}}
\newenvironment{enumarab}{\setcounter{kirshb}{1}
\begin{list}{(\arabic{kirshb})}{\usecounter{kirshb}} }{\end{list}}
\newenvironment{athm}[1]{\vskip3mm\par\noindent
{\bf #1 }. \slshape }
{\upshape\par\vskip10pt minus3pt}
\newtheorem{theorem}{Theorem}[section]

\newtheorem{lemma}[theorem]{Lemma}
\newtheorem{corollary}[theorem]{Corollary}
\newenvironment{demo}[1]{\noindent{\bf #1.}\upshape\mdseries}
{\nopagebreak{\hfill\rule{2mm}{2mm}\nopagebreak}\par\normalfont}
\theoremstyle{definition}
\newtheorem{remark}[theorem]{Remark}

\newtheorem{example}[theorem]{Example}
\newtheorem{definition}[theorem]{Definition}

\def\C{{\mathfrak{C}}}
\def\Fm{{\mathfrak{Fm}}}

\def\N{{\cal N}}

\def\Nr{{\mathfrak{Nr}}}
\def\Fr{{\mathfrak{Fr}}}
\def\Sg{{\mathfrak{Sg}}}
\def\Fm{{\mathfrak{Fm}}}
\def\A{{\mathfrak{A}}}
\def\B{{\mathfrak{B}}}
\def\C{{\mathfrak{C}}}
\def\D{{\mathfrak{D}}}

\def\Rd{{\mathfrak{Rd}}}
\def\(R)RA{{\bf (R)RA}}

\def\B{{\sf B}}

\def\Nr{{\mathfrak{Nr}}}

\def\Nr{{\mathfrak{Nr}}}

\def\A{{\mathfrak{A}}}
\def\B{{\mathfrak{B}}}
\def\C{{\mathfrak{C}}}
\def\D{{\mathfrak{D}}}

\def\Nr{{\mathfrak{Nr}}}
\def\Fr{{\mathfrak{Fr}}}
\def\Sg{{\mathfrak{Sg}}}
\def\Fm{{\mathfrak{Fm}}}
\def\Rd{{\mathfrak{Rd}}}

\def\(R)RA{{\bf (R)RA}}

\def\Nr{{\mathfrak{Nr}}}
\def\Fr{{\mathfrak{Fr}}}
\def\Sg{{\mathfrak{Sg}}}
\def\Fm{{\mathfrak{Fm}}}
\def\Rd{{\mathfrak{Rd}}}

\def\(R)RA{{\bf (R)RA}}

\def\N{\mathbb{N}}

\def\A{{\mathfrak{A}}}
\def\B{{\mathfrak{B}}}
\def\C{{\mathfrak{C}}}
\def\D{{\mathfrak{D}}}

\def\Fm{{\mathfrak{Fm}}}

\def\Nr{{\mathfrak{Nr}}}
\def\F{{\mathfrak{F}}}

\def\Sg{{\mathfrak Sg}}

\title{Representability, and amalgamation for  various reducts of Heyting polyadic algebras}
\author{Tarek Sayed Ahmed}

\begin{document}
\maketitle
\begin{abstract}

We study algebraically infinitely many infinitray extensions of predicate intuitionistic logic. 
We prove several representation theorems 
that reflect a  (weak) Robinson's joint consistency theorem for the extensions studied with and without equality. 
In essence a Henkin-Gabbay construction, our proof uses neat embedding theorems and is purely algebraic.
Neat embedding theorems, are an algebraic version of Henkin constructions that apply to various 
infinitary extensions of predicate first order logics; to the best of our knowledge, 
they were only  implemented in the 
realm of intuitionistic logic in the article 'Amalgamation of polyadic Heyting algebras' 
Studia Math Hungarica, in press.


\end{abstract}

\section{Introduction}

\subsection{Background and History}

It often happens that a theory designed originally as a tool for the
study of a problem, say in computer science, came subsequently to have purely
mathematical interest. When  such a phenomena occurs, the theory is usually
generalized beyond the point needed for applications, the
generalizations make contact with other theories (frequently in
completely unexpected directions), and the subject becomes
established as a new part of pure mathematics. The part of pure
mathematics so created does not (and need not) pretend to solve the
problem from which it arises; it must stand or fall on its
own merits.

A crucial addition to the
collection of mathematical catalysts initiated at the beginning of the 20 century, is formal logic and its study 
using mathematical machinery,
better known as  metamathematical investigations, or simply metamathematics. 
Traced back to the works of Frege, Hilbert, Russel, Tarski, Godel and others; 
one of the branches of pure mathematics that metamathematics  has precipitated to is algebraic logic.

Algebraic logic is an interdisciplinary field; it is the art of tackling problems in formal logic using universal algebraic machinery.
It is similar in this respect to several branches in mathematics, like algebraic geometry, 
where algebraic machinery is used guided by geometric intuition.
In algebraic logic, the intuition underlying its constructions is inspired from (mathematical) logic.

The idea of solving problems in various  branches of logic by first translating them to
algebra, then using the powerful methodology of algebra for solving
them, and then translating the solution back to logic, goes back to
Leibnitz and Pascal. Such a methodology was already fruitfully applied back  in the 19th century
with the work of Boole, De Morgan, Peirce, Schr\"{o}der, and others on classical logic. 
Taking logical equivalence rather than truth as
the primitive logical predicate and exploiting the similarity
between logical equivalence and equality, those pioneers developed logical
systems in which metalogical investigations take on a plainly
algebraic character. The ingenious  transfer of ''logical equivalence" to '' equations" turned out immensely useful and fruitful. 

In particular, Boole's work evolved into the
modern theory of Boolean algebras, and that of De Morgan, Peirce and
Schr\"{o}der into the well-developed theory of relation algebras, 
which is now widely used in such diverse areas, ranging from formalizations of set theory to applications 
in computer science.

From the beginning
of the contemporary era of logic, there were  two approaches to the subject, one
centered on the notion of logical equivalence and the other, reinforced by Hilbert's work on metamathematics, centered
on the notions of assertion and inference.

It was not until much later that logicians started to think about
connections between these two ways of looking at logic. Tarski
gave the precise connection between Boolean algebra
and the classical propositional calculus, inspired by the impressive work of Stone on Boolean algebras.
Tarski's approach builds on
Lindenbaum's idea of viewing the set of formulas as an algebra with
operations induced by the logical connectives.
When the Lindenbaum-Tarski method is applied to the predicate
calculus, it lends itself to cylindric and polyadic algebras rather than
relation algebras.

In the traditional  mid -20th century approach, algebraic logic has focused on
the algebraic investigation of particular classes of algebras 
like cylindric, polyadic and relation algebras. When such a connection could be established, there was
interest in investigating the interconnections  between various
metalogical properties of the logical system in question and the algebraic
properties of the coresponding class of algebras (obtaining what are
sometimes called "bridge theorems"). This branch has now evolved into the relatively new field  of universal algebraic logic, 
in analogy to the well established field 
of universal algebra. 

For example, it was discovered
that there is a natural relation between the interpolation theorems
of classical, intuitionistic, intermediate propositional calculi,
and the amalgamation properties of varieties of Heyting algebras, which constitute the  main focus of  this paper.
The variety  of Heyting algebras is the algebraic counterpart of propositional intuitionistic logic.
We shall deal with Heyting algebras with extra (polyadic) operations reflecting 
quantifiers.  Those algebras are appropriate to study (extensions) of predicate intuitionistic logic. Proving various interpolation theorems 
for such extensions,  we thereby 
extend known amalgamation results of Heyting algebras to  polyadic expansions. 

A historic comment on the development of intuitioinistic logic is in order. 
It was Brouwer who first initiated the programme of intuitionism and intuitionistic logic is its rigorous formalization developed originaly by Arend Heyting. 
Brouwer rejected formalism per se but admitted  the potential usefulness of formulating general logical principles expressing 
intuitionistically correct constructions, such as modus ponens. Heyting realized the importance of formalization, being fashionable at his time, with the rapid development of 
mathematics.
Implementing intuitionistic logic,  turned out   useful for diffrent  forms of mathematical constructivism since it has the existing property.
Philosophically, intuitionism differs from  logicism by treating logic as an independent branch 
of mathematics, rather than as the foundations of mathematics,  
from finitism by permitting intuitionistic reasoning about possibly infinite collections; and from platonism 
by viewing mathematical objects as mental constructs rather than entities with an independent objective existence. 
There is also analogies between logisicm and intuitionism; in fact 
Hilbert's formalist program, aiming to base the whole of classical mathematics on solid foundations by reducing it to a 
huge formal system whose consistency should be established by 
finitistic, concrete (hence constructive) means, was the most powerful contemporary rival to Brouwer's and Heyting's  intuitionism.

\subsection {Subject Matter}

Connections between interpolation theorems
in the predicate calculus and amalgamation results in varieties of
cylindric and polyadic algebras, were initiated mainly by Comer, Pigozzi, Diagneault and Jonsson.

As it happened, during the course of the development of algebraic logic, dating back to the work of Boole, up 
to its comeback in the contemporary era through the pioneering work of Halmos, Tarski, Henkin, Monk, Andr\'eka, and N\'emeti, 
it is now established that 
the two most famous widely used algebraisations of first order logic are Tarski's cylindric algebras \cite{HMT1}, \cite{HMT2}, 
and Halmos' polyadic algebras \cite{Halmos}. Each has its advantages 
and disadvantages.
For example, the class of representable cylindric algebras, though a variety, is not finitely axiomatizable, 
and this class exhibits an inevitable degree of complexity in any of its axiomatizations \cite{Andreka}. 
However, its equational theory is recursive. 
On the other hand, the variety 
of (representable) polyadic 
algebras is axiomatized by a finite schema of equations but its  
equational theory is not recursively enumerable \cite{NS}. 
There have been investigations 
to find a class of algebras that enjoy the positive properties of both. The key idea behind such investigations is to look at (the continuum 
many) reducts of polyadic algebras \cite{AUamal}, \cite{S} searching for the desirable finitely axiomatizable variety among them. 

Indeed, it is folkore in algebraic logic that cylindric algebras and polyadic algebras belong to different paradigms, frequently manifesting contradictory
behaviour.
The paper \cite{S} is a unification of the positive properties  of those two paradigms in the Boolean case, and one of the results of 
this paper can be interpreted as a unification of 
those paradigms when the propositional reducts are Heyting algebras.

A polyadic algebra is typically an instance of a transformation system.
A transformation system can be defined to be a quadruple of the form $(\A, I, G, {\sf S})$ where $\A$ is an algebra of any similarity type, 
$I$ is a non empty set (we will only be concerned with infinite sets),
$G$ is a subsemigroup of $(^II,\circ)$ (the operation $\circ$ denotes composition of maps) 
and ${\sf S}$ is a homomorphism from $G$ to the semigroup of endomorphisms of $\A$ $(End(\A))$. 
Elements of $G$ are called transformations. 

The set $I$ is called the dimension of the algebra, for a transformation $\tau$ on $I$, ${\sf S}({\tau})\in End(\A)$ is called a substitution operator, or simply a substitution. 
Polyadic algebras arise when $\A$ is a Boolean algebra endowed with quantifiers and $G={}^II$. 
There is an extensive literature for polyadic algebras dating back to the fifties and sixties of the last century, 
\cite{Halmos}, \cite{J70}, \cite{D}, \cite{DM}, \cite{ AUamal}, \cite{S}. 
Introduced by Halmos, the theory of polyadic algebras is now picking up again; indeed it's
regaining momentum with pleasing progress and a plathora of results, see the references 
 \cite{MLQ}, \cite{Fer1}, \cite{Fer2}, \cite{Fer3}, \cite{Fer4}, \cite{ANS}, \cite{trans}, to name just a few.

In recent times reducts of polyadic algebras of dimension $I$ were studied \cite{S}, \cite{AUamal}; these reducts 
are obtained by restricting quantifiers to involve only quantification on finitely many 
variables and to study 
(proper) subsemigroups of $^II$ The two extremes are the semigroup of finite transformations 
(a finite transformation is one that moves only finitely many points)
and all of $^II$ but there are infinitely many semigroups in between.

In this paper, we study reducts of polyadic algebras by allowing (proper) subsemigroups of $^II$, 
but we also weaken the Boolean structure  to be a Heyting algebra. Thus we approach  the realm of intuitionistic logic. 
We shall study the cases when $G$ consists of all finite transformations, when $G$ 
is a proper subsemigroup satisfying certain properties but essentially containing infinitary transformations, that is, 
transformations that move infinitely many points
(this involves infinitely many cases), 
and when $G$ is the semigroup of all transformations. Our investigations 
will address the representation of such algebras in a concrete sense
where the operations are interpreted as set-theoretic operations on sets of sequences, and will also address
the amalgamation property and variants thereof of the classes in question. 

In all the cases we study, the scope of quantifiers are finite, so in this respect our algebras also 
resemble cylindric algebras. The interaction between the theories of Boolean cylindric algebras and Boolean polyadic algebras is extensively 
studied in algebraic logic, see e.g
\cite{ANS}, with differences and similarities illuminating both theories.
In fact, the study of $G$ Boolean polyadic algebras ($G$ a semigroup) by Sain in her pioneering paper \cite{S}, 
and its follow up \cite{AUamal}, is an outcome, 
or rather a culmination, of such research; 
it's a typical situation in which the positive properties of both theories amalgamate.

Boolean polyadic algebras, when $G$ is the set of finite transformations of $I$ or $G={}^II$ are old \cite{Halmos}, \cite{D} 
\cite{DM}. In the former case such algebras are known as quasipolyadic 
algebras, and those  are substantially different from full  polyadic algebras (in the infinite dimensional case), as is commonly accepted, 
quasipolyadic algebras
belong to the cylindric paradigm; they share a lot of properties of cylindric algebras. 
While the substitution operators in full Heyting polyadic algebras are uncountable, even if both  the algebra and its dimension 
are  countable, the substitution operators for quasipolyadic equality algebras of countable dimension are countable. 
Unlike full polyadic algebras, quasipolyadic algebras can be formulated as what is known in the literature as a system of varieties 
definable by schemes making them  akin to universal algebraic investigations 
in the spirit of cylindric algebras. Though polyadic algebras can be viewed as a system of varieties, this system cannot be definable by schemes due to the presence of
infinitary substitutions.  
Studying reducts of polyadic algebras by allowing only 
those substitutions coming from an arbitrary  subsemigroup of $^II$ is relatively recent starting at the turn of the last century
\cite{S}.

Such  algebras (of which we study their Heyting reducts)  also provide a possible  solution to a central 
problem in algebraic logic, better known as the finitizability problem,  which asks for a simple (hopefully) finite axiomatization 
for several classes of representable algebras that abound in algebraic logic.  
\footnote {The class of representable algebras is given by specifying the universes of the algebras in the 
class, as sets of certain sets
endowed with set theoretic concrete operations; thus representable algebras are completely determined once one specifies their universes.}
The finitizability problem is not easy, and has been 
discussed at length  in the literature \cite{Bulletin}. Being rather a family of problems, the finitizability problem  
has several scattered reincarnations in the lierature, and in some sense is still open. 
The finitizability problem also has philosophical implications, repercussions, connotations, concerning reasoning about reasoning, 
and can, in so many respects, be likened to Hilbert's programe
of proving mathematics consistent by concrete finitistic methods. 

In fact, our results show that, 
when $G$ satisfies some conditions that are not particularly complicated,
provides us with an algebraisable extension of predicate first order intuitionistic logic, 
whose algebraic counterpart is a variety that is finitely axiomatizable.
 An algebraisable extension is an extension of ordinary predicate intuitionistic logic (allowing formulas of infinite length), whose algebraic counterpart 
consisting of subdirect products of set algebras based on (Kripke) models,
is a finitely based variety (equational class).
This gives a clean cut  solution to the analogue of the finitizability problem for ordinary predicate intuitionistic logic.

Formal systems for intuitionistic propositional and predicate logic and arithmetic were developed by Heyting \cite{H},\cite{Hy} 
Gentzen \cite{G} and Kleene \cite{K}. Godel \cite{Godel} proved the equiconsistency of intuitionistic and classical theories. 
Kripke \cite{Kripke} provided a semantics with respect to which intuitionistic logic is sound and complete. 
We shall use a modified version of 
Kripke semantics below to  prove our representability
results.

The algebraic counterparts of predicate intuitionistic logic, namely, Heyting polyadic algebras 
were studied by Monk \cite{Monk}, Georgescu \cite{G} and the present author \cite{Hung}. 
Algebraically, we shall prove that certain reducts of algebras of full polyadic Heyting algebras (studied in \cite{Hung}) 
consist solely of representable algebras (in a concrete sense) and have 
the superamalgamation property (a strong form of amalgamation). Such results are essentially 
proved in Part 1, with the superamalgmation property deferred to part 3.
We also present some negative 
results for other infinitary intiutionistic logics, based on non finite axiomatizability results proved in part 2, using bridge
theorems. Indeed, in part 3, among other things, 
we show that the minimal algebraisable extension of predicate intuitionistic logic, in a sense to be made precise, 
is essentially incomplete, 
and fail to have the interpolation property. 

Roughly, minimal extension here means this (algebraisable) logic corresponding to the variety generated by the class of algebras 
arising from ordinary intuitionistic predicate logic. 
Such  algebras are locally finite, reflecting the fact that formulas contain only finitely many variables. 
This correspondence is taken in the sense of Blok and Pigozzi associating quasivarieties to algebraisable logics. 
Algebraising here  essentially means that we drop the condition of local finiteness, 
(hence alowing formulas of infinite length); this property is not warranted from
the algebraic point of view because it is a poperty that cannot be expressed by first order formulas, let alone equations or quasiequations

In fact, we show that all positive results in this paper extend to the classical case, 
reproving deep results in \cite{S}, \cite{AUamal}, and many negative results
that conquer the cylindric paradigm, extend in some exact sense, to certain infinitary extensions of predicate intuitionistic logic, 
that arise naturally from the process
of algebraising the intuitionistic predicate logic ( with and without equality). Such results are presented in the context 
of clarifying one facet of 
the finitizability problem for predicate intuionistic logic, 
namely that of drawing a line between positive and negative results.

The techniques used in this  paper intersects those adopted in our recent paper on Heyting polyadic algebras; it uses 
this part of algebraic logic developed essentially by Henkin, Monk, Tarski and Halmos 
- together with deep techniques of Gabbay - 
but there are major differences. 

We mention two
\begin{enumarab}

\item Whereas the results in \cite{Hung} address full Heyting polyadic algebras where infinitary cylindrifications and infinitary substitutions 
are available; this paper, among many other things, shows that the proof survives when we restrict our attention to 
finitely generated semigroups still containing infinitary 
substitutions, and finite cylindrifiers.  The algebras in \cite{Hung} have an axiomatization that is highly complex from the recursion theoretic point of view.
The reducts studied here have recursive axiomatizations.

\item We allow diagonal elements in our algebras (these elements reflect equality), so in fact, 
we are in the realm of infinitary extensions of intuitionistic predicate logic {\it with} equality.




\end{enumarab} 

The interaction between algebraic logic and intuitionistic logic was developed in the monumental work of the Polish logicians 
Rasiowa and Sikorski, and the Russian logician Maksimova, 
but apart from that work, to the best of our knowledge, the surface of predicate intuitionistic logic was barely scratched by algebraic machinery. 
While Maksimova's work \cite{b} is  more focused on propositional intuitionistic logic, 
Rasiowa and Sikorski did deal with expansions of Heyting algebras, to reflect quantifiers, 
but not with polyadic algebras per se. Besides, Rasiowa and Sikorski, dealt only with classical predicate
intuitionistic logic. 

In this paper, we continue the trend initiated in \cite{Hung}, by studying strict reducts of full fledged infinitary logics, 
which are still infinitary, together with their expansions 
by the equality symbol, proving completenes theorems and interpolation properties, and we also maintain the borderline where such 
theorems cease to hold.

\subsection*{Organization}

In the following section we prepare for our algebraic proof, 
by formulating and proving the necessary algebraic preliminaries (be it concepts or theorems) 
addressing various reducts of Heyting polyadic algebras, possibly endowed with diagonal elements.

Our algebraic proof of the interpolation property for infinitary extensions of predicate intuitionistic logic (with and without
equality) are proved in section 3, which is the soul and heart of this part of the paper. This is accomplished
using the well developed methodology of algebraic logic; 
particularly so-called neat embedding theorems, which are  algebraic generalizations of 
Henkin constructions.

\subsection*{On the notation}
 
Throughout the paper, our notation is fairly standard, or self explanatory. 
However, we usually  distinguish notationally between algebras and their domains, though there will be occasions 
when we do not distinguish between the two.
Algebras will be denoted by Gothic letters, and when we write $\A$ for an algebra, then it is to be tacitly
assumed  that the corresponding Roman letter $A$ denotes 
its domain. Unfamiliar notation will be introduced at its first occurrence in the text. We extensively use the axiom of choice 
(any set can be well ordered, so that in many places we deal with ordinals or order types, that is, we impose well orders on arbitrary sets).
For a set $X$, $Id_X$, or simply $Id$, when $X$ is clear from context, denotes the identity function on $X$. 
The set of all functions from $X$ to $Y$ is denoted by $^XY$. If $f\in {}^XY$, then we write $f:X\to Y$. 
The domain $X$ of $f$, will be denoted by $Dof$, and the range of $f$ will be denoted by $Rgf$.
Composition of functions $f\circ g$ is defined so that the function at the right acts first, that is $(f\circ g)(x)=f(g(x))$, for $x\in Dog$ such that 
$g(x)\in Dof$.

\section{Algebraic Preliminaries}

In this section, we define our algebras and state and prove certain algebraic notions and properties that we shall 
need in our main (algebraic) proof implemented in the following section. 

Other results, formulated in lemmata \ref{dl} and \ref{cylindrify} 
in this section are non-trivial modifications of existing theorems for both cylindric algebras 
and polyadic algebras;  we give detailed  proofs of such results, 
skipping those parts that can be found in the literature referring to the necessary references instead. These lemmata address
a very important and key concept in both cylindric and polyadic theories, namely, 
that of forming dilations and neat reducts (which are, in fact, dual operations.)

\subsection{The algebras}

For an algebra $\A$, $End(\A)$ denotes the set of endomorphisms of $\A$ (homomorphisms of $\A$ into itself), 
which is a semigroup under the operation $\circ$ of composition of maps.

\begin{definition} A transformation system is a quadruple $(\A, I, G, {\sf S})$ where 
$\A$ is an algebra, $I$ is a set, $G$ is a subsemigroup of $(^II,\circ)$ 
and ${\sf S}$ is a homomorphism from
$G$ into $End(\A).$
\end{definition}
Throughout the paper, $\A$ will always be a Heyting algebra.
If we want to study predicate intuitionistic logic, then we are naturally led to expansions of Heyting algebras allowing quantification.
But we do not have negation in the classical sense, so we have to deal with existential and universal quantifiers each separately.

\begin{definition} Let $\A=(A, \lor, \land,\rightarrow,0)$ be a Heyting algebra. An existential quantifier $\exists$ on $A$ is a mapping
$\exists:\A\to \A$ such that the following hold for all $p,q\in A$: 
\begin{enumarab}
\item $\exists(0)=0,$
\item $p\leq \exists p,$
\item $\exists(p\land \exists q)=\exists p\land \exists q,$
\item $\exists(\exists p\rightarrow \exists q)=\exists p\rightarrow \exists q,$
\item $\exists(\exists p\lor \exists q)=\exists p\lor \exists q,$
\item $\exists\exists p=\exists p.$
\end{enumarab}
\end{definition}
\begin{definition}  Let $\A=(A, \lor, \land,\rightarrow,0)$ be a Heyting algebra. A universal quantifier  $\forall$ on $A$ is a mapping
$\forall:\A\to \A$ such that the following hold for all $p,q\in A$: 
\begin{enumarab}
\item $\forall 1=1,$
\item $\forall p\leq p,$
\item $\forall(p\rightarrow q)\leq \forall p\rightarrow \forall q,$
\item $\forall \forall p=\forall p.$
\end{enumarab}
\end{definition}

Now we define our algebras. Their similarity type depends on a fixed in advance semigroup. 
We write $X\subseteq_{\omega} Y$ to denote that $X$ is a finite subset 
of $Y$.
\begin{definition} Let $\alpha$ be an infinite set. Let $G\subseteq {}^{\alpha}\alpha$ be a semigroup under the operation of composition of maps. 
An $\alpha$ dimensional polyadic Heyting $G$ algebra, a $GPHA_{\alpha}$ for short, is an algebra of the following
form
$$(A,\lor,\land,\rightarrow, 0, {\sf s}_{\tau}, {\sf c}_{(J)}, {\sf q}_{(J)})_{\tau\in G, J\subseteq_{\omega} \alpha}$$
where $(A,\lor,\land, \rightarrow, 0)$ is a Heyting algebra, ${\sf s}_{\tau}:\A\to \A$ is an endomorphism of Heyting algebras,
${\sf c}_{(J)}$ is an existential quantifier, ${\sf q}_{(J)}$ is a universal quantifier, such that the following hold for all 
$p\in A$, $\sigma, \tau\in [G]$ and $J,J'\subseteq_{\omega} \alpha:$
\begin{enumarab}
\item ${\sf s}_{Id}p=p.$
\item ${\sf s}_{\sigma\circ \tau}p={\sf s}_{\sigma}{\sf s}_{\tau}p$ (so that ${\sf s}:\tau\mapsto {\sf s}_{\tau}$ defines a homomorphism from $G$ to $End(\A)$; 
that is $(A, \lor, \land, \to, 0, G, {\sf s})$ is a transformation system).
\item ${\sf c}_{(J\cup J')}p={\sf c}_{(J)}{\sf c}_{(J')}p , \ \  {\sf q}_{(J\cup J')}p={\sf q}_{(J)}{\sf c}_{(J')}p.$
\item ${\sf c}_{(J)}{\sf q}_{(J)}p={\sf q}_{(J)}p , \ \  {\sf q}_{(J)}{\sf c}_{(J)}p={\sf c}_{(J)}p.$
\item If $\sigma\upharpoonright \alpha\sim J=\tau\upharpoonright \alpha\sim J$, then
${\sf s}_{\sigma}{\sf c}_{(J)}p={\sf s}_{\tau}{\sf c}_{(J)}p$ and ${\sf s}_{\sigma}{\sf q}_{(J)}p={\sf s}_{\tau}{\sf q}_{(J)}p.$
\item If $\sigma\upharpoonright \sigma^{-1}(J)$ is injective, then
${\sf c}_{(J)}{\sf s}_{\sigma}p={\sf s}_{\sigma}{\sf c}_{\sigma^{-1}(J)}p$
and ${\sf q}_{(J)}{\sf s}_{\sigma}p={\sf s}_{\sigma}{\sf q}_{\sigma^{-1}(J)}p.$
\end{enumarab}
\end{definition}

\begin{definition} Let $\alpha$ and $G$ be as in the prevoius definition. 
By a $G$ polyadic equality algebra, a $GPHAE_{\alpha}$ for short, 
we understand an algebra  of the form  
$$(A,\lor,\land,\rightarrow, 0, {\sf s}_{\tau}, {\sf c}_{(J)}, {\sf q}_{(J)},  {\sf d}_{ij})_{\tau\in G, J\subseteq_{\omega} \alpha, i,j\in \alpha}$$
where $(A,\lor,\land,\rightarrow, 0, {\sf s}_{\tau}, {\sf c}_{(J)}, {\sf q}_{(J)})_{\tau\in G\subseteq {}^{\alpha}\alpha, J\subseteq_{\omega} \alpha}$
is a $GPHA_{\alpha}$ and ${\sf d}_{ij}\in A$ for each $i,j\in \alpha,$ such that
the following identities hold for all $k,l\in \alpha$
and all $\tau\in G:$
\begin{enumarab}
\item ${\sf d}_{kk}=1$

\item ${\sf s}_{\tau}{\sf d}_{kl}={\sf d}_{\tau(k), \tau(l)}.$

\item $x\cdot {\sf d}_{kl}\leq {\sf s}_{[k|l]}x$

\end{enumarab}
\end{definition}

Here $[k|l]$ is the replacement that sends $k$ to $l$ and otherwise is the identity.
In our definition of algebras, we depart from \cite{HMT2} by defining polyadic algebras on sets rather than on ordinals. 
In this manner, we follow the tradition of Halmos.
We refer to $\alpha$ as the dimension of $\A$ and we write $\alpha=dim\A$.
Borrowing terminology from cylindric algebras, we refer to ${\sf c}_{(\{i\})}$ by ${\sf c}_i$ and ${\sf q}_{(\{i\})}$ by ${\sf q}_i.$
However, we will have occasion to impose a well order on dimensions thereby dealing with ordinals.

\begin{remark}

When $G$ consists of all finite transformations, 
then any algebra with a Boolean reduct satisfying the above identities relating cylindrifications, diagonal elements and substitutions,
will be a quasipolyadic equality algebra of infinite dimension.
\end{remark}
Besides dealing with the two extremes when $G$ consists only of finite transformations, supplied with an additional condition,
and when $G$ is $^{\alpha}\alpha$, 
we also consider cases when $G$ is a possibly  proper subsemigroup of $^{\alpha}\alpha$ (under the operation of composition).
We need some preparations to define such semigroups.

\begin{athm}{Notation.} For  a set $X$, $|X|$ stands for the cardinality of $X$.
For functions $f$ and $g$ and a set $H$, $f[H|g]$ is the function
that agrees with $g$ on $H$, and is otherwise equal to $f$. Recall that $Rgf$ denotes 
the range of $f$. For a transformation $\tau$ on $\alpha$, 
the support of $\tau$, or $sup(\tau)$ for short, is the set:
$$sup(\tau)=\{i\in \alpha: \tau(i)\neq i\}.$$
Let $i,j\in \omega$, then $\tau[i|j]$ is the transformation on $\alpha$ defined as follows:
$$\tau[i|j](x)=\tau(x)\text { if } x\neq i \text { and }\tau[i|j](i)=j.$$ 
Recall that the map $[i|j]$ is the transformation that sends $i$ to $j$ and is the equal to the identity elsewhere. On the other hand, 
the map denoted by $[i,j]$ is the transpostion that interchanges $i$ and $j$.

For a function $f$, $f^n$ denotes the composition $f\circ f\ldots \circ f$
$n$ times.
\end{athm}
We extend the known definition of (strongly) rich semigroups \cite{S}, \cite{AUamal},  allowing possibly uncountable sets and semigroups. 
This will be needed when $G={}^\alpha\alpha$, cf. lemma \ref{cylindrify}. However, throughout when we 
mention  rich semigroups, then we will be tacitly assuming that both the dimension of the algebra involved 
and the semigroup are countable, {\it unless} otherwise explicity mentioned.

\begin{definition}\label{rich}Let $\alpha$ be any set. Let $T\subseteq \langle {}^{\alpha}\alpha, \circ \rangle$ be a semigroup.
We say that $T$ is {\it rich } if $T$ satisfies the following conditions:
\begin{enumerate}
\item $(\forall i,j\in \alpha)(\forall \tau\in T) \tau[i|j]\in T.$
\item There exist $\sigma,\pi\in T$ such that
$(\pi\circ \sigma=Id,\  Rg\sigma\neq \alpha), $ satisfying 
$$ (\forall \tau\in T)(\sigma\circ \tau\circ \pi)[(\alpha\sim Rg\sigma)|Id]\in T.$$
\end{enumerate}
\end{definition}
\begin{definition}\label{stronglyrich}
Let $T\subseteq \langle {}^{\alpha}\alpha, \circ\rangle$ be a rich semigroup. 
Let $\sigma$ and $\pi$ be as in the previous definition. 
If $\sigma$ and $\pi$ satisfy:

\begin{enumerate}
\item $(\forall n\in \omega) |supp(\sigma^n\circ \pi^n)|<\alpha, $  
\item $(\forall n\in \omega)[supp(\sigma^n\circ \pi^n)\subseteq 
\alpha\smallsetminus Rng(\sigma^n)];$
\end{enumerate}
then we say that $T$ is  {\it a strongly rich} semigroup. 
\end{definition}

\begin{example}
Examples of  rich semigroups of $\omega$ are $(^{\omega}\omega, \circ)$ and its  semigroup generated by
$\{[i|j], [i,j], i, j\in \omega, suc, pred\}$. Here $suc$ abbreviates the successor function 
on $\omega$ and $pred$ acts as its right inverse, the predecessor
function, defined by 
$pred(0)=0$ and for other $n\in \omega$, $pred(n)=n-1$. In fact, both semigroups are strongly rich, 
in the second case $suc$ plays the role of $\sigma$ while $pred$ plays the role of 
$\pi$. 
\end{example}
Rich semigroups were introduced in \cite{S} (to prove a representability result) 
and those that are strongly rich were intoduced in \cite{AUamal} (to prove an amalgamation result).

Next, we collect some properties of $G$ algebras that are more handy to use in our subsequent work.
In what follows, we will be writing $GPHA$ ($GPHAE$)  for all algebras considered.
\begin{theorem}\label{axioms} Let $\alpha$ be an infinite set and $\A\in GPHA_{\alpha}$. 
Then $\A$ satisfies the following identities for $\tau,\sigma\in G$
and all $i,j,k\in \alpha$.
\begin{enumerate}

\item $x\leq {\sf c}_ix={\sf c}_i{\sf c}_ix,\ {\sf c}_i(x\lor y)={\sf c}_ix\lor {\sf c}_iy,\ {\sf c}_i{\sf c}_jx={\sf c}_j{\sf c}_ix$.

That is  ${\sf c}_i$ is an additive operator (a modality)  and ${\sf c}_i,{\sf c}_j$ commute.

\item ${\sf s}_{\tau}$ is a Heyting algebra  endomorphism.

\item  ${\sf s}_{\tau}{\sf s}_{\sigma}x={\sf s}_{\tau\circ \sigma}x$
and ${\sf s}_{Id}x=x$.

\item ${\sf s}_{\tau}{\sf c}_ix={\sf s}_{\tau[i|j]}{\sf c}_ix$.

Recall that $\tau[i|j]$ is the transformation that agrees with $\tau$ on 
$\alpha\smallsetminus\{i\}$ and $\tau[i|j](i)=j$.

\item ${\sf s}_{\tau}{\sf c}_ix={\sf c}_j{\sf s}_{\tau}x$ if $\tau^{-1}(j)=\{i\}$, 
${\sf s}_{\tau}{\sf q}_ix={\sf q}_j{\sf s}_{\tau}x$ 
if $\tau^{-1}(j)=\{i\}$.  

\item  ${\sf c}_i{\sf s}_{[i|j]}x={\sf s}_{[i|j]}x$,\ \ ${\sf q}_i{\sf s}_{[i|j]}x={\sf s}_{[i|j]}x$

\item  ${\sf s}_{[i|j]}{\sf c}_ix={\sf c}_ix$, \ \ ${\sf s}_{[i|j]}{\sf q}_ix={\sf q}_ix$.

\item ${\sf s}_{[i|j]}{\sf c}_kx={\sf c}_k{\sf s}_{[i|j]}x$,\ \ ${\sf s}_{[i|j]}{\sf q}_kx={\sf q}_k{\sf s}_{[i|j]}x$
whenever $k\notin \{i,j\}$.

\item  ${\sf c}_i{\sf s}_{[j|i]}x={\sf c}_j{\sf s}_{[i|j]}x$,\ \  ${\sf q}_i{\sf s}_{[j|i]}x={\sf q}_j{\sf s}_{[i|j]}x$.
\end{enumerate}
\end{theorem}
\begin{demo}{Proof} The proof is tedious but fairly straighforward.
\end{demo}
Obviously the previous equations hold in $GPHAE_{\alpha}$.
Following cylindric algebra tradition and terminology, we will be often writing ${\sf s}_j^i$ for ${\sf s}_{[i|j]}$.

\begin{remark} For $GPHA_{\alpha}$ when $G$ is rich or $G$ consists only of finite transformation it is enough to restrict our attenstion to replacements. Other substitutions are 
definable from those.
\end{remark}

\subsection{Neat reducts and dilations}

Now we recall the important notion of neat reducts, a central concept in cylindric algebra theory, strongly related to representation theorems.
This concept also occurs in polyadic algebras, but unfortunately under a different name, that of compressions.

Forming dilations of an algebra, is basically  an algebraic reflection of a Henkin 
construction; in fact, the dilation of an algebra is another algebra that has an infinite number of new dimensions (constants) 
that potentially eliminate cylindrifications (quantifiers). Forming neat reducts has to do with 
restricting or compressing dimensions (number of variables) rather than increasing them.
(Here the duality has a precise categorical sense which will be formulated in the part 3 of this paper as an adjoint situation).

\begin{definition} 
\begin{enumarab}
\item  Let $ \alpha\subseteq \beta$ be infinite sets. Let $G_{\beta}$ be a semigroup of transformations on $\beta$, 
and let $G_{\alpha}$ be a semigroup of transformations on $\alpha$ such that for all $\tau\in G_{\alpha}$, one has $\bar{\tau}=\tau\cup Id\in G_{\beta}$.
Let $\B=(B, \lor, \land, \to, 0,  {\sf c}_i, {\sf s}_{\tau})_{i\in \beta, \tau\in G_{\beta}}$ be a $G_{\beta}$ algebra.
\begin{enumroman}

\item  We denote by $\Rd_{\alpha}\B$ the $G_{\alpha}$ algebra obtained by dicarding operations in $\beta\sim \alpha$. That is 
$\Rd_{\alpha}\B=(B, \lor,  \land, \to, 0, {\sf c}_i, {\sf s}_{\bar{\tau}})_{i\in \alpha, \tau\in G_{\alpha}}$. Here ${\sf s}_{\bar{\tau}}$ is evaluated 
in $\B$.
\item For $x\in B$,  then $\Delta x,$ the dimension set of $x$, 
is defined by $\Delta x=\{i\in \beta: {\sf c}_ix\neq x\}.$
Let $A=\{x\in B: \Delta x\subseteq \alpha\}$. If $A$ is
a subuniverse of $\Rd_{\alpha}\B$, then $\A$ (the algebra with universe $A$) is a subreduct of $\B$, it is called the {\it neat $\alpha$ 
reduct} of $\B$ and is denoted by $\Nr_{\alpha}\B$.
\end{enumroman}
\item If $\A\subseteq \Nr_{\alpha}\B$, then $\B$ is called a {\it dilation} of $\A$, and we say that $\A$ {\it neatly embeds} in $\B$.
if $A$ generates $\B$ (using all operations of $\B$), then $\B$ is called a {\it minimal dilation} of $\A$.
\end{enumarab}
\end{definition}
The above definition applies equally well to $GPHAE_{\alpha}$. 

\begin{remark}
In certain contexts minimal dilations may not be unique (up to isomorphism), but what we show next is that in all the cases 
we study, they are unique, so for a given algebra $\A$, we may safely say {\it the} minimal dilation of $\A$.
\end{remark}

For an algebra $\A$, and $X\subseteq \A$, $\Sg^{\A}X$ or simply $\Sg X$, when $\A$ is clear from context,
denotes the subalgebra of $\A$ generated by $X.$
The next theorems apply equally well to $GPHAE_{\alpha}$ with easy modifications which we state as we go along. 

\begin{lemma}\label{dl} 
\begin{enumarab}
\item Let $\alpha\subseteq \beta$ be countably infinite sets. If $G$ is a strongly rich semigroup on $\alpha$ 
and $\A\in GPHA_{\alpha}$, then there exists a strongly rich semigroup $T$ on $\beta$ and $\B\in TPHA_{\beta},$ such that
$\A\subseteq \Nr_{\alpha}\B$ and for all $X\subseteq A,$ one has  $\Sg^{\A}X=\Nr_{\alpha}\Sg^{\B}X$. 

\item Let $G_{I}$ be the semigroup of finite transformations on $I$.
Let $\A\in G_{\alpha}PHA_{\alpha}$ be such that $\alpha\sim \Delta x$ is infinite for every $x\in A$. 
Then for any set $\beta$, such that $\alpha\subseteq \beta$, there exists $\B\in G_{\beta}PHA_{\beta},$ 
such that $\A\subseteq \Nr_{\alpha}\B$ and for all $X\subseteq 
A$, one has $\Sg^{\A}X=\Nr_{\alpha}\Sg^{\B}X.$

\item Let $G_I$ be the semigroup of all transformations on $I$. Let $\A\in G_{\alpha}PHA_{\alpha}$. 
Then for any set $\beta$ such that $\alpha\subseteq \beta$, 
there exists $\B\in G_{\beta}PHA_{\beta},$ 
such that $\A\subseteq \Nr_{\alpha}\B$ and for all $X\subseteq 
A,$ one has $\Sg^{\A}X=\Nr_{\alpha}\Sg^{\B}X.$

\end{enumarab}

\end{lemma}

\begin{demo}{Proof}
\begin{enumarab}  
\item cf. \cite{AUamal}. We assume that $\alpha$ is an ordinal; in fact without loss of generality we can assume that it 
is the least infinite ordinal $\omega.$
We also assume a particular strongly rich semigroup, that namely that generated by finite transformations together with 
$suc$, $pred$. The general case is the same \cite{AUamal} Remark 2.8 p. 327.
We follow \cite {AUamal} p. 323-336. 
For $n\leq \omega$, let  $\alpha_n=\omega+n$ 
and $M_n=\alpha_n\sim \omega$.
Note that when $n\in \omega$, then $M_n=\{\omega,\ldots,\omega+n-1\}$.
Let $\tau\in G$. Then $\tau_n=\tau\cup Id_{M_n}$. $T_n$ denotes the 
subsemigroup of $\langle {}^{\alpha_n}\alpha_n,\circ \rangle$ generated by
$\{\tau_n:\tau\in G\} \cup \cup_{i,j\in \alpha_n}\{[i|j],[i,j]\}$.
For $n\in \omega$, we let $\rho_n:\alpha_n\to \omega$ 
be the bijection defined by 
$\rho_n\upharpoonright \omega=suc^n$ and $\rho_n(\omega+i)=i$ for all $i<n$.
Let $n\in \omega$. For $v\in T_n,$ let $v'=\rho_n\circ v\circ \rho_n^{-1}$.
Then $v'\in G$.
For $\tau\in T_{\omega}$, let 
$D_{\tau}=\{m\in M_{\omega}:\tau^{-1}(m)=\{m\}=\{\tau(m)\}\}$.
Then $|M_{\omega}\sim D_{\tau}|<\omega.$ 
Let $\A$ be a given  countable $G$ algebra.  
Let $\A_n$ be the algebra defined as follows:
$\A_n=\langle A,\lor, \land, \to, 0, {\sf c}_i^{\A_n},{\sf s}_v^{\A_n}\rangle_{i\in \alpha_n,v\in T_n}$
where for each $i\in \alpha_n$ and $v\in T_n$, 
${\sf c}_i^{\A_n}:= {\sf c}_{\rho_n(i)}^{\A} \text { and }{\sf s}_v^{\A_n}:= {\sf s}_{v'}^{\A}.$
Let $\Rd_{\omega}\A_n$ be the following reduct of 
$\A_n$ obtained by restricting the type of $\A_n$ to the first 
$\omega$ dimensions:
$\Rd_{\omega}\A_n=\langle A_n,\lor,\land, \to,0, {\sf c}_i^{\A_n}, {\sf s}_{\tau_n}^{\A_n}\rangle_{i\in \omega,\tau\in G}.$
For $x\in A$, let $e_n(x)={\sf s}_{suc^n}^{\A}(x)$. Then $e_n:A\to A_n$ and  $e_n$ is an isomorphism 
from $\A$ into $\Rd_{\omega}\A_n$
such that $e_n(\Sg^{\A}Y)=\Nr_{\omega}(\Sg^{\A_n}e_n(Y))$ for all
$Y\subseteq A$, cf. \cite{AUamal} claim 2.7. While $\sigma$ and condition (2) in the definition of \ref{rich} are needed to implement the neat embedding, the left inverse $\pi$ of $\sigma$ is needed to
show that forming neat reducts commute with froming subalgebras; in particular $\A$ is the full $\omega$ neat reduct
of $\A_n$. To extend the neat embedding part to infinite dimensions, we use a fairly straightforward construction
involving an ultraproduct of exapansions of the algebras $\A_n$, on any cofinite ultrafilter on $\omega$.
For the sake of brevity, let $\alpha=\alpha_{\omega}=\omega+\omega$.
Let $T_{\omega}$ is the semigroup generated by the set
$\{\tau_{\omega}: \tau\in G\}\cup_{i,j\in \alpha}\{[i|j],[i,j]\}.$
For $\sigma\in T_{\omega}$, and $n\in \omega$, 
let $[\sigma]_n=\sigma\upharpoonright \omega+n$.
For each $n\in \omega,$ let 
$\A_n^+=\langle A,\lor,\land, \to, 0, {\sf c}_i^{\A_n^+}, {\sf s}_{\sigma}^{\A_n^+}\rangle_{i\in \alpha, \sigma\in T_{\omega}}$
be an expansion of $\A_n$ 
such that there Heyting  reducts coincide  
and for each 
$\sigma\in T_{\omega}$ and $i\in \alpha,$ 
${\sf s}_{\sigma}^{\A_n^+}:={\sf s}_{[\sigma]_n}^{\A_n} 
\text { iff } [\sigma]_n\in T_n,$ 
and 
${\sf c}_i^{\A_n^+}:={\sf c}_i^{\A_n}\text { iff }i<\omega+n.$
Let $F$ to be any non-principal ultrafilter on $\omega$. Now 
forming the ultraproduct of the $\A_n^+$'s relative to $F$, let
$\A^+=\prod_{n\in \omega}\A_n^+/F.$
For $x\in A$, let 
$e(x)=\langle e_n(x):n\in \omega\rangle/F.$
Let $\Rd_{\omega}A^+=\langle A^+, \lor, \land, \to, 0, {\sf c}_i^{\A^+}, {\sf s}_{\tau_{\omega}}^{\A^+}
\rangle_{i<\omega,\tau\in T}.$
Then $e$ is an isomorphism from $\A$ into $\Rd_{\omega}\A^+$
such that  $e(\Sg^{\A}Y)=\Nr_{\omega}\Sg^{\A^+}e(Y)$ for all $Y\subseteq A.$

We have shown that $\A$ neatly embeds in algebras in finite extra dimensions and in $\omega$ extra dimension. 
An iteration of this embedding yields the required result. 

In the presence of diagonals one has to check that homomorphisms defined preserve diagonal elements. 
But this is completely straightforward using properties of substitutions when applied to diagonal elements.

\item  Let $\alpha\subseteq \beta$. We assume, loss of generality, that $\alpha$ and $\beta$ are ordinals with $\alpha<\beta$.
The proof is a direct adaptation of the proof of Theorem 2.6.49(i)
in \cite{HMT1}. First we show that there exists $\B\in G_{\alpha+1}PHA_{\alpha+1}$ 
such that $\A$ embeds into $\Nr_{\alpha}\B,$ then we proceed inductively. 
Let $$R = Id\upharpoonright (\alpha\times A)
 \cup \{  ((k,x), (\lambda, y)) : k, \lambda <
\alpha, x, y \in A, \lambda \notin \Delta x, y = {\mathsf s}_{[k|\lambda]} x \}.$$
It is easy to see  that $R$ is an equivalence relation on $\alpha
\times A$.
Define the following operations on $(\alpha\times A)/R$ with $\mu, i, k\in \alpha$ and $x,y\in A$ :
\begin{equation*}\label{l5}
\begin{split}
(\mu, x)/R \lor (\mu, y)/R = (\mu, x \lor y)/R, 
\end{split}
\end{equation*}
\begin{equation*}\label{l6}
\begin{split}
(\mu, x)/R\land  (\mu, y)/R = (\mu, x\land  y)/R, 
\end{split}
\end{equation*}
 \begin{equation*}\label{l7}
\begin{split}
(\mu, x)/R\to  (\mu, y)/R = (\mu, x\to  y)/R, 
\end{split}
\end{equation*}
\begin{equation*}\label{l10}
\begin{split}
{\mathsf c}_i ((\mu, x)/R)  = (\mu, {\mathsf c}_i x )/R, \quad
\mu \in \alpha \smallsetminus
\{i\},
\end{split}
\end{equation*}
\begin{equation*}\label{l11}
\begin{split}
{\mathsf s}_{[j|i]} ((\mu, x)/R)  = (\mu, {\mathsf s}_{[j|i]} x )/R, \quad \mu \in \alpha
\smallsetminus \{i, j\}.
\end{split}
\end{equation*}
It can be checked that these operations are well defined.
Let $$\C=((\alpha\times A)/R, \lor, \land, \to, 0, {\sf c_i}, {\sf s}_{i|j]})_{i,j\in \alpha},$$
and
let $$h=\{(x, (\mu,x)/R): x\in A, \mu\in \alpha\sim \Delta x\}.$$
Then $h$ is an isomorphism from $\A$ into $\C$. 
Now to show that $\A$ neatly embeds into $\alpha+1$ extra dimensions, we define the operations ${\sf c}_{\alpha}, {\sf s}_{[i|\alpha]}$
and ${\sf s}_{[\alpha|i]}$ on $\C$ as follows:
$${\mathsf c}_\alpha = \{ ((\mu, x)/R, (\mu, {\mathsf c}_\mu x)/R) :
\mu \in \alpha, x \in B \},$$
$${\mathsf s}_{[i|\alpha]} = \{ ((\mu, x)/R, (\mu, {\mathsf s}_{[i|\mu]}
x)/R) : \mu \in \alpha \smallsetminus \{i\}, x \in B \},$$
$${\mathsf s}_{[\alpha|i]} = \{ ((\mu, x)/R, (\mu, {\mathsf s}_{[\mu|i]}
x)/R) : \mu \in \alpha \smallsetminus \{i\}, x \in B \}.$$
Let $$\B=((\alpha\times A)/R, \lor,\land, \to, {\sf c}_i, {\sf s}_{[i|j]})_{i,j\leq \alpha}.$$
Then $$\B\in G_{\alpha+1}PA_{\alpha+1}\text{ and }h(\A)\subseteq \Nr_{\alpha}\B.$$
It is not hard to check that the defined operations are as desired. We have our result when $G$ consists only of replacements.
But since $\alpha\sim \Delta x$ is infinite one can show that substitutions corresponding to all finite transformations are term definable.
For a finite transformation $\tau\in {}^{\alpha}\alpha$ we write $[u_0|v_0, u_1|v_1,\ldots,
u_{k-1}|v_{k-1}]$ if $sup\tau=\{u_0,\ldots ,u_{k-1}\}$, $u_0<u_1
\ldots <u_{k-1}$ and $\tau(u_i)=v_i$ for $i<k$.
Let $\A\in GPHA_{\alpha}$ be such that $\alpha\sim \Delta x$ is
infinite for every $x\in A$. If $\tau=[u_0|v_0, u_1|v_1,\ldots,
u_{k-1}|v_{k-1}]$ is a finite transformation, if $x\in A$ and if
$\pi_0,\ldots ,\pi_{k-1}$ are in this order the first $k$ ordinals
in $\alpha\sim (\Delta x\cup Rg(u)\cup Rg(v))$, then
$${\mathsf s}_{\tau}x={\mathsf s}_{v_0}^{\pi_0}\ldots
{\mathsf s}_{v_{k-1}}^{\pi_{k-1}}{\mathsf s}_{\pi_0}^{u_0}\ldots
{\mathsf s}_{\pi_{k-1}}^{u_{k-1}}x.$$
The ${\sf s}_{\tau}$'s so defined satisfy the polyadic axioms, cf \cite{HMT1} Theorem 1.11.11.
Then one proceeds by a simple induction to show that for all $n\in \omega$ there exists $\B\in G_{\alpha+n}PHA_{\alpha+n}$ 
such that $\A\subseteq \Nr_{\alpha}\B.$ For the transfinite, one uses ultraproducts \cite{HMT1} theorem 2.6.34. 
For the second part, let $\A\subseteq \Nr_{\alpha}\B$ and $A$ generates $\B$ then $\B$ consists of all elements ${\sf s}_{\sigma}^{\B}x$ such that 
$x\in A$ and $\sigma$ is a finite transformation on $\beta$ such that
$\sigma\upharpoonright \alpha$ is one to one \cite{HMT1} lemma 2.6.66. 
Now suppose $x\in \Nr_{\alpha}\Sg^{\B}X$ and $\Delta x\subseteq
\alpha$, then there exist $y\in \Sg^{\A}X$ and a finite transformation $\sigma$
of $\beta$ such that $\sigma\upharpoonright \alpha$ is one to one
and $x={\sf s}_{\sigma}^{\B}y.$  
Let $\tau$ be a finite
transformation of $\beta$ such that $\tau\upharpoonright  \alpha=Id
\text { and } (\tau\circ \sigma) \alpha\subseteq \alpha.$ Then
$x={\sf s}_{\tau}^{\B}x={\sf s}_{\tau}^{\B}{\sf s}_{\sigma}y=
{\sf s}_{\tau\circ \sigma}^{\B}y={\sf s}_{\tau\circ
\sigma\upharpoonright \alpha}^{\A}y.$
In the presence of diagonal elements, one defines them in the bigger algebra (the dilation) precisely  as in \cite{HMT1}, theorem 2.6.49(i).

\item Here we extensively use the techniques in \cite{DM}, but we have to watch out, for we only have finite cylindrifications.
Let $(\A, \alpha,S)$ be a transformation system. 
That is to say, $\A$ is a Heyting algebra and $S:{}^\alpha\alpha\to End(\A)$ is a homomorphism. For any set $X$, let $F(^{\alpha}X,\A)$ 
be the set of all functions from $^{\alpha}X$ to $\A$ endowed with Heyting operations defined pointwise and for 
$\tau\in {}^\alpha\alpha$ and $f\in F(^{\alpha}X, \A)$, ${\sf s}_{\tau}f(x)=f(x\circ \tau)$. 
This turns $F(^{\alpha}X,\A)$ to a transformation system as well. 
The map $H:\A\to F(^{\alpha}\alpha, \A)$ defined by $H(p)(x)={\sf s}_xp$ is
easily checked to be an isomorphism. Assume that $\beta\supseteq \alpha$. Then $K:F(^{\alpha}\alpha, \A)\to F(^{\beta}\alpha, \A)$ 
defined by $K(f)x=f(x\upharpoonright \alpha)$ is an isomorphism. These facts are straighforward to establish, cf. theorem 3.1, 3.2 
in \cite{DM}. 
$F(^{\beta}\alpha, \A)$ is called a minimal dilation of $F(^{\alpha}\alpha, \A)$. Elements of the big algebra, or the cylindrifier free 
dilation, are of form ${\sf s}_{\sigma}p$,
$p\in F(^{\beta}\alpha, \A)$ where $\sigma$ is one to one on $\alpha$, cf. \cite{DM} theorem 4.3-4.4.
We say that $J\subseteq I$ supports an element $p\in A,$ if whenever $\sigma_1$ and  $\sigma_2$ are 
transformations that agree on $J,$ then  ${\sf s}_{\sigma_1}p={\sf s}_{\sigma_2}p$.
$\Nr_JA$, consisting of the elements that $J$ supports, is just the  neat $J$ reduct of $\A$; 
with the operations defined the obvious way as indicated above. 
If $\A$ is an $\B$ valued $I$ transformaton system with domain $X$, 
then the $J$ compression of $\A$ is isomorphic to a $\B$ valued $J$ transformation system
via $H: \Nr_J\A\to F(^JX, \A)$ by setting for $f\in\Nr_J\A$ and $x\in {}^JX$, $H(f)x=f(y)$ where $y\in X^I$ and $y\upharpoonright J=x$, 
cf. \cite{DM} theorem 3.10.
Now let $\alpha\subseteq \beta.$ If $|\alpha|=|\beta|$ then the the required algebra is defined as follows. 
Let $\mu$ be a bijection from $\beta$ onto $\alpha$. For $\tau\in {}^{\beta}\beta,$ let ${\sf s}_{\tau}={\sf s}_{\mu\tau\mu^{-1}}$ 
and for each $i\in \beta,$ let 
${\sf c}_i={\sf c}_{\mu(i)}$. Then this defined $\B\in GPHA_{\beta}$ in which $\A$ neatly embeds via ${\sf s}_{\mu\upharpoonright\alpha},$
cf. \cite{DM} p.168.  Now assume that $|\alpha|<|\beta|$.
Let $\A$ be a  given polyadic algebra of dimension $\alpha$; discard its cylindrifications and then take its minimal dilation $\B$, 
which exists by the above.
We need to define cylindrifications on the big algebra, so that they agree with their values in $\A$ and to have $\A\cong \Nr_{\alpha}\B$. We let (*):
$${\sf c}_k{\sf s}_{\sigma}^{\B}p={\sf s}_{\rho^{-1}}^{\B} {\sf c}_{\rho(\{k\}\cap \sigma \alpha)}{\sf s}_{(\rho\sigma\upharpoonright \alpha)}^{\A}p.$$
Here $\rho$ is a any permutation such that $\rho\circ \sigma(\alpha)\subseteq \sigma(\alpha.)$
Then we claim that the definition is sound, that is, it is independent of $\rho, \sigma, p$. 
Towards this end, let $q={\sf s}_{\sigma}^{\B}p={\sf s}_{\sigma_1}^{\B}p_1$ and 
$(\rho_1\circ \sigma_1)(\alpha)\subseteq \alpha.$
We need to show that (**)
$${\sf s}_{\rho^{-1}}^{\B}{\sf c}_{[\rho(\{k\}\cap \sigma(\alpha)]}^{\A}{\sf s}_{(\rho\circ \sigma\upharpoonright \alpha)}^{\A}p=
{\sf s}_{\rho_1{^{-1}}}^{\B}{\sf c}_{[\rho_1(\{k\}\cap \sigma(\alpha)]}^{\A}{\sf s}_{(\rho_1\circ \sigma\upharpoonright \alpha)}^{\A}p.$$
Let $\mu$ be a permutation of $\beta$ such that
$\mu(\sigma(\alpha)\cup \sigma_1(\alpha))\subseteq \alpha$.
Now applying ${\sf s}_{\mu}$ to the left hand side of (**), we get that 
$${\sf s}_{\mu}^{\B}{\sf s}_{\rho^{-1}}^{\B}{\sf c}_{[\rho(\{k\})\cap \sigma(\alpha)]}^{\A}{\sf s}_{(\rho\circ \sigma|\alpha)}^{\A}p
={\sf s}_{\mu\circ \rho^{-1}}^{\B}{\sf c}_{[\rho(\{k\})\cap \sigma(\alpha)]}^{\A}{\sf s}_{(\rho\circ \sigma|\alpha)}^{\A}p.$$
The latter is equal to
${\sf c}_{(\mu(\{k\})\cap \sigma(\alpha))}{\sf s}_{\sigma}^{\B}q.$
Now since $\mu(\sigma(\alpha)\cap \sigma_1(\alpha))\subseteq \alpha$, we have
${\sf s}_{\mu}^{\B}p={\sf s}_{(\mu\circ \sigma\upharpoonright \alpha)}^{\A}p={\sf s}_{(\mu\circ \sigma_1)\upharpoonright \alpha)}^{\A}p_1\in A$.
It thus follows that 
$${\sf s}_{\rho^{-1}}^{\B}{\sf c}_{[\rho(\{k\})\cap \sigma(\alpha)]}^{\A}{\sf s}_{(\rho\circ \sigma\upharpoonright \alpha)}^{\A}p=
{\sf c}_{[\mu(\{k\})\cap \mu\circ \sigma(\alpha)\cap \mu\circ \sigma_1(\alpha))}{\sf s}_{\sigma}^{\B}q.$$ 
By exactly the same method, it can be shown that 
$${\sf s}_{\rho_1{^{-1}}}^{\B}{\sf c}_{[\rho_1(\{k\})\cap \sigma(\alpha)]}^{\A}{\sf s}_{(\rho_1\circ \sigma\upharpoonright \alpha)}^{\A}p
={\sf c}_{[\mu(\{k\})\cap \mu\circ \sigma(\alpha)\cap \mu\circ \sigma_1(\alpha))}{\sf s}_{\sigma}^{\B}q.$$ 
By this we have proved (**).

Furthermore, it defines the required algebra $\B$. Let us check this.
Since our definition is slightly different than that in \cite{DM}, by restricting cylindrifications to be olny finite, 
we need to check the polyadic axioms which is tedious but routine. The idea is that every axiom can be pulled back to 
its corresponding axiom holding in the small algebra 
$\A$.
We check only the axiom $${\sf c}_k(q_1\land {\sf c}_kq_2)={\sf c}_kq_1\land {\sf c}_kq_2.$$
We follow closely \cite{DM} p. 166. 
Assume that $q_1={\sf s}_{\sigma}^{\B}p_1$ and $q_2={\sf s}_{\sigma}^{\B}p_2$. 
Let $\rho$ be a permutation of $I$ such that $\rho(\sigma_1I\cup \sigma_2I)\subseteq I$ and let 
$$p={\sf s}_{\rho}^{\B}[q_1\land {\sf c}_kq_2].$$
Then $$p={\sf s}_{\rho}^{\B}q_1\land {\sf s}_{\rho}^{\B}{\sf c}_kq_2
={\sf s}_{\rho}^{\B}{\sf s}_{\sigma_1}^{\B}p_1\land {\sf s}_{\rho}^{\B}{\sf c}_k {\sf s}_{\sigma_2}^{\B}p_2.$$
Now we calculate ${\sf c}_k{\sf s}_{\sigma_2}^{\B}p_2.$
We have by (*)
$${\sf c}_k{\sf s}_{\sigma_2}^{\B}p_2= {\sf s}^{\B}_{\sigma_2^{-1}}{\sf c}_{\rho(\{k\}\cap \sigma_2I)} {\sf s}^{\A}_{(\rho\sigma_2\upharpoonright I)}p_2.$$
Hence $$p={\sf s}_{\rho}^{\B}{\sf s}_{\sigma_1}^{\B}p_1\land {\sf s}_{\rho}^{\B}{\sf s}^{\B}_{\sigma^{-1}}{\sf c}_{\rho(\{k\}\cap \sigma_2I)} 
{\sf s}^{\A}_{(\rho\sigma_2\upharpoonright I)}p_2.$$
\begin{equation*}
\begin{split}
&={\sf s}^{\A}_{\rho\sigma_1\upharpoonright I}p_1\land {\sf s}_{\rho}^{\B}{\sf s}^{\A}_{\sigma^{-1}}{\sf c}_{\rho(\{k\}\cap \sigma_2I)} 
{\sf s}^{\A}_{(\rho\sigma_2\upharpoonright I)}p_2,\\
&={\sf s}^{\A}_{\rho\sigma_1\upharpoonright I}p_1\land {\sf s}_{\rho\sigma^{-1}}^{\A}
{\sf c}_{\rho(\{k\}\cap \sigma_2I)} {\sf s}^{\A}_{(\rho\sigma_2\upharpoonright I)}p_2,\\
&={\sf s}^{\A}_{\rho\sigma_1\upharpoonright I}p_1\land {\sf c}_{\rho(\{k\}\cap \sigma_2I)} {\sf s}^{\A}_{(\rho\sigma_2\upharpoonright I)}p_2.\\
\end{split}
\end{equation*} 
Now $${\sf c}_k{\sf s}_{\rho^{-1}}^{\B}p={\sf c}_k{\sf s}_{\rho^{-1}}^{\B}{\sf s}_{\rho}^{\B}(q_1\land {\sf c}_k q_2)={\sf c}_k(q_1\land {\sf c}_kq_2)$$
We next calculate ${\sf c}_k{\sf s}_{\rho^{-1}}p$.
Let $\mu$ be a permutation of $I$ such that $\mu\rho^{-1}I\subseteq I$. Let $j=\mu(\{k\}\cap \rho^{-1}I)$.
Then applying (*), we have:
\begin{equation*}
\begin{split}
&{\sf c}_k{\sf s}_{\rho^{-1}}p={\sf s}^{\B}_{\mu^{-1}}{\sf c}_{j}{\sf s}_{(\mu\rho^{-1}|I)}^{\A}p,\\
&={\sf s}^{\B}_{\mu^{-1}}{\sf c}_{j}{\sf s}_{(\mu\rho^{-1}|I)}^{\A}
{\sf s}^{\A}_{\rho\sigma_1\upharpoonright I}p_1\land {\sf c}_{(\rho\{k\}\cap \sigma_2I)} {\sf s}^{\B}_{(\rho\sigma_2\upharpoonright I)}p_2,\\
 &={\sf s}^{\B}_{\mu^{-1}}{\sf c}_{j}[{\sf s}_{\mu \sigma_1\upharpoonright I}p_1\land r].\\
\end{split}
\end{equation*}

where 
$$r={\sf s}_{\mu\rho^{-1}}^{\B}{\sf c}_j {\sf s}_{\rho \sigma_2\upharpoonright I}^{\A}p_2.$$
Now ${\sf c}_kr=r$. Hence, applying the axiom in the small algebra, we get: 
$${\sf s}^{\B}_{\mu^{-1}}{\sf c}_{j}[{\sf s}_{\mu \sigma_1\upharpoonright I}^{\A}p_1]\land {\sf c}_k q_2
={\sf s}^{\B}_{\mu^{-1}}{\sf c}_{j}[{\sf s}_{\mu \sigma_1\upharpoonright I}^{\A}p_1\land r].$$
But
$${\sf c}_{\mu(\{k\}\cap \rho^{-1}I)}{\sf s}_{(\mu\sigma_1|I)}^{\A}p_1=
{\sf c}_{\mu(\{k\}\cap \sigma_1I)}{\sf s}_{(\mu\sigma_1|I)}^{\A}p_1.$$
So 
$${\sf s}^{\B}_{\mu^{-1}}{\sf c}_{k}[{\sf s}_{\mu \sigma_1\upharpoonright I}^{\A}p_1]={\sf c}_kq_1,$$ and 
we are done.
To show that neat reducts commute with forming subalgebras, we proceed as 
in the previous proof replacing finite transformation by transformation.

When we have diagonal elements, we first discard them, obtaining a $GPHA_{\alpha}$ then form the diagonal free dilation of this algebra, and 
finally define
the diagonal elements in the dilation as in \cite{HMT2}, theorem 5.4.17, p.233.
\end{enumarab}

\end{demo}

The next lemma formulated only for $GPHA_{\alpha}$ will be used in proving our main (algebraic) result.
The proof works without any modifications when we add diagonal elements. The lemma says, roughly, that if we have an $\alpha$ dimensional  
algebra $\A$, and a set $\beta$ containing $\alpha$, then we can find an extension $\B$ of $\A$ in $\beta$ dimensions, specified by a 
carefully chosen  subsemigroup of $^{\beta}\beta$, such that $\A=\Nr_{\alpha}\B$ and  for all $b\in B$, $|\Delta b\sim \alpha|<\omega$.
$\B$ is not necessarily the minimal dilation of $\A$, because the large subsemigroup chosen maybe 
smaller than the semigroup used to form 
the unique dilation. It can happen that this extension is the minimal dilation, 
but in the case we consider all transformations, the constructed algebra 
is only a proper subreduct of the dilation obtained basically by discarding 
those elements $b$ in the original dilation for which $\Delta b\sim \alpha$ 
is infinite.
\begin{lemma}\label{cylindrify}
\begin{enumarab}
\item  For a set $I,$ let $G_I$ be the semigroup of all finite transformations on $I$. Let $\alpha\subseteq \beta$ be infinite sets.
Let $\A\in G_{\alpha}PHA_{\alpha}$  and $\B\in G_{\beta}PHA_{\beta}.$ 
If $\A\subseteq \Nr_{\alpha}\B$ and $X\subseteq A$, then for any $b\in \Sg^{\B}X,$ one has $|\Delta b\sim \alpha|<\omega.$ In particular, 
the cylindrifier ${\sf c}_{(\Delta\sim\alpha)}b$, for any such $b$ is meaningful.
\item Let $\alpha<\beta$ be countable ordinals and let $G_{\alpha}$ and $G_{\beta}$ 
be strongly rich semigroups on $\alpha$ and $\beta$, respectivey.
Let  $\A\in G_{\alpha}PHA_{\alpha}$  and $\B\in G_{\beta}PHA_{\beta}.$ 
If $\A\subseteq \Nr_{\alpha}\B$ and $X\subseteq A$, then for any $b\in \Sg^{\B}X,$ we have $|\Delta b\sim \alpha|<\omega.$ 
\item For a set $I$, let $G_I$ denote the set of all transformations on $I$. 
Let $\alpha\subseteq \beta$ be infinite sets,   such that $|\alpha|<|\beta|$. Let $\A\in G_{\alpha}PHA_{\alpha}$. 
Then there exist a semigroup $S$ of $G_{\beta}$ and $\B\in SPHA_{\beta},$ 
such that $\A=\Nr_{\alpha}\B$, $S$ contains elements $\pi$, $\sigma$ as in definition \ref{stronglyrich}, 
and for all $X\subseteq A$, one has $\Sg^{\A}X=\Nr_{\alpha}\Sg^{\B}X$. Furthermore,
for all $b\in B$, $|\Delta b\sim \alpha|<\omega.$ 

In this case we say that $\B$ is  a minimal extension of $\A$. 

\item Let $\alpha\subseteq \beta$ be infinite sets, and assume that $|\alpha|=|\beta|$. Let $S\subseteq {}^{\alpha}\alpha$ 
be a semigroup that contains all finite transformations,
and two infinitary ones $\pi$ and $\sigma$ as in the definition \ref{stronglyrich}. 
Let $\A\in SPHA_{\alpha}$. Then there exist a semigroup $T\subseteq {}^{\beta}\beta$, 
such that $\B\in TPHA_{\beta}$, $\A=\Nr_{\alpha}\B,$ and for all $X\subseteq A$, one has $\Sg^{\A}X=\Nr_{\alpha}\Sg^{\B}X$. Furthermore, 
for all $b\in B$, $|\Delta b\sim \alpha|<\omega.$
\end{enumarab}
\end{lemma}
\begin{demo}{Proof} 
\begin{enumarab}
\item This trivially holds for elements of $\A$. 
The rest follows by an easy inductive argument, since substitution can move only finitely many points.
\item This part is delicate because we have infinitary substitutions, so, in principal, it can happen that $|\Delta x\sim \alpha|<\omega$ and 
$|\Delta({\sf s}_{\tau}x)\sim \alpha|\geq \omega$, when $\tau$ moves infinitely many points. 
We show that in this particular case, this cannot happen.
Let $M=\beta\sim \alpha$. We can well assume that $\beta=\omega+\omega$ and $\alpha=\omega$.
Then since $M\cap \Delta x=\emptyset$ for all $x\in \A$, it suffices to show inductively that for any $x\in \B$ and
any (unary) operation $f$ of $\B$, the following condition holds:
$$ \text {If }|M\cap \Delta x|<\omega\text {  then }|M\cap \Delta (fx)|<\omega.$$
Of course,  we should check that the above holds for the Heyting  operations as well, 
but this is absolutely straightforward.  Assume that $f$ is a substitution. 
So let $\tau\in G_{\beta}$, such that $f={\sf s}_{\tau}$. Let 
$D_{\tau}=\{m\in M:\tau^{-1}(m)=\{m\}=\{\tau(m)\}\}.$ 
Then, it is easy to check that  $|M\sim D_{\tau}|<\omega$. 
For the sake of brevity, let $C_{\tau}$ denote the finite set $M\sim D_{\tau}.$
By $|M\cap \Delta x|<\omega$, we have 
$|(M\cap \Delta x)\cup C_{\tau}|<\omega$.
We will show that $M\cap \Delta ({\sf s}_{\tau}x)\subseteq (M\cap \Delta x)\cup C_{\tau}$ 
by which we will be done. So assume that $i\in M\sim\omega$, and that
$i\notin (M\cap \Delta x)\cup C_{\tau}$. Then 
$i\in D_{\tau}\sim \Delta x$,
so $\{\tau(i)\}=\{i\}=\tau^{-1}(i)$. Thus we get that 
${\sf c}_i{\sf s}_{\tau}x={\sf s}_{\tau}{\sf c}_ix$ by item (5) in theorem \ref{axioms}, proving that $i\notin M\cap \Delta {\sf s}_{\tau}x.$
Now assume that $f={\sf c}_j$ with $j\in \alpha$. 
If $i\in M$ and $i\notin \Delta x,$ then   
we have ${\sf c}_i{\sf c}_jx={\sf c}_j{\sf c}_ix={\sf c}_jx$ and we are done in this case, too.

We note that the condition (2) in the definition of richness suffices to implement the neat embeding, 
while strong richness is needed so that $\A$ exhausts the ful neat reduct.

\item Let $\A$ and $\beta$ be given. Choose $\pi$ and $\sigma$ in ${}^{\beta}\beta$ satisfying (3) and (4) in definition \ref{rich}. 
Let $H_{\beta}=\{\rho\in {}^{\beta}\beta: |\rho(\alpha)\cap (\beta\sim \alpha)|<\omega\}\cup \{\sigma, \pi\}$.
Let $S$ be the semigroup generated by $H_{\beta}.$ 
Let $\B'\in G_{\beta}PHA_{\beta}$ be an ordinary  dilation of $\A$ where all transformations in $^{\beta}\beta$ are used. Exists by \ref{dl}. 
Then $\A=\Nr_{\alpha}\B'$. We take a suitable reduct of $\B'$. Let $\B$ be the subalgebra of 
$\B'$ generated from $A$ be all operations except for substitutions indexed by transformations not in $S$.
Then, of course $A\subseteq \B$; in fact, $\A=\Nr_{\alpha}\B$, since for each $\tau\in {}^{\alpha}\alpha$, $\tau\cup Id\in S.$ 
We check that $\B$ is as required. It suffices to show inductively that for $b\in B$, if 
$|\Delta b\sim \alpha|<\omega$, and $\rho\in S$, then $|\rho(\Delta b)\sim \alpha|<\omega$. 
For $\rho \in H_{\beta}\sim \{\pi, \sigma\}$, this easily follows from how $\rho$ is defined, otherwise the proof is as in the previous item.

\item We can obviously write $\beta$ as a sum of ordinals $\alpha+\omega$, so that $\beta$ itself is an ordinal,  
and iterate $\sigma$ as in theorem \ref{dl} (1), by noting that the proof does not depend on the 
countability of $\A,$ but rather on that of $\beta\sim \alpha$. 
In more detail, for $n\leq \omega$, let $\alpha_n=\alpha+n$ 
and $M_n=\beta\sim \alpha_n$.
For $\tau\in S$,  let $\tau_n=\tau\cup Id_{M_n}$. 
Let $T_n$ be subsemigroup of $\langle {}^{\alpha_n}\alpha_n,\circ \rangle$ generated by
$\{\tau_n:\tau\in G\} \cup \cup_{i,j\in \alpha_n}\{[i|j],[i,j]\}$.
For $n\in \omega$, we let $\rho_n:\alpha_n\to \alpha$ 
be the bijection defined by 
$\rho_n\upharpoonright \alpha=\sigma^n$ and $\rho_n(\alpha+i)=i$ for all $i<n$.
(Here $\sigma$ is as in the fdefinition of \ref{rich}. For $n\in \omega$, for $v\in T_n$ let $v'=\rho_n\circ v\circ \rho_n^{-1}$.
Then $v'\in S$.
For $\tau\in T_{\omega}$, let  $D_{\tau}=\{m\in M_{\omega}:\tau^{-1}(m)=\{m\}=\{\tau(m)\}\}$.
Then $|M_{\omega}\sim D_{\tau}|<\omega.$ 
Let $\A$ is an  $S$ algebra.   Let $\A_n$ be the algebra defined as follows:
$\A_n=\langle A,\lor, \land, \to, 0, {\sf c}_i^{\A_n},{\sf s}_v^{\A_n}\rangle_{i\in \alpha_n,v\in T_n}$
where for each $i\in \alpha_n$ and $v\in T_n$, 
${\sf c}_i^{\A_n}:= {\sf c}_{\rho_n(i)}^{\A} \text { and }{\sf s}_v^{\A_n}:= {\sf s}_{v'}^{\A.}$
Then continue as in the proof of the above theorem \ref{dl},  by taking the ultraproduct of the $\A_n$'s 
relative to a cofinite ultrafilter, one then  gets a dilation in $\beta$ dimensions in which $\A$ neatly embeds satisfying the required.
 \end{enumarab}
\end{demo}
\begin{remark} If $\A \in GPHA_{\alpha}$ where $G_{\alpha}$ is the semigroup of all transformations 
on $\alpha$ and $\alpha\subseteq \beta$, there are two kinds of extensions of $\A$ to $\beta$ dimensions. 
The minimal dilation of $\A$ which uses all substitutions
in $G_{\beta}$, and a minimal extension of $\A$ which is can be a proper subreduct of the 
minimal dilation, using operations in a rich subsemigroup
of $G_{\beta}.$
\end{remark}  

\section{Algebraic Proofs of main theorems}

Henceforth,  when we write $GPHA_{\alpha}$ without further specification, we understand that we simultaneously  dealing with all possibilities of $G$,
and that whatever we are saying applies equaly well to all cases considered.  
We could also say $\A$ is a $G$ algebra without further notice; 
the same is to be understood. Throughout the paper dimensions will be specified by {\it infinite} sets or ordinals. 

Our work in this section is closely related to that in 
\cite{Hung}. Our main theorem is a typical representabilty result, where we start with an abstract (free) algebra, 
and we find a non-trivial homomorphism from this algebra
to a concrete algebra based on Kripke systems (an algebraic version of Kripke frames).

The idea (at least for the equality-free case) is that we start with a theory 
(which is defined as a pair of sets of formulas, as is the case with classical intuitionistic logic), extend 
it to a saturated one in enough spare dimensions, or an appropraite dilation (lemma \ref{t2}),  
and then iterate this process countably many times forming  consecutive (countably many) dilations
in enough spare dimensions, using pairs of pairs (theories), cf. lemma \ref{t3}; finally forming an
 extension that will be used  to construct desired Kripke models (theorem \ref{main}). 
The extensions constructed are essentially conservative extensions, 
and they will actually constitute the set of worlds of our desired Kripke model.

The iteration is done by a subtle 
zig-zag process, a technique due to 
Gabbay \cite{b}. When we have diagonal elements (equality), constructing desired Kripke model, 
is substantialy different, and much more intricate.

All definitions and results up to lemma \ref{main1}, 
though formulated only for the diagonal-free case, applies equally well to the case when there are diagonal elements, 
with absolutely no modifications. (The case when diagonal elements are present will be dealt with in part 2).

\begin{definition} Let $\A\in GPHA_{\alpha}$. 
\begin{enumarab}
\item  A theory in $\A$ is a pair $(\Gamma, \Delta)$ such that $\Gamma, \Delta\subseteq \A$.
\item A theory $(\Gamma, \Delta)$ is consistent if there are no $a_1,\ldots a_n\in \Gamma$ and 
$b_1,\ldots b_m\in \Delta$ ($m,n\in \omega$) such that
$$a_1\land\ldots a_n\leq b_1\lor\ldots b_m.$$
Not that in this case, we have $\Gamma\cap \Delta=\emptyset$. Also if $F$ is a filter 
(has the finite intersection property), then it is always the case that
$(F, \{0\})$ is consistent.
\item A theory $(\Gamma, \Delta)$ is complete if for all $a\in A,$ either $a\in \Gamma$ or $a\in \Delta$.
\item A theory $(\Gamma, \Delta)$ is saturated if for all $a\in A$ and $j\in \alpha$, 
if ${\sf c}_ja\in \Gamma$,
then there exists $k\in \alpha\sim \Delta a$, such that ${\sf s}^j_ka\in \Gamma$.
Note that a saturated theory depends only on $\Gamma$.
\end{enumarab}
\end{definition}

\begin{lemma} \label{t1}Let $\A\in GPHA_{\alpha}$ and $(\Gamma,\Delta)$ be a consistent theory. 
\begin{enumroman}
\item For any $a\in A,$ either $(\Gamma\cup \{a\}, \Delta)$ or $(\Gamma, \Delta\cup\{a\})$ is consistent.
\item $(\Gamma,\Delta)$ can be extended to a complete theory in $\A.$
\end{enumroman}
\end{lemma}
\begin{demo}{Proof} \begin{enumroman}
\item Cf.  \cite{Hung}. Suppose for contradiction that both theories are inconsistent.
Then we have $\mu_1\land a\leq \delta_1$ and $\mu_2\leq a\land \delta_2$ where 
$\mu_1$ and $\mu_2$ are some conjunction of elements of $\Gamma$
and $\delta_1$, $\delta_2$ are some disjunction of elements of $\Delta$.
But from 
$(\mu_1\land a\to \delta_1)\land (\mu_2\to a\lor \delta_2)\leq (\mu_1\land \mu_2\to \delta_1\lor \delta_2),$
we get
$\mu_1\land \mu_2\leq \delta_1\lor \delta_2,$ which contradicts the consistency of $(\Gamma, \Delta)$.

\item Cf. \cite{Hung}. Assume that $|A|=\kappa$. Enumerate the elements of $\A$ as $(a_i:i<\kappa)$. 
Then we can extend $(\Gamma, \Delta)$ consecutively by adding $a_i$ either
to $\Gamma$ or $\Delta$ while preserving consistency. In more detail, we define by transfinite induction a sequence of 
theories $(\Gamma_i,\Delta_i)$ for 
$i\in \kappa$ as follows. Set $\Gamma_0=\Gamma$ and $\Delta_0=\Delta$. If $\Gamma_i,\Delta_i$ are defined for all $i<\mu$ 
where $\mu$ is a limit ordinal, let $\Gamma_{\mu}=(\bigcup_{i\in \mu} \Gamma_i, \bigcup_{i\in \mu} \Delta_i)$. Now for successor ordinals. 
Assume that $(\Gamma_i, \Delta_i)$ are defined.
Set $\Gamma_{i+1}=\Gamma_i\cup \{a_i\}, \Delta_{i+1}=\Delta_i$ in case this is consistent, else 
set $\Gamma_{i+1}=\Gamma_i$ and $\Delta_{i+1}=\Delta_i\cup \{a_i\}$.
Let $T=\bigcup_{i\in \kappa}T_i$ and  $F= \bigcup_{i\in \kappa} F_i$, then $(T, F)$ is as desired.
\end{enumroman}
\end{demo}

\begin{lemma}\label{t2} Let $\A\in GPHA_{\alpha}$ and $(\Gamma,\Delta)$  be a consistent theory of $\A$. 
Let $I$ be a set such that $\alpha\subseteq I$  and  let $\beta=|I\sim \alpha|=\max(|A|, |\alpha|).$ Then there exists a minimal dilation
$\B$ of $\A$ of dimension $I$, and a theory $(T,F)$ in $\B$, 
extending $(\Gamma,\Delta)$ such that $(T,F)$ is saturated and complete.
\end{lemma}

\begin{demo}{Proof}
Let $I$ be provided as in the statement of the lemma. By lemma \ref{dl}, 
there exists $\B\in GPHA_I$ such that $\A\subseteq \Nr_{\alpha}\B$ and $\A$ generates $\B$.
We also have for all $X\subseteq \A$, $\Sg^{\A}X=\Nr_{\alpha}\Sg^{\B}X$.
Let $\{b_i:i<\kappa\}$ be an enumeration of the elements of $\B$; here $\kappa=|B|.$ 
Define by transfinite recursion a sequence $(T_i, F_i)$ for  $i<\kappa$ of
theories as follows. Set $T_0=\Gamma$ and $F_0=\Delta$. We assume inductively that 
$$|\beta\sim \bigcup_{x\in T_i} \Delta x\cup \bigcup_{x\in F_i}\Delta x|\geq \omega.$$
This is clearly satisfied for $F_0$ and $T_0$. Now we need to worry only about successor ordinals.
Assume that $T_i$ and $F_i$ are defined. 
We distinguish between two cases:
\begin{enumerate}
\item $(T_i, F_i\cup \{b_i\})$ is consistent.
Then set $T_{i+1}=T_i$ and $F_{i+1}=F_i\cup \{b_i\}.$
\item  If not, that is if $(T_i, F_i\cup \{b_i\})$ is inconsistent.
In this case, we distinguish between two subcases:

(a) $b_i$ is not of the form ${\sf c}_jp.$
Then set
$T_{i+1}=T_i\cup \{b_i\}$ and $F_{i+1}=F_i$.

(b) $b_i={\sf c}_jp$ for some $j\in I$.
Then set
$T_{i+1}=T_i\cup \{{\sf c}_jp, {\sf s}_u^jp\}$ where $u\notin \Delta p\cup \bigcup_{x\in T_i}\cup \bigcup_{x\in F_i}\Delta x$ and $F_{i+1}=F_i$.
\end{enumerate}
Such a $u$ exists by the inductive assumption. Now we check by induction that each $(T_i, F_i)$ is consistent.
The only part that needs checking, in view of the previous lemma,  is subcase (b).
So assume that $(T_i,F_i)$ is consistent and $b_i={\sf c}_jp.$
If $(T_{i+1}, F_{i+1})$ is inconsistent, then we would have for some $a\in T_i$ and some $\delta\in F_i$ that
$a\land {\sf c}_jp\land {\sf s}_u^jp\leq \delta.$ From this we get
$a\land {\sf c}_jp\leq \delta,$ because ${\sf s}_u^jp\leq {\sf c}_jp.$
But this contradicts the consistency of 
$(T_i\cup \{{\sf c}_jp\}, F_i)$.
Let $T=\bigcup_{i\in \kappa}T_i$ and $F=\bigcup_{i\in \kappa} F_i$, then $(T,F)$ is consistent. 
We show that it is saturated. If ${\sf c}_jp\in T$, then ${\sf c}_jp\in T_{i+1}$ for some $i$,
hence ${\sf s}_u^jp\in T_{i+1}\subseteq T$ and $u\notin \Delta p$.
Now by lemma \ref{t1}, we can extend $(T,F)$ is $\B$ to a complete theory, 
and this will not affect saturation, since the process of completion does not take us out of 
$\B$.
\end{demo}

The next lemma constitutes the core of our construction; involving a zig-zag Gabbay construction, 
it will be used repeatedly, to construct our desired representation
via a set algebra based on a Kripke system defined in  \ref{Kripke}

\begin{lemma}\label{t3} Let $\A\in GPHA_{\alpha}$ be generated by $X$ and let $X=X_1\cup X_2$.  
Let $(\Delta_0, \Gamma_0)$, $(\Theta_0, \Gamma_0^*)$ be two consistent theories in $\Sg^{\A}X_1$ and $\Sg^{\A}X_2,$ respectively
such that $\Gamma_0\subseteq \Sg^{\A}(X_1\cap X_2)$, $\Gamma_0\subseteq \Gamma_0^*$. Assume further that 
$(\Delta_0\cap \Theta_0\cap \Sg^{\A}X_1\cap \Sg^{\A}X_2, \Gamma_0)$ is complete in $\Sg^{\A}X_1\cap \Sg^{\A}X_2$. 
Suppose that $I$ is a set such that $\alpha\subseteq I$ and $|I\sim \alpha|=max (|A|,|\alpha|)$. 
Then there exist a dilation $\B\in GPHA_I$ of $\A$, and theories $T_1=(\Delta_{\omega}, \Gamma_{\omega})$, 
$T_2=(\Theta_{\omega}, \Gamma_{\omega}^*)$ extending 
$(\Delta_0, \Gamma_0)$, $(\Theta_0, \Gamma_0^*)$, such that $T_1$ and $T_2$ are consistent and saturated in 
$\Sg^{\B}X_1$ and $\Sg^{\B}X_2,$ respectively,  
$(\Delta_{\omega}\cap \Theta_{\omega}, \Gamma_{\omega})$ is complete in $\Sg^{\B}X_1\cap \Sg^{\B}X_2,$ and
$\Gamma_{\omega}\subseteq \Gamma_{\omega}^*$.
\end{lemma}

\begin{demo}{Proof}  Like the corresponding proof in \cite{Hung}, we will build the desired theories 
in a step-by-step zig-zag manner in a large enough dilation whose dimension is specified 
by $I$. The spare dimensions play a role of added witnesses, that will allow us to eliminate quantifiers, in a sense.
Let $\A=\A_0\in GPHA_{\alpha}$. The proof consists of an iteration of lemmata \ref{t1} and \ref{t2}.
Let $\beta=max(|A|, |\alpha|)$, and let $I$ be such that $|I\sim \alpha|=\beta$.

We distinguish between two cases:

\begin{enumarab}

\item Assume  that $G$ is strongly rich or $G$ contains consists of all finite transformations. In this case we only deal with minimal dilations.
We can write  
$\beta = I\sim \alpha$ as $\bigcup_{n=1}^{\infty}C_n$ where $C_i\cap C_j=\emptyset$ for distinct $i$ and $j$ and $|C_i|=\beta$
for all $i$. Then iterate first two items in lemma \ref{dl}.
Let $\A_1=\A(C_1)\in G_{\alpha\cup C_1}PHA_{\alpha\cup C_1}$ be a minimal dilation of $\A$, so that $\A=\Nr_{\alpha}\A_1$.
Let $\A_2=\A(C_1)(C_2)$ be a minimal dilation of $\A_1$ so that $\A_1=\Nr_{\alpha\cup C_1}\A_2$. 
Generally, we define inductively $\A_n=\A(C_1)(C_2)\ldots (C_n)$ to be  a minimal dilation of $\A_{n-1}$, so that 
$\A_{n-1}=\Nr_{\alpha\cup C_1\cup \ldots C_{n-1}}\A_n$.
Notice that for $k<n$, $\A_n$ is a minimal dilation of $\A_k$.
So we have a sequence of algebras 
$\A_0\subseteq \A_1\subseteq \A_2\ldots.$ Each element in the sequence is the minimal dilation of its preceding one.

\item $G$ contains all transformations. Here we shall have to use minimal extensions at the start, i.e at the first step of the iteration.  
We iterate lemma \ref{dl}, using items (3) and (4) in lemma \ref{cylindrify} by taking $|C_1|=\beta$, 
and $|C_i|=\omega$ for all $i\geq 2$; this will yield the desired sequence of extensions.

\end{enumarab}
Now that we have a sequence of extensions $\A_0\subseteq \A_1\ldots$ in different increasing dimensions,  
we now form a limit of this sequence
in $I$ dimensions. We can use ultraproducts, but instead we use products, and quotient algebras. First form the Heyting algebra, that is the product of the Heyting reducts of the constructed algebras, that is take   
$\C=\prod_{n=0}^{\infty}\Rd A_n$, where $\Rd \A_n$ denotes the Heyting reduct of $\A_n$ obtained by discarding substitutions and cylindrifiers. 
Let 
$$M=\{f\in C: (\exists n\in \omega)(\forall k\geq n) f_{k}=0\}.$$
Then $M$ is a Heyting ideal of $\C$. Now form the quotient Heyting  algebra $\D=\C/M.$
We want to expand this Heyting algebra algebra by cylindrifiers and substitutions, i.e to an algebra in $GPHA_{I}$.
Towards this aim, for $\tau\in {}G,$ define
$\phi({\tau})\in {} ^CC$ as follows:
$$(\phi(\tau)f)_n={\sf s}_{\tau\upharpoonright dim \A_n}^{\A_n}f_n$$
if $\tau(dim(\A_n))\subseteq dim (\A_n)$.
Otherwise $$(\phi(\tau)f)_n=f_n.$$
For $j\in I$,  define
$${\sf c}_jf_n={\sf c}_{(dim \A_n\cap \{j\})}^{\A_n}f_n,$$
and
$${\sf q}_jf _n={\sf q}_{(dim \A_n\cap \{j\})}^{\A_n}f_n.$$
Then for $\tau\in G$ and $j\in I$,  set
$${\sf s}_{\tau}(f/M)=\phi({\tau})f/M,$$
$${\sf c}_{j}(f/M)=({\sf c}_j f)/M,$$
and
$${\sf q}_{j}(f/M)=({\sf q}_j f)/M.$$
Then, it can be easily checked that,  $\A_{\infty}=(\D, {\sf s}_{\tau}, {\sf c}_{j}, {\sf q}_{j})$ is a $GPHA_I$,
in which every $\A_n$ neatly embeds. We can and will assume that $\A_n=\Nr_{\alpha\cup C_1\ldots \cup C_n}\A_{\infty}$. 
Also $\A_{\infty}$ is a minimal dilation of $\A_n$ for all $n$.
During our 'zig-zagging' 
we shall be extensively using lemma \ref{cylindrify}.

From now on, fix $\A$ to be as in the statement of lemma \ref{t3} for some time to come. 
So $\A\in GPHA_{\alpha}$ is generated by $X$ and $X=X_1\cup X_2$.  
$(\Delta_0, \Gamma_0)$, $(\Theta_0, \Gamma_0^*)$ are two consistent theories in $\Sg^{\A}X_1$ and $\Sg^{\A}X_2,$ respectively
such that $\Gamma_0\subseteq \Sg^{\A}(X_1\cap X_2)$, $\Gamma_0\subseteq \Gamma_0^*$. Finally  
$(\Delta_0\cap \Theta_0\cap \Sg^{\A}X_1\cap \Sg^{\A}X_2, \Gamma_0)$ is complete in $\Sg^{\A}X_1\cap \Sg^{\A}X_2.$ 
Now we have:
$$\Delta_0\subseteq \Sg^{\A}X_1\subseteq \Sg^{\A(C_1)}X_1\subseteq \Sg^{\A(C_1)(C_2)}X_1\subseteq \Sg^{\A(C_1)(C_2)(C_3)}X_1 \ldots\subseteq 
\Sg^{\A_{\infty}}X_1.$$ 
$$\Theta_0\subseteq \Sg^{\A}X_2\subseteq \Sg^{\A(C_1)}X_2\subseteq 
\Sg^{\A(C_1)(C_2)}X_2\subseteq \Sg^{\A(C_1)(C_2)(C_3)}X_2 \ldots\subseteq \Sg^{\A_{\infty}}X_2. $$
In view of lemmata  \ref{t1}, \ref{t2}, extend $(\Delta_0, \Gamma_0)$ 
to a complete and saturated theory $(\Delta_1, \Gamma_1')$ in $\Sg^{\A(C_1)}X_1$. Consider $(\Delta_1, \Gamma_0)$. Zig-zagging away,
we extend our theories in a step by step manner. The proofs of the coming Claims, 1, 2 and 3, 
are very similar to the proofs of the corresponding claims in \cite{Hung}, which are in turn an algebraic version of lemmata 4.18-19-20 in \cite {b},
with one major difference from the former. 
In our present situation, we can cylindrify on only finitely many indices, so we have to be 
careful, when talking about dimension sets, and in forming neat reducts (or compressions). 
Our proof then becomes substantially more involved. In the course of our proof we use extensively lemmata 
\ref{dl} and \ref{cylindrify} which are not formulated
in \cite{Hung} because we simply did not need them when we had cylindrifications on possibly 
infinite sets.

\begin{athm}{Claim 1} The theory $T_1=(\Theta_0\cup (\Delta_1\cap \Sg^{\A(C_1)}X_2), \Gamma_0^*)$ is consistent in $\Sg^{\A(C_1)}X_2.$
\end{athm}
\begin{demo}{Proof of Claim 1}
Assume that $T_1$ is inconsistent. Then 
for some conjunction $\theta_0$ of elements in $\Theta_0$, some $E_1\in \Delta_1\cap \Sg^{\A(C_1)}X_2,$
and some 
disjunction $\mu_0^*$ in $\Gamma_0^*,$ we have
$\theta_0\land E_1\leq \mu_0^*,$
and so 
$E_1\leq \theta_0\rightarrow \mu_0^*.$
Since $\theta_0\in \Theta_0\subseteq \Sg^{\A}X_2$ 
and $\mu_0^*\in \Gamma_0^*\subseteq \Sg^{\A}X_2\subseteq \Nr_{\alpha}^{\A(C_1)}\A$, 
therefore, for any finite 
set $D\subseteq C_1\sim \alpha$, we have ${\sf c}_{(D)}\theta_0=\theta_0$ and ${\sf c}_{(D)}\mu_0^*=\mu_0^*$.
Also for any finite set $D\subseteq C_1\sim \alpha,$ we have 
${\sf c}_{(D)}E_1\leq {\sf c}_{(D)}(\theta_0\to \mu_0^*)=\theta_0\to \mu^*.$
Now $E_1\in \Delta_1$, hence $E_1\in \Sg^{\A(C_1)}X_1$. 
By definition, we also have $E_1\in \Sg^{\A(C_1)}X_2.$
By lemma \ref{cylindrify} there exist finite sets $D_1$ and $D_2$ contained in $C_1\sim \alpha,$ such that
$${\sf c}_{(D_1)}E_1\in \Nr_{\alpha}\Sg^{\A(C_1)}X_1$$ and 
$${\sf c}_{(D_2)}E_1\in \Nr_{\alpha}\Sg^{\A(C_1)}X_2.$$
Le $D=D_1\cup D_2$. Then $D\subseteq  C_1\sim \alpha$ and we have: 
$${\sf c}_{(D)}E_1\in \Nr_{\alpha}\Sg^{\A(C_1)}X_1=\Sg^{\Nr_{\alpha}\A(C_1)}X_1=\Sg^{\A}X_1$$ and
$${\sf c}_{(D)}E_1\in \Nr_{\alpha}\Sg^{\A(C_1)}X_2=\Sg^{\Nr_{\alpha}\A(C_1)}X_2=\Sg^{\A}X_2,$$ 
that is to say
$${\sf c}_{(D)}E_1\in \Sg^{\A}X_1\cap \Sg^{\A}X_2.$$
Since  
$(\Delta_0\cap \Theta_0\cap \Sg^{\A}X_1\cap \Sg^{\A}X_2, \Gamma_0)$ is complete in $\Sg^{\A}X_1\cap \Sg^{\A}X_2,$
we get that ${\sf c}_{(D)}E_1$ is either in $\Delta_0\cap \Theta_0$ or $\Gamma_0$.
We show that either way leads to a contradiction, by which we will be done. Suppose it is in $\Gamma_0$. 
Recall that we extended $(\Delta_0, \Gamma_0)$ to a complete saturated extension 
$(\Delta, \Gamma')$ in $\Sg^{\A(C_1)}X_1$. Since $\Gamma_0\subseteq \Gamma_1',$ we get that ${\sf c}_{(D)}E_1\in \Gamma_1'$ hence
${\sf c}_{(D)}E_1\notin \Delta_1$  because $(\Delta_1,\Gamma_1')$ is saturated and consistent. But this 
contradicts that $E_1\in \Delta_1$ because $E_1\leq {\sf c}_{(D)}E_1.$
Thus, we can infer that  ${\sf c}_{(D)}E_1\in \Delta_0\cap \Theta_0$. In particular, it is in $\Theta_0,$ and so 
$\theta_0\rightarrow \mu_0^*\in \Theta_0$. 
But again  this contradicts the consistency of $(\Theta_0, \Gamma_0^*)$.
\end{demo}

Now we extend $T_1$ to a complete and saturated theory $(\Theta_2, \Gamma_2^*)$ in $\Sg^{\A(C_1)(C_2)}X_2$. Let 
$\Gamma_2=\Gamma_2^*\cap \Sg^{\A(C_1)(C_2)}X_1$. 

\begin{athm}{Claim 2} The theory $T_2=(\Delta_1\cup (\Theta_2\cap \Sg^{\A(C_1)(C_2))}X_1), \Gamma_2)$ 
is consistent in $\Sg^{\A(C_1)(C_2)}X_1$.
\end{athm}
\begin{demo}{Proof of Claim 2}  If the Claim fails to hold, then we would have some 
$\delta_1\in \Delta_1$, $E_2\in \Theta_2\cap \Sg^{\A(C_1)(C_2)}X_1,$  and a disjunction $\mu_2\in \Gamma_2$  such that  
$\delta_1\land E_2\rightarrow \mu_2,$
and so 
$\delta_1\leq (E_2\rightarrow \mu_2)$
since  $\delta_1\in \Delta_1\subseteq \Sg^{\A(C_1)}X_1$. But $\Sg^{\A(C_1)}X_1\subseteq \Nr_{\alpha\cup C_1}^{\A(C_1)(C_2)}X_1$, 
therefore for any finite set $D\subseteq C_2\sim C_1,$ 
we have ${\sf q}_{(D)}\delta_1=\delta_1.$
The following holds for any finite set $D\subseteq C_2\sim C_1,$
$$\delta_1\leq {\sf q}_{(D)}(E_2\rightarrow \mu_2).$$
Now, by lemma \ref{cylindrify}, there is a finite set $D\subseteq C_2\sim C_1,$ satisfying 
\begin{equation*}
\begin{split}
\delta_1\to {\sf q}_{(D)}(E_2\rightarrow \mu_2)
&\in \Nr_{\alpha\cup C_1}\Sg^{\A(C_1)(C_2)}X_2,\\
&=\Sg^{\Nr_{\alpha\cup C_1}\A(C_1)(\A(C_2)}X_2,\\
&=\Sg^{\A(C_1)}X_2.\\
\end{split}
\end{equation*}
Since $\delta_1\in \Delta_1$, and $\delta_1\leq {\sf q}_{(D)}(E_2\to \mu_2)$,  
we get that ${\sf q}_{(D)}(E_2\rightarrow \mu_2)$ is in $\Delta_1\cap \Sg^{\A(C_1)}X_2$. 
We proceed as in the previous claim replacing $\Theta_0$ by $\Theta_2$ and the existental quantifier by the universal one.
Let $E_1= {\sf q}_{(D)}(E_2\to \mu_2)$. Then $E_1\in \Sg^{\A(C_1)}X_1\cap \Sg^{\A(C_2)}X_2$.
By lemma \ref{cylindrify} there exist finite sets $D_1$ and $D_2$ contained in $C_1\sim \alpha$ such that
$${\sf q}_{(D_1)}E_1\in \Nr_{\alpha}\Sg^{\A(C_1)}X_1,$$ and 
$${\sf q}_{(D_2)}E_1\in \Nr_{\alpha}\Sg^{\A(C_1)}X_2.$$
Le $J=D_1\cup D_2$. Then $J\subseteq  C_1\sim \alpha,$ and we have: 
$${\sf q}_{(J)}E_1\in \Nr_{\alpha}\Sg^{\A(C_1)}X_1=\Sg^{\Nr_{\alpha}\A(C_1)}X_1=\Sg^{\A}X_1$$ and
$${\sf q}_{(J)}E_1\in \Nr_{\alpha}\Sg^{\A(C_1)}X_2=\Sg^{\Nr_{\alpha}\A(C_1)}X_2=\Sg^{\A}X_2.$$ 
That is to say,
$${\sf q}_{(J)}E_1\in \Sg^{\A}X_1\cap \Sg^{\A}X_2.$$
Now $(\Delta_0\cap \Theta_2\cap \Sg^{\A}X_1\cap \Sg^{\A}X_2, \Gamma_0)$ is complete in $\Sg^{\A}X_1\cap \Sg^{\A}X_2,$
we get that ${\sf q}_{(J)}E_1$ is either in $\Delta_0\cap \Theta_2$ or $\Gamma_0$.
Suppose it is in $\Gamma_0$. Since $\Gamma_0\subseteq \Gamma_1'$, we get that ${\sf q}_{(J)}E_1\in \Gamma_1'$, hence
${\sf q}_{(J)}E_1\notin \Delta_1,$  because $(\Delta_1,\Gamma_1')$ is saturated and consistent. 
Here, recall that,  $(\Delta, \Gamma')$ is a saturated complete extension of 
$(\Gamma, \Delta)$.  But this 
contradicts that $E_1\in \Delta_1$.
Thus, we can infer that  ${\sf q}_{(J)}E_1\in \Delta_0\cap \Theta_2$.  In particular, it is in $\Theta_2$. 
Hence ${\sf q}_{(D\cup J)}(E_2\to \mu_2)\in \Theta_2$, and so $E_2\to \mu_2\in \Theta_2$ since 
${\sf q}_{(D\cup J)}(E_2\to \mu_2)\leq E_1\to \mu_2$. But this is a contradiction,
since $E_2\in \Theta_2$, $\mu_2\in \Gamma_2^*$ and $(\Theta_2,\Gamma_2^*)$ is consistent. 
\end{demo}
Extend $T_2$ to a complete and saturated theory $(\Delta_3, \Gamma_3')$ in $\Sg^{\A(C_1)(C_2)(C_3)}X_1$
such that $\Gamma_2\subseteq \Gamma_3'$. Again we are interested only in $(\Delta_3, \Gamma_2)$.

\begin{athm}{Claim 3 } The theory $T_3=(\Theta_2\cup \Delta_3\cap \Sg^{\A(C_1)(C_2)(C_3)}X_2, \Gamma_2^*)$ is consistent in 
$\Sg^{\A(C_1)(C_2)(C_3)}X_2.$
\end{athm}
\begin{demo}{Proof of Claim 3}  Seeking a contradiction, assume that the Claim does not hold. Then we would get 
for some $\theta_2\in \Theta_2$, $E_3\in \Delta_3\cap \Sg^{\A(C_1)(C_2)(C_3)}X_2$ and some disjunction $\mu_2^*\in \Gamma_2^*,$  that
$\theta_2\land E_3\leq \mu_2^*.$
Hence $E_3\leq \theta_2\rightarrow \mu_2^*.$
For any finite set $D\subseteq C_3\sim (C_1\cup C_2),$ we have ${\sf c}_{(D)}E_3\leq \theta_2\rightarrow \mu_2^*$.
By lemma \ref{cylindrify},  there is a  finite set $D_3\subseteq C_3\sim (C_1\cup C_2),$ satisfying  
\begin{equation*}
\begin{split}
{\sf c}_{(D_3)}E_3
&\in \Nr_{\alpha\cup C_1\cup C_2}\Sg^{\A(C_1)(C_2)(C_3)}X_1\\
&=\Sg^{\Nr_{\alpha\cup C_1\cup C_2}\A(C_1)C_2)(C_3)}X_1\\
&= \Sg^{\A(C_1)(C_2)}X_1.\\ 
\end{split}
\end{equation*}
If ${\sf c}_{(D_3)}E_3\in \Gamma_2^*$, then it in $\Gamma_2$, 
and since $\Gamma_2\subseteq \Gamma_3'$, it cannot be in $\Delta_3$.
But this contradicts that $E_3\in \Delta_3$. So ${\sf c}_{(D_3)}E_3\in \Theta_2,$ because $E_3\leq {\sf c}_{(D_3)}E_3,$ 
and so $(\theta_2\rightarrow \mu_2^*)\in \Theta_2,$ which contradicts 
the consistency of $(\Theta_2, \Gamma_2^*).$
\end{demo}
 Likewise, now extend $T_3$ to a complete and saturated theory $(\Delta_4, \Gamma_4')$ in $\Sg^{\A(C_1)(C_2)(C_3)(C_4)}X_2$
such that $\Gamma_3\subseteq \Gamma_4'.$ As before the theory 
$(\Delta_3, \Theta_4\cap \Sg^{\A(C_1)(C_2)(C_3)(C_4)}X_1, \Gamma_4)$ is consistent in
$\Sg^{\A(C_1)(C_2)(C_3)(C_4)}X_1$.
Continue, inductively,  to construct $(\Delta_5, \Gamma_5')$, $(\Delta_5, \Gamma_4)$
and so on.
We obtain, zigzaging along, the following sequences: 
$$(\Delta_0, \Gamma_0), (\Delta_1, \Gamma_0), (\Delta_3, \Gamma_2)\ldots$$
$$(\Theta_0, \Gamma_0^*), (\Theta_2, \Gamma_2^*), (\Theta_4, \Gamma_4^*)\ldots $$
such that
\begin{enumarab}
\item $(\theta_{2n}, \Gamma_{2n}^*)$ is complete and saturated in $\Sg^{\A(C_1)\ldots (C_{2n})}X_2,$
\item $(\Delta_{2n+1}, \Gamma_{2n})$ is a saturated theory in $\Sg^{\A(C_1)\ldots (C_{2n+1})}X_1,$
\item $\Theta_{2n}\subseteq \Theta_{2n+2}$, $\Gamma_{2n}^*\subseteq \Gamma_{2n+2}^*$ and 
$\Gamma_{2n}=\Gamma_{2n}^*\cap \Sg^{\A(C_1)\ldots \A(C_{2n})}X_1,$
\item $\Delta_0\subseteq \Delta_1\subseteq \Delta_3\subseteq \ldots .$
\end{enumarab}
Now let $\Delta_{\omega}=\bigcup_{n}\Delta_n$, $\Gamma_{\omega}=\bigcup_{n}\Gamma_n$, $\Gamma_{\omega}^*=\bigcup_{n}\Gamma_n^*$ and
$\Theta_{\omega}=\bigcup_n\Theta_n$. Then we have
 $T_1=(\Delta_{\omega}, \Gamma_{\omega})$, 
$T_2=(\Theta_{\omega}, \Gamma_{\omega}^*)$ extend 
$(\Delta, \Gamma)$, $(\Theta, \Gamma^*)$, such that $T_1$ and $T_2$ are consistent and saturated in $\Sg^{\B}X_1$ and $\Sg^{\B}X_2,$ respectively, 
$\Delta_{\omega}\cap \Theta_{\omega}$ is complete in $\Sg^{\B}X_1\cap \Sg^{\B}X_2,$ and
$\Gamma_{\omega}\subseteq \Gamma_{\omega}^*$.
We check that $(\Delta_{\omega}\cap \Theta_{\omega},\Gamma_{\omega})$
is complete in $\Sg^{\B}X_1\cap \Sg^{\B}X_2$.
Let $a\in \Sg^{\B}X_1\cap \Sg^{\B}X_2$. Then there exists $n$ such that
$a\in \Sg^{\A_n}X_1\cap \Sg^{\A_n}X_2$. Now $(\Theta_{2n}, \Gamma_{2n}^*)$ is complete and so either $a\in \Theta_{2n}$ or $a\in \Gamma_{2n}^*$.
If $a\in \Theta_{2n}$ it will be in $\Delta_{2n+1}$ and if $a\in \Gamma_{2n}^*$ it will be in $\Gamma_{2n}$. In either case,
$a\in \Delta_{\omega}\cap \Theta_{\omega}$ or $a\in \Gamma_{\omega}$. 
\end{demo}
\begin{definition}
\begin{enumarab}
\item  Let $\A$ be an algebra generated by $X$ and assume that $X=X_1\cup X_2$. A pair $((\Delta,\Gamma)$  $(T,F))$ of theories in 
$\Sg^{\A}X_1$ and $\Sg^{\A}X_2$ is a matched pair of theories if 
$(\Delta\cap T\cap \Sg^{\A}X_1\cap \Sg^{\A}X_2, \Gamma\cap F\cap \Sg^{\A}X_1\cap \Sg^{\A}X_2)$
is complete in $\Sg^{\A}X_1\cap \Sg^{\A}X_2$.
\item A theory $(T, F)$ extends a theory $(\Delta, \Gamma)$ if $\Delta\subseteq T$ and $\Gamma\subseteq F$.
\item A pair $(T_1, T_2)$ of theories extend another pair  $(\Delta_1, \Delta_2)$ if $T_1$ 
extends $\Delta_1$ and $T_2$ extends $\Delta_2.$ 
\end{enumarab}
\end{definition}
The following Corollary follows directly from the proof of lemma \ref{t3}.

\begin{corollary}\label{main1} Let $\A\in GPHA_{\alpha}$ be generated by $X$ and let $X=X_1\cup X_2$.  
Let $((\Delta_0, \Gamma_0)$, $(\Theta_0, \Gamma_0^*))$ be a matched pair in $\Sg^{\A}X_1$ and $\Sg^{\A}X_2,$ respectively.
Let $I$ be  a set such that $\alpha \subseteq I$, and $|I\sim \alpha|=max(|A|, |\alpha|)$.
Then there exists a dilation $\B\in GPHA_I$ of $\A$, and a matched pair, $(T_1, T_2)$ extending 
$((\Delta_0, \Gamma_0)$, $(\Theta_0, \Gamma_0^*))$, such that $T_1$ and $T_2$ are saturated in $\Sg^{\B}X_1$ and $\Sg^{\B}X_2$, respectively.
\end{corollary}
We next define set algebras based on Kripke systems. We stipulate that ubdirect products (in the univerasl algebraic sense) 
are the representable algebras, which the abstract axioms formulated in ? aspire to capture.
Here Kripke systems (a direct generalization of Kripke frames) are defined differently 
than those defined in \cite{Hung}, because we allow {\it relativized} semantics. In the clasical case, such algebras reduce to products of set algebras.
\footnote{The idea of relativization, similar to Henkin's semantics for second order logic, has proved a very 
fruitful idea in the theory of cylindric algebras.} 

\begin{definition}
Let $\alpha$ be an infinite set. 
A Kripke system of dimension $\alpha$ is a quadruple $\mathfrak{K}=(K, \leq \{X_k\}_{k\in K}, \{V_k\}_{k\in K}),$ such that
$V_k\subseteq {}^{\alpha}X_k,$
and
\begin{enumarab}
\item $(K,\leq)$ is preordered set,
\item For any $k\in K$, $X_k$ is a non-empty set such that
$$k\leq k'\implies X_k\subseteq  X_{k'}\text { and } V_k\subseteq V_{k'}.$$
\end{enumarab}
\end{definition}\label{Kripke} 
Let $\mathfrak{O}$ be the Boolean algebra $\{0,1\}$. Now Kripke systems define concrete polyadic Heyting algebras as follows.
Let $\alpha$ be an infinite set and $G$ be a 
semigroup of transformations on $\alpha$. Let $\mathfrak{K}=(K,\leq \{X_k\}_{k\in K}, \{V_k\}_{k\in K})$ 
be a Kripke system.
Consider the set
$$\mathfrak{F}_{\mathfrak{K}}=\{(f_k:k\in K); f_k:V_k\to \mathfrak{O}, k\leq k'\implies f_k\leq f_{k'}\}.$$
If $x,y\in {}^{\alpha}X_k$ and $j\in \alpha$ we write $x\equiv_jy$ if $x(i)=y(i)$ for all $i\neq j$.
We write $(f_k)$ instead of $(f_k:k\in K)$.
In $\mathfrak{F}_{\mathfrak{K}}$ we introduce the following operations:
$$(f_k)\lor (g_k)=(f_k\lor g_k)$$
$$(f_k)\land (g_k)=(f_k\land g_k.)$$
For any $(f_k)$ and $(g_k)\in \mathfrak{F}$, define
$$(f_k)\rightarrow (g_k)=(h_k),$$
where
$(h_k)$ is given for $x\in V_k$ by $h_k(x)=1$ if and only if for any $k'\geq k$ if $f_{k'}(x)=1$ then $g_{k'}(x)=1$.
For any $\tau\in G,$  define
$${\sf s}_{\tau}:\mathfrak{F}\to \mathfrak{F}$$ by
$${\sf s}_{\tau}(f_k)=(g_k)$$
where
$$g_k(x)=f_k(x\circ \tau)\text { for any }k\in K\text { and }x\in V_k.$$
For any $j\in \alpha$ and $(f_k)\in \mathfrak{F},$
define
$${\sf c}_{j}(f_k)=(g_k),$$
where for $x\in V_k$
$$g_k(x)=\bigvee\{f_k(y): y\in V_k,\  y\equiv_j x\}.$$
Finally, set
$${\sf q}_{j}(f_k)=(g_k)$$
where for $x\in V_k,$
$$g_k(x)=\bigwedge\{f_l(y): k\leq l, \ y\in V_k, y\equiv_j x\}.$$

The diagonal element ${\sf d}_{ij}$ is defined to be the tuple $(f_k:k\in K)$ where for $x\in V_k$, $f_k(x)=1$ iff $x_i=x_j.$ 
 
The algebra $\F_{\bold K}$ is called the set algebra based on the Kripke system $\bold K$.

\subsection{Diagonal Free case}

Our next theorem addresses the cases of $GPHA_{\alpha}$ with $G$ a rich semigroup, 
and everything is countable, and the case when $G={}^{\alpha}\alpha$ with no restrictions on cardinality.
It is an algebraic version of a  version of Robinson's joint consistency theorem:
A pair of consistent theories that agree on their 
common part can be amalgamated by taking their union to form a consistent extension of both; 
however, we stipulate that the second component of the first theory is included in the second component of the second theory.
We will provide examples showing that we cannot omit this condition.
The case when $G$ consists of finite transformations will be dealt with separately.

It also says, that the results in \cite{Hung} proved for full polyadic Heyting algebras
remains valid when we restrict cylindrifications to be finite, possibly add  diagonal elements,  
and consider semigroups that could be finitely generated, showing that the presence of 
all  infinitary substitutions and infinitary cylindrifications 
is somewhat of an overkill. 

Indeed, the axiomatization of full polyadic Heyting algebras studied in \cite{Hung}
is extremely complex from the recursion theoretic point of view \cite{NS}, while the axiomatizations studied here are far less 
complex; indeed they are recursive.  This is definitely an acet from the algebraic point of view.

\begin{theorem}\label{main}
Let $\alpha$ be an infinite set. Let $G$ be a semigroup on $\alpha$ containing at least one infinitary transformation. 
Let $\A$ be the free  $G$ algebra generated by $X$, and suppose that $X=X_1\cup X_2$.
Let $(\Delta_0, \Gamma_0)$, $(\Theta_0, \Gamma_0^*)$ be two consistent theories in $\Sg^{\A}X_1$ and $\Sg^{\A}X_2,$ respectively.
Assume that $\Gamma_0\subseteq \Sg^{\A}(X_1\cap X_2)$ and $\Gamma_0\subseteq \Gamma_0^*$.
Assume, further, that  
$(\Delta_0\cap \Theta_0\cap \Sg^{\A}X_1\cap \Sg^{\A}X_2, \Gamma_0)$ is complete in $\Sg^{\A}X_1\cap \Sg^{\A}X_2$. 
Then there exist a Kripke system $\mathfrak{K}=(K,\leq \{X_k\}_{k\in K}\{V_k\}_{k\in K}),$ a homomorphism $\psi:\A\to \mathfrak{F}_{\mathfrak K},$
$k_0\in K$, and $x\in V_{k_0}$,  such that for all $p\in \Delta_0\cup \Theta_0$ if $\psi(p)=(f_k)$, then $f_{k_0}(x)=1$
and for all $p\in \Gamma_0^*$ if $\psi(p)=(f_k)$, then $f_{k_0}(x)=0$.
\end{theorem}
\begin{demo}{Proof} We use lemma \ref{t3}, extensively. Assume that $\alpha$, $G$, $\A$ and $X_1$, $X_2$ and everything 
else in the hypothesis
are given.  Let $I$ be  a set containing $\alpha$ such that 
$\beta=|I\sim \alpha|=max(|A|, |\alpha|).$
If $G$ is strongly rich, let $(K_n:n\in \omega)$ be a family of pairwise disjoint sets such that $|K_n|=\beta.$
Define a sequence of algebras
$\A=\A_0\subseteq \A_1\subseteq \A_2\subseteq \A_2\ldots \subseteq \A_n\ldots,$
such that
$\A_{n+1}$ is a minimal dilation of $\A_n$ and $dim(\A_{n+1})=\dim\A_n\cup K_n$.

If $G={}^{\alpha}\alpha$, then let $(K_n:n\in \omega\}$ be a family of pairwise disjoint sets, such that $|K_1|=\beta$ 
and $|K_n|=\omega$ for $n\geq 1$, and define a sequence of algebras
$\A=\A_0\subseteq \A_1\subseteq \A_2\subseteq \A_2\ldots \subseteq \A_n\ldots,$
such that $\A_1$ is a minimal extension of $\A$, and 
$\A_{n+1}$ is a minimal dilation of $\A_n$ for $n\geq 2$, with  $dim(\A_{n+1})=\dim\A_n\cup K_n$.

We denote $dim(\A_n)$ by $I_n$ for $n\geq 1$. Recall that $dim(\A_0)=\dim\A=\alpha$.

We interrupt the main stream of the proof by two consecutive claims. Not to digress, it might be useful that the reader 
at first reading, only memorize their statements, skip their proofs, go on with the main proof, and then get back to them.
The proofs of Claims 1 and 2 to follow are completely analogous to the corresponding claims in \cite{Hung}. 
The only difference is that we deal with only finite 
cylindrifiers, and in this respect they are closer to the proofs of lemmata 4.22-23 in \cite{b}. Those two claims are essential 
in showing that the maps that will be defined shortly into concrete set algebras based on appropriate 
Kripke systems, defined via pairs of theories, in increasing extensions (dimensions),
are actually homomorphisms. In fact, they have to do with the preservation of 
the operations of implication and universal quantification. The two claims use lemma \ref{t3}.

\begin{athm}{Claim 1} Let $n\in \omega$. If $((\Delta, \Gamma), (T,F))$ 
is a matched pair of saturated theories in $\Sg^{\A_n}X_1$ and $\Sg^{\A_n}X_2$, then the following hold.
For any $a,b\in \Sg^{\A_n}X_1$ if $a\rightarrow b\notin \Delta$, then there is a matched pair $((\Delta',\Gamma'), (T', F'))$ of saturated theories in
$\Sg^{\A_{n+1}}X_1$ and $\Sg^{\A_{n+1}}X_2,$ respectively, such that $\Delta\subseteq \Delta'$, $T\subseteq T',$  
$a\in \Delta'$ and $b\notin \Delta'$.
\end{athm}
\begin{demo}{Proof of Claim 1}  Since $a\rightarrow b\notin \Delta,$
 we have $(\Delta\cup\{a\}, b)$ is consistent in $\Sg^{\A_n}X_1$.
Then by lemma \ref{t1}, it can be extended to a complete theory $(\Delta', T')$ in $\Sg^{\A_n}X_1$.
Take $$\Phi=\Delta'\cap \Sg^{\A_n}X_1\cap \Sg^{\A_n}X_2,$$ 
and $$\Psi=T'\cap \Sg^{\A_n}X_1\cap \Sg^{\A_n}X_2.$$
Then $(\Phi, \Psi)$ is complete in $\Sg^{\A_n}X_1\cap \Sg^{\A_n}X_2.$
We shall now show that $(T\cup \Phi, \Psi)$ is consistent in $\Sg^{\A_n}X_2$. If not, then there is $\theta\in T$, 
$\phi\in \Phi$ and $\psi\in \Psi$ such that
$\theta\land \phi\leq \psi$. So $\theta\leq \phi\rightarrow \psi$. Since $T$ is saturated, we get that $\phi\rightarrow \psi$ is in $T$.
Now $\phi\rightarrow \psi\in \Delta\cap \Sg^{\A_n}X_1\cap \Sg^{\A_n}X_2\subseteq \Delta'\cap \Sg^{\A}X_1\cap \Sg^{\A}X_2=\Phi$.
Since $\phi\in \Phi$ and $\phi\rightarrow \psi\in \Phi,$ we get that $\psi\in \Phi\cap \Psi$. But this means that 
$(\Phi, \Psi)$ is inconsistent which is impossible. Thus $(T\cup \Phi, \Psi)$ is consistent.
Now the pair $((\Delta', T') (T\cup \Phi, \Psi))$ satisfy the conditions of lemma \ref{t3}. 
Hence this pair can be extended to a matched pair of 
saturated theories in $\Sg^{\A_{n+1}}X_1$ and $\Sg^{\A_{n+1}}X_2$. 
This pair is as required by the conclusion of lemma \ref{t3}.
\end{demo}

\begin{athm}{Claim 2} Let $n\in \omega$. If $((\Delta, \Gamma), (T,F))$ is a matched pair of saturated theories in $\Sg^{\A_n}X_1$ and $\Sg^{\A_n}X_2$, 
then the following hold.
For $x\in \Sg^{\A_n}X_1$ and $j\in I_n=dim\A_n$,  if ${\sf q}_{j}x\notin  \Delta$, then there is a matched pair $((\Delta',\Gamma'), (T', F'))$ 
of saturated theories in
$\Sg^{\A_{n+2}}X_1$ and $\Sg^{\A_{n+2}}X_2$ respectively,  $u\in I_{n+2}$
 such that $\Delta\subseteq \Delta'$, $T\subseteq T'$ and ${\sf s}_u^j x\notin \Delta'$. 
\end{athm}
\begin{demo}{Proof}
Assume that $x\in \Sg^{\A_n}X_1$ and $j\in I_n$ such that ${\sf q}_{j}x\notin  \Sg^{\A_n}X_1$. Then there exists 
$u\in I_{n+1}\sim I_n$ such that $(\Delta, {\sf s}_u^j x)$ is consistent in $\Sg^{\A_{n+1}}X_1$.
So $(\Delta, {\sf s}_u^jx)$ can be extended to a complete theory $(\Delta', T')$ in $\Sg^{\A_{n+1}}X_1$.
Take $$\Phi=\Delta'\cap \Sg^{\A_{n+1}}X_1\cap \Sg^{\A_{n+1}}X_2,$$
and
$$\Psi=T'\cap \Sg^{\A_{n+1}}X_1\cap \Sg^{\A_{n+1}}X_2.$$
Then $(\Phi,\Psi)$ is complete in $\Sg^{\A_{n+1}}X_1\cap \Sg^{\A_{n+1}}X_2$.
We shall show that $(T\cup \Phi,\Psi)$ is consistent in $\Sg^{\A_{n+1}}X_2$. If not, then there exist $\theta\in T,$ 
$\phi\in \Phi$ and $\psi\in \Psi,$ such that
$\theta\land \phi\leq \psi$. Hence, $\theta\leq \phi\rightarrow \psi$. Now
$$\theta={\sf q}_j(\theta)\leq  {\sf q}_{j}(\phi\rightarrow \psi).$$ 
Since $(T,F)$ is saturated in $\Sg^{\A_n}X_2,$ it thus follows that 
$${\sf q}_{j}(\phi\rightarrow \psi) \in T\cap \Sg^{\A_n}X_1\cap \Sg^{\A_n}X_2=\Delta\cap \Sg^{\A_n}X_1\cap \Sg^{\A_n}X_2.$$
So ${\sf q}_{j}(\phi\rightarrow \psi)\in \Delta'$ and consequently we get ${\sf q}_{j}(\phi\rightarrow \psi)\in \Phi$. Also, we have, $\phi\in \Phi$.
But $(\Phi, \Psi)$ is complete, we get $\psi\in \Phi$ and this contradicts that $\psi\in \Psi$.
Now the pair $((\Delta', \Gamma'), (T\cup \Phi, \Psi))$ satisfies the hypothesis of lemma \ref{t3} applied to 
$\Sg^{\A_{n+1}}X_1, \Sg^{\A_{n+1}}X_2$. 
The required now follows from the concusion of lemma \ref{t3}.
\end{demo}

Now that we have proved our claims, we go on with the proof.
We prove the theorem when $G$ is a strongly rich semigroup, because in this case we deal with relativized semantics,
and during the proof we state the necessary modifications for the case when $G$ is the semigroup of all transformations. 
Let $$K=\{((\Delta, \Gamma), (T,F)): \exists n\in \omega \text { such that } (\Delta, \Gamma), (T,F)$$
$$\text { is a a matched pair of saturated theories in }
\Sg^{\A_n}X_1, \Sg^{\A_n}X_2\}.$$
We have $((\Delta_0, \Gamma_0)$, $(\Theta_0, \Gamma_0^*))$ is a matched pair but the theories are not saturated. But by lemma \ref{t3}
there are $T_1=(\Delta_{\omega}, \Gamma_{\omega})$, 
$T_2=(\Theta_{\omega}, \Gamma_{\omega}^*)$ extending 
$(\Delta_0, \Gamma_0)$, $(\Theta_0, \Gamma_0^*)$, such that $T_1$ and $T_2$ are saturated in $\Sg^{\A_1}X_1$ and $\Sg^{\A_1}X_2,$ respectively.
Let $k_0=((\Delta_{\omega}, \Gamma_{\omega}), (\Theta_{\omega}, \Gamma_{\omega}^*)).$ Then $k_0\in K.$

If $i=((\Delta, \Gamma), (T,F))$ is a matched pair of saturated theories in $\Sg^{\A_n}X_1$ and $\Sg^{\A_n}X_2$, let $M_i=dim \A_n$, 
where $n$ is the least such number, so $n$ is unique to $i$.
Before going on we introduce a piece of notation. For a set $M$ and a sequence $p\in {}^{\alpha}M$, $^{\alpha}M^{(p)}$ is the following set
$$\{s\in {}^{\alpha}M: |\{i\in \alpha: s_i\neq p_i\}|<\omega\}.$$
Let $${\mathfrak{K}}=(K, \leq, \{M_i\}, \{V_i\})_{i\in \mathfrak{K}}$$
where $V_i=\bigcup_{p\in G_n}{}^{\alpha}M_i^{(p)}$, and $G_n$ is the strongly rich semigroup determining the similarity type of $\A_n$, with $n$ 
the least number such $i$ is a saturated matched pair in $\A_n$.
The order $\leq $ is defined as follows:
If $i_1=((\Delta_1, \Gamma_1)), (T_1, F_1))$ and $i_2=((\Delta_2, \Gamma_2), (T_2,F_2))$ are in $\mathfrak{K}$, then define
$$i_1\leq i_2\Longleftrightarrow  M_{i_1}\subseteq M_{i_2}, \Delta_1\subseteq \Delta_2, T_1\subseteq T_2.$$
This is, indeed as easily checked,  a preorder on $K$.

We  define two maps on $\A_1=\Sg^{\A}X_1$ and $\A_2=\Sg^{\A}X_2$ respectively, then those will be pasted using the freeness of $\A$
to give the required single homomorphism, by noticing that they agree on the common part, that is on $\Sg^{\A}(X_1\cap X_2).$

Set $\psi_1: \Sg^{\A}X_1\to \mathfrak{F}_{\mathfrak K}$ by
$\psi_1(p)=(f_k)$ such that if $k=((\Delta, \Gamma), (T,F))\in K$ is a matched pair of saturated theories in 
$\Sg^{\A_n}X_1$ and $\Sg^{\A_n}X_2$,
and $M_k=dim \A_n$, then for $x\in V_k=\bigcup_{p\in G_n}{}^{\alpha}M_k^{(p)}$,
$$f_k(x)=1\Longleftrightarrow {\sf s}_{x\cup (Id_{M_k\sim \alpha)}}^{\A_n}p\in \Delta\cup T.$$
To avoid tiresome notation, we shall denote the map $x\cup Id_{M_k\sim \alpha}$ simply by $\bar{x}$ when $M_k$ is clear from context.
It is easily verifiable that $\bar{x}$ is in the semigroup determining the similarity type of $\A_n$ hence the map is well defined.
More concisely, we we write $$f_k(x)=1\Longleftrightarrow {\sf s}_{\bar{x}}^{\A_n}p\in \Delta\cup T.$$
The map $\psi_2:\Sg^{\A}X_2\to \mathfrak{F}_{\mathfrak K}$ is defined in exactly the same way.
Since the theories are matched pairs, $\psi_1$ and $\psi_2$ agree on the common part, i.e. on $\Sg^{\A}(X_1\cap X_2).$
Here we also make the tacit assumption that if $k\leq k'$ then $V_k\subseteq V_{k'}$ 
via the embedding $\tau\mapsto \tau\cup Id$.

When $G$ is the semigroup of all transformations, with no restrictions on cardinalities, 
we need not relativize since $\bar{\tau}$ is in the big semigroup. In more detail, 
in this case, we take for $k=((\Delta,\Gamma), (T,F))$ a matched pair of saturated theories in 
$\Sg^{\A_n}X_1,\Sg^{\A_n}X_2$,
$M_k=dim\A_n$ and $V_k={}^{\alpha}M_k$
and for $x\in {}^{\alpha}M_k$, we set 
$$f_k(x)=1\Longleftrightarrow {\sf s}_{x\cup (Id_{M_k\sim \alpha)}}^{\A_n}p\in \Delta\cup T.$$

Before proving that $\psi$ is a homomorphism, we show that 
$$k_0=((\Delta_{\omega},\Gamma_{\omega}), (\Theta_{\omega}, \Gamma^*_{\omega}))$$ 
is as desired. Let $x\in V_{k_0}$ be the identity map. Let $p\in \Delta_0\cup \Theta_0$, then 
${\sf s}_xp=p\in \Delta_{\omega}\cup \Theta_{\omega},$ 
and so if $\psi(p)=(f_k)$ then $f_{k_0}(x)=1$.
On the other hand if $p\in \Gamma_0^*$, then $p\notin \Delta_{\omega}\cup \Theta_{\omega}$, and so $f_{k_0}(x)=0$.
Then the union $\psi$ of $\psi_1$ and $\psi_2$, $k_0$ and $Id$ are as required, modulo proving that 
$\psi$ is a homomorphism from $\A$, to the set algebra based on the above defined
Kripke system, which we proceed to show. We start by $\psi_1$.
Abusing notation, we denote $\psi_1$ by $\psi$, and we write a matched pair in $\A_n$ 
instead of a matched pair of saturated theories in $\Sg^{\A_n}X_1$,
$\Sg^{\A_n}X_2$, since $X_1$ and $X_2$ are fixed. The proof that the postulated map is a homomorphism is 
similar to the proof in \cite{Hung} baring in mind that it is far from being identical because 
cylindrifiers and their duals are only finite.
\begin{enumroman}
\item We prove that $\psi$ preserves $\land$. Let $p,q\in A$. Assume that
$\psi(p)=(f_k)$ and $\psi(q)=(g_k)$. Then $\psi(p)\land \psi(q)=(f_k\land g_k)$.
We now compute $\psi(p\land q)=(h_k)$
Assume that  $x\in V_k$, where $k=((\Delta,\Gamma), (T, F))$ is a matched pair in $\A_n$ and  $M_k=dim\A_n$.
Then
$$h_k(x)=1\Longleftrightarrow {\sf s}_{\bar{x}}^{\A_n}(p\land q)\in \Delta\cup T$$
$$\Longleftrightarrow {\sf s}_{\bar{x}}^{\A_n}p\land {\sf s}_{\bar{x}}^{\A_n}q\in \Delta\cup T$$
$$\Longleftrightarrow {\sf s}_{\bar{x}}^{\A_n}p\in T\cup \Delta\text { and }{\sf s}_{\bar{x}}^{\A_n}q\in \Delta\cup T$$
$$\Longleftrightarrow f_k(x)=1 \text { and } g_k(x)=1$$
$$\Longleftrightarrow (f_k\land g_k)(x)=1$$
$$\Longleftrightarrow (\psi(p)\land \psi(q))(x)=1.$$

\item $\psi$ preserves $\rightarrow.$ (Here we use Claim 1).  Let $p,q\in A$. Let $\psi(p)=(f_k)$ and $\psi(q)=(g_k)$.
Let $\psi(p\rightarrow q)=(h_k)$ and $\psi(p)\rightarrow \psi(q)=(h'_k)$.
We shall prove that for any $k\in \mathfrak{K}$ and any $x\in V_k$, we have
$$h_k(x)=1\Longleftrightarrow h'_k(x)=1.$$
Let $x\in V_k$. Then $k=((\Delta,\Gamma),(T,F))$ is a matched pair in $\A_n$ and $M_k=dim\A_n$.
Assume that $h_k(x)=1$. Then we have $${\sf s}_{\bar{x}}^{\A_n}(p\rightarrow q)\in \Delta\cup T,$$
from which we get that $$(*) \ \ \ {\sf s}_{\bar{x}}^{\A_n}p\rightarrow {\sf s}_{\bar{x}}^{\A_n}q\in \Delta\cup T.$$
Let $k'\in K$ such that $k\leq k'$. Then $k'=((\Delta', \Gamma'), (T', F'))$ is a matched pair in $\A_m$ with $m\geq n$.
Assume that $f_{k'}(x)=1$. Then, by definition we have (**)
$${\sf s}_{\bar{x}}^{\A_m}p\in \Delta'\cup T'.$$
But $\A_m$ is a dilation of $\A_n$ and so 
$${\sf s}_{\bar{x}}^{\A_m}p={\sf s}_{\bar{x}}^{\A_n}p\text { and } {\sf s}_{\bar{x}}^{\A_m}q={\sf s}_{\bar{x}}^{\A_n}q.$$
From (*) we get that, 
$$ {\sf s}_{\bar{x}}^{\A_m}p\rightarrow {\sf s}_{\bar{x}}^{\A_m}q\in \Delta'\cup T'.$$
But, on the other hand,  from (**), we have ${\sf s}_{\bar{x}}^{\A_m}q\in \Delta'\cup T',$
so $$f_{k'}(x)=1\Longrightarrow g_{k'}(x)=1.$$
That is to say, we have  $h_{k'}(x)=1$.
Conversely, assume that $h_k(x)\neq 1,$
then
$$ {\sf s}_{\bar{x}}^{\A_n}p\rightarrow {\sf s}_{\bar{x}}^{\A_n}q\notin \Delta\cup T,$$
and consequently
$$ {\sf s}_{\bar{x}}^{\A_n}p\rightarrow {\sf s}_{\bar{x}}^{\A_n}q\notin \Delta.$$
From Claim 1, we get that there exists a matched pair $k'=((\Delta',\Gamma')((T',F'))$ in $\A_{n+2},$ such that
$$ {\sf s}_{\bar{x}}^{\A_{n+2}}p\in \Delta'\text { and } {\sf s}_{\bar{x}}^{\A_{n+2}}q\notin \Delta'.$$
We claim that ${\sf s}_{\bar{x}}^{\A_{n+2}}q\notin T'$, for otherwise, if it is in $T'$, then we would get that
$${\sf s}_{\bar{x}}^{\A_{n+2}}q\in \Sg^{\A_{n+2}}X_1\cap \Sg^{\A_{n+2}}X_2.$$
But $$(\Delta'\cap T'\cap \Sg^{\A_{n+2}}X_1\cap \Sg^{\A_{n+2}}X_2, \Gamma'\cap F'\cap\Sg^{\A_{n+2}}X_1\cap \Sg^{\A_{n+2}}X_2)$$ is complete 
in $\Sg^{\A_{n+2}}X_1\cap \Sg^{\A_{n+2}}X_2,$
and ${\sf s}_{\bar{x}}^{\A_{n+2}}q\notin \Delta'\cap T'$, hence it must be the case that 
 $${\sf s}_{\bar{x}}^{\A_{n+2}}q\in \Gamma'\cap F'.$$
In particular, we have
$${\sf s}_{\bar{x}}^{\A_{n+2}}q\in F',$$
which contradicts the consistency of $(T', F'),$ since by assumption ${\sf s}_x^{\A_{n+2}}q\in T'$.
Now we have
$${\sf s}_{\bar{x}}^{\A_{n+2}}q\notin \Delta'\cup T',$$
and
$${\sf s}_{\bar{x}}^{\A_{n+2}}p\in \Delta'\cup T'.$$
Since $\Delta'\cup T'$ extends $\Delta\cup T$, we get that $h_k'(x)\neq 1$.
\item $\psi$ preserves substitutions. Let $p\in \A$. Let $\sigma\in {}G$.
Assume that $\psi(p)=(f_k)$ and $\psi({\sf s}_{\sigma}p)=(g_k).$ 
Assume that $M_k=\dim\A_n$ where $k=((\Delta,\Gamma),(T,F))$ is a matched pair 
in $\A_n$.
Then, for $x\in V_k$, we have 
$$g_k(x)=1\Longleftrightarrow {\sf s}_{\bar{x}}^{\A_n}{\sf s}_{\sigma}^{\A}p\in \Delta\cup T$$
$$\Longleftrightarrow  {\sf s}_{\bar{x}}^{\A_n}{\sf s}_{\bar{\sigma}}^{\A_n}p\in \Delta\cup T$$
$$\Longleftrightarrow  {\sf s}_{\bar{x}\circ {\bar{\sigma}}}^{\A_n}p\in \Delta\cup T$$
$$\Longleftrightarrow {\sf s}_{\overline{x\circ \sigma}}^{\A_n}p\in \Delta\cup T$$
$$\Longleftrightarrow f_k(x\circ \sigma)=1.$$

\item $\psi$ preserves cylindrifications. Let $p\in A.$ 
Assume that $m\in I$  and assume that  $\psi({\sf c}_{m}p)=(f_k)$ and ${\sf c}_m\psi(p)=(g_k)$. 
Assume that $k=((\Delta,\Gamma),(T,F))$ is a matched pair in 
$\A_n$ and that $M_k=dim\A_n$.  Let $x\in V_k$. Then 
$$f_k(x)=1\Longleftrightarrow {\sf s}_{\bar{x}}^{\A_n}{\sf c}_{m}p\in \Delta\cup T.$$
We can assume that  $${\sf s}_{\bar{x}}^{\A_n}{\sf c}_{m}p\in \Delta.$$
For if not, that is if  $${\sf s}_{\bar{x}}^{\A_n}{\sf c}_{m}p\notin \Delta\text { and }  {\sf s}_{\bar{x}}^{\A_n}{\sf c}_{(m)}p\in T,$$
then  $${\sf s}_{\bar{x}}^{\A_n}{\sf c}_{m}p\in \Sg^{\A_n}X_1\cap \Sg^{\A_n}X_2, $$
but  $$(\Delta\cap T\cap \Sg^{\A_n}X_1\cap \Sg^{\A_n}X_2, \Gamma\cap F\cap \Sg^{\A_n}X_1\cap \Sg^{\A_n}X_2)$$ is complete in 
$\Sg^{\A_n}X_1\cap \Sg^{\A_n}X_2$, 
and $${\sf s}_{\bar{x}}^{\A_n}{\sf c}_{m}p\notin \Delta\cap T,$$ it must be the case that
$${\sf s}_{\bar{x}}^{\A_n}{\sf c}_{m}p\in \Gamma\cap F.$$
In particular, $${\sf s}_{\bar{x}}^{\A_n}{\sf c}_{m}p\in F.$$
But this contradicts the consistency of $(T,F)$.

Assuming that ${\sf s}_x{\sf c}_mp\in \Delta,$ we proceed as follows.
Let  $$\lambda\in \{\eta\in I_n: x^{-1}\{\eta\}=\eta\}\sim \Delta p.$$
Let $$\tau=x\upharpoonright I_n\sim\{m, \lambda\}\cup \{(m,\lambda)(\lambda, m)\}.$$
Then, by item (5) in theorem \ref{axioms}, we have
$${\sf c}_{\lambda}{\sf s}^{\A_n}_{\bar{\tau}}p={\sf s}_{\bar{\tau}}^{\A_n}{\sf c}_{m}p={\sf s}_{\bar{x}}^{\A_n}{\sf c}_mp\in \Delta.$$
We introduce a piece of helpful notation. For a function $f$, let $f(m\to u)$ is the function that agrees with $f$ except at $m$ 
where its value is $u$.
Since $\Delta$ is saturated, there exists $u\notin \Delta x$ such that ${\sf s}_u^{\lambda}{\sf s}_xp\in \Delta$, and so ${\sf s}_{(x(m\to u))}
p\in \Delta$. 
This implies that $x\in {\sf c}_mf(p)$ and so $g_k(x)=1$. Conversely, assume that $g_k(x)=1$ with 
$k=((\Gamma,\Delta))$, $(T,F))$ a matched pair in $\A_n$.
Let $y\in V_k$ such that $y\equiv_m x$ and  $\psi(p)y=1$. Then ${\sf s}_{\bar{y}}p\in \Delta\cup T$. 
Hence ${\sf s}_{\bar{y}}{\sf c}_mp\in \Delta\cup T$ and so ${\sf s}_{\bar{x}}{\sf c}_mp\in \Delta\cup T$, 
thus  $f_k(x)=1$ and we are done.

\item $\psi$ preserves universal quantifiers. (Here we use Claim 2). Let $p\in A$ and $m\in I$. 
Let $\psi(p)=(f_k)$, ${\sf q}_{m}\psi(p)=(g_k)$ and $\psi({\sf q}_{m}p)=(h_k).$
Assume that $h_k(x)=1$. We have $k=((\Delta,\Gamma), (T,F))$ is a matched pair in $\A_n$ and $x\in V_k$.
Then $${\sf s}_{\bar{x}}^{\A_n}{\sf q}_{m}p\in \Delta\cup T,$$
and so 
$${\sf s}_{\bar{y}}^{\A_n}{\sf q}_{m}p\in \Delta\cup T \text{ for all } y\in {}^IM_k, y\equiv_m x.$$ 
Let $k'\geq k$. Then $k'=((\Delta',\Gamma'), (T',F'))$ is a matched pair in $\A_l$ $l\geq n$, 
$\Delta\subseteq \Delta'$ and $T\subseteq T'.$ 
Since $p\geq {\sf q}_{m}p$ it follows that
$${\sf s}_{\bar{y}}^{\A_n}p\in \Delta'\cup T' \text{ for all } y\in {}^IM_k, y\equiv_mx.$$ 
Thus $g_k(x)=1$.
Now conversely, assume that $h_k(x)=0$, $k=((\Delta,\Gamma), (T,F))$ is a matched pair in $\A_n,$
then, we have
$${\sf s}_{\bar{x}}^{\A_n}{\sf q}_{m}p\notin \Delta\cup T,$$
and so $${\sf s}_{\bar{x}}^{\A_n}{\sf q}_{m}p\notin \Delta.$$
Let  $$\lambda\in \{\eta\in I_n: x^{-1}\{\eta\}=\eta\}\sim \Delta p.$$
Let $$\tau=x\upharpoonright I_n\sim\{m, \lambda\}\cup \{(m,\lambda)(\lambda, m)\}.$$
Then, like in the existential case, using polyadic axioms, 
we get
$${\sf q}_{\lambda}{\sf s}_{\tau}p={\sf s}_{\tau}{\sf q}_{m}p={\sf s}_{x}{\sf q}_mp\notin \Delta$$
Then there exists $u$ such that ${\sf s}_u^{\lambda}{\sf s}_xp\notin \Delta.$
So ${\sf s}_u^{\lambda}{\sf s}_xp\notin T$, for if it is, then by the previous reasoning since it is an element of 
$\Sg^{\A_{n+2}}X_1\cap \Sg^{\A_{n+2}}X_2$ and by completeness
of $(\Delta\cap T, \Gamma\cap F)$ we would reach a contradiction.
The we get that ${\sf s}_{(x(m\to u))}p\notin \Delta\cup T$
which means that $g_k(x)=0,$
and we are done. 
\end{enumroman}
\end{demo}

We now deal with the case when $G$ is the semigroup of all finite transformations on $\alpha$. In this case, 
we stipulate that $\alpha\sim \Delta x$ is infinite for all $x$ in algebras considered.
To deal with such a case, we need to define certain free algebras, called dimension restricted. 
Those algebras were introduced by Henkin, 
Monk and Tarski. 
The free algebras defined the usual way, will have
the dimensions sets of their elements equal to their dimension, but we do not want that. 
For a class $K$, ${\bf S}$ 
stands for the operation of forming subalgebras of $K$, ${\bf P}K$ that of forming direct products, 
and ${\bf H}K$ stands for the operation of taking homomorphic images.
In particular, for a class $K$, ${\bf HSP}K$ stands for the variety generated by $K$.    

Our dimension restricted free algebbras, are an instance of certain independently generated algebras, 
obtained by an appropriate relativization of the universal algebraic 
concept of free algebras. For an algebra $\A,$ we write
$R\in Con\A$ if $R$ is a congruence relation on $\A.$

\begin{definition} Assume that $K$ is a class of algebras of similarity $t$ and $S$ is any set of ordered pairs of words of $\Fr_{\alpha}^t,$
the absolutely free algebra of type $t$. Let
$$Cr_{\alpha}^{(S)}K=\cap \{R\in Con \Fr_{\alpha}^t, \Fr_{\alpha}^t/R\in SK, S\subseteq R\}$$
and let
$$\Fr_{\alpha}^{(S)}K=\Fr_{\alpha}^t/Cr_{\alpha}^{(S)}K.$$
$\Fr_{\alpha}^{(S)}K$ is called the free algebra over $K$ with $\alpha$ generators subject to the defining relations $S$.
\end{definition}
As a special case, we obtain dimension restricted free algebra, defined next.
\begin{definition}
\begin{enumarab}
\item Let $\delta$ be a cardinal. Let $\alpha$ be an ordinal, and let $G$ be the semigroup of finite transformations on $\alpha$.
Let$_{\alpha} \Fr_{\delta}$ be the absolutely free algebra on $\delta$
generators and of type $GPHA_{\alpha}$.  Let $\rho\in
{}^{\delta}\wp(\alpha)$. Let $L$ be a class having the same
similarity type as $GPHA_{\alpha}.$ Let
$$Cr_{\delta}^{(\rho)}L=\bigcap\{R: R\in Con_{\alpha}\Fr_{\delta},
{}_{\alpha}\Fr_{\delta}/R\in \mathbf{SP}L, {\mathsf
c}_k^{_{\alpha}\Fr_{\delta}}{\eta}/R=\eta/R \text { for each }$$
$$\eta<\delta \text
{ and each }k\in \alpha\smallsetminus \rho(\eta)\}$$ and
$$\Fr_{\delta}^{\rho}L={}_{\alpha}\Fr_{\beta}/Cr_{\delta}^{(\rho)}L.$$

The ordinal $\alpha$ does not figure out in $Cr_{\delta}^{(\rho)}L$
and $\Fr_{\delta}^{(\rho)}L$ though it is involved in their
definition. However, $\alpha$ will be clear from context so that no
confusion is likely to ensue.


\item Assume that $\delta$ is a cardinal, $L\subseteq GPHA_{\alpha}$, $\A\in L$,
$x=\langle x_{\eta}:\eta<\beta\rangle\in {}^{\delta}A$ and $\rho\in
{}^{\delta}\wp(\alpha)$. We say that the sequence $x$ $L$-freely
generates $\A$ under the dimension restricting function $\rho$, or
simply $x$ freely generates $\A$ under $\rho,$ if the following two
conditions hold:
\begin{enumroman}
\item $\A=\Sg^{\A}Rg(x)$ and $\Delta^{\A} x_{\eta}\subseteq \rho(\eta)$ for all $\eta<\delta$.
\item Whenever $\B\in L$, $y=\langle y_{\eta}, \eta<\delta\rangle\in
{}^{\delta}\B$ and $\Delta^{\B}y_{\eta}\subseteq \rho(\eta)$ for
every $\eta<\delta$, then there is a unique homomorphism from $\A$ to $\B$, such
that $h\circ x=y$.
\end{enumroman}
\end{enumarab}
\end{definition}
The second item says that dimension restricted free algebras has the universal property of free algebras with respect to algebras whose dimensions 
are also restricted. The following theorem can be easily distilled from the literature of cylindic algebra.

\begin{theorem}
Assume that $\delta$ is a cardinal, $L\subseteq GPHA_{\alpha}$,
$\A\in L$, $x=\langle x_{\eta}:\eta<\delta\rangle\in {}^{\delta}A$
and $\rho\in {}^{\delta}\wp(\alpha).$ Then the following hold:
\end{theorem}
\begin{enumroman}
\item $\Fr_{\delta}^{\rho}L\in GPHA_{\alpha}$
and $x=\langle \eta/Cr_{\delta}^{\rho}L: \eta<\delta \rangle$
$\mathbf{SP}L$- freely generates $\A$ under $\rho$.
\item In order that
$\A\cong \Fr_{\delta}^{\rho}L$ it is necessary and sufficient that
there exists a sequence $x\in {}^{\delta}A$ which $L$ freely
generates $\A$ under $\rho$.
\end{enumroman}
\begin{demo}{Proof} \cite{HMT1} theorems 2.5.35, 2.5.36, 2.5.37.
\end{demo}
Note that when $\rho(i)=\alpha$ for all $i$ then $\rho$ is not restricting the dimension, and we recover the notion of ordinary free algebras.
That is for such a $\rho$, we have $\Fr_{\beta}^{\rho}GPHA_{\alpha}\cong \Fr_{\beta}GPHA_{\alpha}.$

Now we formulate the analogue of  theorem \ref{main} for dimension restricted agebras, 
which adresses infinitely many cases, 
because we have infinitely many dimension restricted free algebras having the same number of generators.

\begin{theorem}\label{main2}
Let $G$ be the semigroup of finite transformations on an infinite set 
$\alpha$ and let $\delta$ be a cardinal $>0$. Let $\rho\in {}^{\delta}\wp(\alpha)$ be such that
$\alpha\sim \rho(i)$ is infinite for every 
$i\in \delta$. Let $\A$ be the free  $G$ algebra generated by $X$ restristed by $\rho$;
 that is $\A=\Fr_{\delta}^{\rho}GPHA_{\alpha},$
and suppose that $X=X_1\cup X_2$.
Let $(\Delta_0, \Gamma_0)$, $(\Theta_0, \Gamma_0^*)$ be two consistent theories in $\Sg^{\A}X_1$ and $\Sg^{\A}X_2,$ respectively.
Assume that $\Gamma_0\subseteq \Sg^{\A}(X_1\cap X_2)$ and $\Gamma_0\subseteq \Gamma_0^*$.
Assume, further, that  
$(\Delta_0\cap \Theta_0\cap \Sg^{\A}X_1\cap \Sg^{\A}X_2, \Gamma_0)$ is complete in $\Sg^{\A}X_1\cap \Sg^{\A}X_2$. 
Then there exist a Kripke system $\mathfrak{K}=(K,\leq \{X_k\}_{k\in K}\{V_k\}_{k\in K}),$ a homomorphism $\psi:\A\to \mathfrak{F}_K,$
$k_0\in K$, and $x\in V_{k_0}$,  such that for all $p\in \Delta_0\cup \Theta_0$ if $\psi(p)=(f_k)$, then $f_{k_0}(x)=1$
and for all $p\in \Gamma_0^*$ if $\psi(p)=(f_k)$, then $f_{k_0}(x)=0$.
\end{theorem}
\begin{demo}{Proof} We state the modifications in the above proof of theorem \ref{main}.
Form the sequence of minimal dilations $(\A_n:n\in \omega)$ built on the sequence $(K_n:n\in \omega)$, with $|K_n|=\beta$,
$\beta=|I\sim \alpha|=max(|A|, \alpha)$ with $I$ is a superset of $\alpha.$ 
If $i=((\Delta, \Gamma), (T,F))$ is a matched pair of saturated theories in $\Sg^{\A_n}X_1$ and $\Sg^{\A_n}X_2$, let $M_i=dim \A_n$, 
where $n$ is the least such number, so $n$ is unique to $i$. Define $K$ as in in the proof of theorem \ref{main}, that is, 
let $$K=\{((\Delta, \Gamma), (T,F)): \exists n\in \omega \text { such that } (\Delta, \Gamma), (T,F)$$
$$\text { is a a matched pair of saturated theories in }
\Sg^{\A_n}X_1, \Sg^{\A_n}X_2\}.$$
Let $${\mathfrak{K}}=(K, \leq, \{M_i\}, \{V_i\})_{i\in \mathfrak{K}},$$
where now $V_i={}^{\alpha}M_i^{(Id)}=\{s\in {}^{\alpha}M: |\{i\in \alpha: s_i\neq i\}|<\omega\},$ and the order $\leq $ is defined by:
If $i_1=((\Delta_1, \Gamma_1)), (T_1, F_1))$ and $i_2=((\Delta_2, \Gamma_2), (T_2,F_2))$ are in $\mathfrak{K}$, then
$$i_1\leq i_2\Longleftrightarrow  M_{i_1}\subseteq M_{i_2}, \Delta_1\subseteq \Delta_2, T_1\subseteq T_2.$$
This is a preorder on $K$.
Set $\psi_1: \Sg^{\A}X_1\to \mathfrak{F}_{\mathfrak K}$ by
$\psi_1(p)=(f_k)$ such that if $k=((\Delta, \Gamma), (T,F))\in \mathfrak{K}$ 
is a matched pair of saturated theories in $\Sg^{\A_n}X_1$ and $\Sg^{\A_n}X_2$,
and $M_k=dim \A_n$, then for $x\in V_k={}^{\alpha}M_k^{(Id)}$,
$$f_k(x)=1\Longleftrightarrow {\sf s}_{x\cup (Id_{M_k\sim \alpha)}}^{\A_n}p\in \Delta\cup T.$$
Define $\psi_2$ analogously. The rest of the proof is identical to the previous one.
\end{demo}
It is known that the condition $\Gamma\subseteq \Gamma^*$ cannot be omitted.
On the other hand, to prove our completeness theorem, we need the following weaker 
version of theorem \ref{main}, with a slight modification in the proof, which is still a step-by-step technique,
though, we do not `zig-zag'.  


\begin{lemma}\label{rep}
Let $\A\in GPHA_{\alpha}$. Let $(\Delta_0, \Gamma_0)$ be consistent. 
Suppose that $I$ is a set such that $\alpha\subseteq I$  and $|I\sim \alpha|=max (|A|,|\alpha|)$. 
\begin{enumarab}
\item Then there exists a dilation $\B\in GPHA_I$ of $\A$, and theory $T=(\Delta_{\omega}, \Gamma_{\omega})$, 
extending  $(\Delta_0, \Gamma_0)$, such that $T$  is consistent and saturated in $\B$. 
\item There exists $\mathfrak{K}=(K,\leq \{X_k\}_{k\in K}\{V_k\}_{k\in K}),$ a homomorphism $\psi:\A\to \mathfrak{F}_K,$
$k_0\in K$, and $x\in V_{k_0}$,  such that for all $p\in \Delta_0$ if $\psi(p)=(f_k)$, then $f_{k_0}(x)=1$ and for all $p\in \Gamma_0$
if $\psi(p)=(g_k)$, then $g_{k_0}(x)=0.$
\end{enumarab}
\end{lemma}
\begin{demo}{Proof} We deal only with the case when $G$ is strongly rich. The other cases can be dealt with in a similar manner 
by undergoing the obvious modifications, as indicated above. 
As opposed to theorem \ref{main}, we use theories rather than pairs of theories, since we are not dealing with two subalgebras simultaneously.
(i) follows from \ref{t2}. Now we prove (ii). The proof is a simpler version of the proof of \ref{main}.
Let $I$ be  a set such that 
$\beta=|I\sim \alpha|=max(|A|, |\alpha|).$
Let $(K_n:n\in \omega)$ be a family of pairwise disjoint sets such that $|K_n|=\beta.$
Define a sequence of algebras
$\A=\A_0\subseteq \A_1\subseteq \A_2\subseteq \A_2\ldots \subseteq \A_n\ldots$
such that
$\A_{n+1}$ is a minimal dilation of $\A_n$ and $dim(\A_{n+1})=\dim\A_n\cup K_n$.
We denote $dim(\A_n)$ by $I_n$ for $n\geq 1$. 
If $(\Delta, \Gamma)$ 
is saturated in $\A_n$ then the following analogues of Claims 1 and 2 in theorem \ref{main} hold:
For any $a,b\in \A_n$ if $a\rightarrow b\notin \Delta$, then there is a saturated  theory $(\Delta',\Gamma')$ in $\A_{n+1}$ 
such that $\Delta\subseteq \Delta'$ $a\in \Delta'$ and $b\notin \Delta'$.
If $(\Delta, \Gamma)$ is saturated in $\A_n$ then for all  $x\in \A_n$ and $j\in I_n$,  if ${\sf q}_{j}x\notin  \Delta,$ 
then there $(\Delta',\Gamma')$ 
of saturated theories in
$\A_{n+2}$, $u\in I_{n+2}$
such that $\Delta\subseteq \Delta'$, and ${\sf s}_j^u x\notin \Delta'$. Now let
$$K=\{(\Delta, \Gamma): \exists n\in \omega \text { such that } (\Delta,\Gamma) \text { is saturated in }\A_n.\}$$
If $i=(\Delta, \Gamma)$ is a saturated theory in $\A_n$, let $M_i=dim \A_n$, 
where $n$ is the least such number, so $n$ is unique to $i$.
If $i_1=(\Delta_1, \Gamma_1)$ and $i_2=(\Delta_2, \Gamma_2)$ are in $K$, then set
$$i_1\leq i_2\Longleftrightarrow  M_{i_1}\subseteq M_{i_2}, \Delta_1\subseteq \Delta_2. $$
This is a preorder on $K$; define the kripke system ${\mathfrak K}$ based on the set of worlds $K$ as before.
Set $\psi: \A\to \mathfrak{F}_{\mathfrak K}$ by
$\psi_1(p)=(f_k)$ such that if $k=(\Delta, \Gamma)\in \mathfrak{K}$ is  saturated in $\A_n$,
and $M_k=dim \A_n$, then for $x\in V_k=\bigcup_{p\in G_n}{}^{\alpha}M_k^{(p)}$,
$$f_k(x)=1\Longleftrightarrow {\sf s}_{x\cup (Id_{M_k\sim \alpha)}}^{\A_n}p\in \Delta.$$
Let  $k_0=(\Delta_{\omega}, \Gamma_{\omega})$ be defined as a complete saturated extension of $(\Delta_0, \Gamma_0)$
in $\A_1$, then $\psi,$ $k_0$ and $Id$ are as desired. 
The analogues of Claims 1 and 2 in theorem \ref{main} 
are used to show that $\psi$ so defined preserves implication and 
universal quantifiers.
\end{demo}

\section{Presence of diagonal elements}

All results, in Part 1, up to the previous theorem,  are proved in the absence of diagonal elements.
Now lets see how far we can go if we have diagonal elements. 
Considering diagonal elements, as we shall see, turn out to be problematic but not hopeless.

Our representation theorem has to respect diagonal elements, 
and this seems to be an impossible task with the presence of infinitary substitutions, 
unless we make a compromise that is, from our point of view, acceptable.
The interaction of substitutions based on infinitary transformations, 
together with the existence of diagonal elements tends to make matters `blow up'; indeed this even happens in the classical case,
when the class of (ordinary) set algebras ceases to be closed under ultraproducts \cite{S}. 
The natural thing to do is to avoid those infinitary substitutions at the start, while finding the interpolant possibly using such substitutions.
We shall also show that in some cases the interpolant has to use infinitary substitutions, even if the original implication uses only finite transformations.

So for an algebra $\A$, we let $\Rd\A$ denote its reduct when we discard infinitary substitutions. $\Rd\A$ satisfies 
cylindric algebra axioms.

\begin{theorem}\label{main3}
Let $\alpha$ be an infinite set. Let $G$ be a semigroup on $\alpha$ containing at least one infinitary transformation. 
Let $\A\in GPHAE_{\alpha}$ be the free $G$ algebra generated by $X$, and suppose that $X=X_1\cup X_2$.
Let $(\Delta_0, \Gamma_0)$, $(\Theta_0, \Gamma_0^*)$ be two consistent theories in $\Sg^{\Rd\A}X_1$ and $\Sg^{\Rd\A}X_2,$ respectively.
Assume that $\Gamma_0\subseteq \Sg^{\A}(X_1\cap X_2)$ and $\Gamma_0\subseteq \Gamma_0^*$.
Assume, further, that  
$(\Delta_0\cap \Theta_0\cap \Sg^{\A}X_1\cap \Sg^{\A}X_2, \Gamma_0)$ is complete in $\Sg^{\Rd\A}X_1\cap \Sg^{\Rd\A}X_2$. 
Then there exist $\mathfrak{K}=(K,\leq \{X_k\}_{k\in K}\{V_k\}_{k\in K}),$ a homomorphism $\psi:\A\to \mathfrak{F}_K,$
$k_0\in K$, and $x\in V_{k_0}$,  such that for all $p\in \Delta_0\cup \Theta_0$ if $\psi(p)=(f_k)$, then $f_{k_0}(x)=1$
and for all $p\in \Gamma_0^*$ if $\psi(p)=(f_k)$, then $f_{k_0}(x)=0$.
\end{theorem}
\begin{demo}{Proof}
The first half of the proof is almost identical to that of  lemma \ref{main}. We highlight the main steps, 
for the convenience of the reader, except that we only deal with the case
when $G$ is strongly rich.
Assume, as usual, that $\alpha$, $G$, $\A$ and $X_1$, $X_2$, and everything else in the hypothesis are given.
Let $I$ be  a set such that  $\beta=|I\sim \alpha|=max(|A|, |\alpha|).$
Let $(K_n:n\in \omega)$ be a family of pairwise disjoint sets such that $|K_n|=\beta.$
Define a sequence of algebras
$\A=\A_0\subseteq \A_1\subseteq \A_2\subseteq \A_2\ldots \subseteq \A_n\ldots$
such that
$\A_{n+1}$ is a minimal dilation of $\A_n$ and $dim(\A_{n+1})=\dim\A_n\cup K_n$.We denote $dim(\A_n)$ by $I_n$ for $n\geq 1$. 
The proofs of Claims 1 and 2 in the proof of \ref{main} are the same.

Now we prove the theorem when $G$ is a strongly rich semigroup.
Let $$K=\{((\Delta, \Gamma), (T,F)): \exists n\in \omega \text { such that } (\Delta, \Gamma), (T,F)$$
$$\text { is a a matched pair of saturated theories in }
\Sg^{\Rd\A_n}X_1, \Sg^{\Rd\A_n}X_2\}.$$
We have $((\Delta_0, \Gamma_0)$, $(\Theta_0, \Gamma_0^*))$ is a matched pair but the theories are not saturated. But by lemma \ref{t3}
there are $T_1=(\Delta_{\omega}, \Gamma_{\omega})$, 
$T_2=(\Theta_{\omega}, \Gamma_{\omega}^*)$ extending 
$(\Delta_0, \Gamma_0)$, $(\Theta_0, \Gamma_0^*)$, such that $T_1$ and $T_2$ are saturated in $\Sg^{\Rd\A_1}X_1$ and $\Sg^{\Rd\A_1}X_2,$ 
respectively. Let $k_0=((\Delta_{\omega}, \Gamma_{\omega}), (\Theta_{\omega}, \Gamma_{\omega}^*)).$ Then $k_0\in K,$
and   $k_0$ will be the desired world and $x$ will be specified later; in fact $x$ will be the identity map on some specified 
domain.

If $i=((\Delta, \Gamma), (T,F))$ is a matched pair of saturated theories in $\Sg^{\Rd\A_n}X_1$ and $\Sg^{\Rd\A_n}X_2$, let $M_i=dim \A_n$, 
where $n$ is the least such number, so $n$ is unique to $i$.
Let $${\bf K}=(K, \leq, \{M_i\}, \{V_i\})_{i\in \mathfrak{K}},$$
where $V_i=\bigcup_{p\in G_n, p\text { a finitary transformation }}{}^{\alpha}M_i^{(p)}$ 
(here we are considering only substitutions that move only finitely many points), 
and $G_n$ 
is the strongly rich semigroup determining the similarity type of $\A_n$, with $n$ 
the least number such $i$ is a saturated matched pair in $\A_n$, and $\leq $ is defined as follows: 
If $i_1=((\Delta_1, \Gamma_1)), (T_1, F_1))$ and $i_2=((\Delta_2, \Gamma_2), (T_2,F_2))$ are in $\bold K$, then set
$$i_1\leq i_2\Longleftrightarrow  M_{i_1}\subseteq M_{i_2}, \Delta_1\subseteq \Delta_2, T_1\subseteq T_2.$$ 
We are not yet there, to preserve diagonal elements we have to factor out $\bold K$ 
by an infinite family equivalence relations, each defined on the dimension of $\A_n$, for some $n$, which will actually turn out to be 
a congruence in an exact sense. 
As usual, using freeness of $\A$, we will  define two maps on $\A_1=\Sg^{\Rd\A}X_1$ and $\A_2=\Sg^{\Rd\A}X_2$, respectively;
then those will be pasted 
to give the required single homomorphism.

Let $i=((\Delta, \Gamma), (T,F))$ be a matched pair of saturated theories in $\Sg^{\Rd\A_n}X_1$ and $\Sg^{\Rd\A_n}X_2$, let $M_i=dim \A_n$, 
where $n$ is the least such number, so $n$ is unique to $i$.
For $k,l\in dim\A_n=I_n$, set $k\sim_i l$ iff ${\sf d}_{kl}^{\A_n}\in \Delta\cup T$. This is well defined since $\Delta\cup T\subseteq \A_n$.
We omit the superscript $\A_n$.
These are infinitely many relations, one for each $i$, defined on $I_n$, with $n$ depending uniquely on $i$, 
we denote them uniformly by $\sim$ to 
avoid complicated unnecessary notation.
We hope that no confusion is likely to ensue. We claim that $\sim$ is an equivalence relation on $I_n.$
Indeed,  $\sim$ is reflexive because ${\sf d}_{ii}=1$ and symmetric 
because ${\sf d}_{ij}={\sf d}_{ji};$
finally $E$ is transitive because for  $k,l,u<\alpha$, with $l\notin \{k,u\}$, 
we have 
$${\sf d}_{kl}\cdot {\sf d}_{lu}\leq {\sf c}_l({\sf d}_{kl}\cdot {\sf d}_{lu})={\sf d}_{ku},$$
and we can assume that $T\cup \Delta$ is closed upwards.
For $\sigma,\tau \in V_k,$ define $\sigma\sim \tau$ iff $\sigma(i)\sim \tau(i)$ for all $i\in \alpha$.
Then clearly $\sigma$ is an equivalence relation on $V_k$. 

Let $W_k=V_k/\sim$, and $\mathfrak{K}=(K, \leq, M_k, W_k)_{k\in K}$, with $\leq$ defined on $K$ as above.
We write $h=[x]$ for $x\in V_k$ if $x(i)/\sim =h(i)$ for all $i\in \alpha$; of course $X$ may not be unique, but this will not matter.
Let $\F_{\mathfrak K}$ be the set algebra based on the new Kripke system ${\mathfrak K}$ obtained by factoring out $\bold K$.

Set $\psi_1: \Sg^{\Rd\A}X_1\to \mathfrak{F}_{\mathfrak K}$ by
$\psi_1(p)=(f_k)$ such that if $k=((\Delta, \Gamma), (T,F))\in K$ 
is a matched pair of saturated theories in $\Sg^{\Rd\A_n}X_1$ and $\Sg^{\Rd\A_n}X_2$,
and $M_k=dim \A_n$, with $n$ unique to $k$, then for $x\in W_k$
$$f_k([x])=1\Longleftrightarrow {\sf s}_{x\cup (Id_{M_k\sim \alpha)}}^{\A_n}p\in \Delta\cup T,$$
with $x\in V_k$ and $[x]\in W_k$ is define as above. 

To avoid cumbersome notation, we 
write ${\sf s}_{x}^{\A_n}p$, or even simply ${\sf s}_xp,$ for 
${\sf s}_{x\cup (Id_{M_k\sim \alpha)}}^{\A_n}p$.  No ambiguity should arise because the dimension $n$ will be clear from context.

We need to check that $\psi_1$ is well defined. 
It suffices to show that if $\sigma, \tau\in V_k$ if $\sigma \sim \tau$ and $p\in \A_n$,  
with $n$ unique to $k$, 
then $${\sf s}_{\tau}p\in \Delta\cup T\text { iff } {\sf s}_{\sigma}p\in \Delta\cup T.$$

This can be proved by induction on the cardinality of 
$J=\{i\in I_n: \sigma i\neq \tau i\}$, which is finite since we are only taking finite substitutions.
If $J$ is empty, the result is obvious. 
Otherwise assume that $k\in J$. We recall the following piece of notation.
For $\eta\in V_k$ and $k,l<\alpha$, write  
$\eta(k\mapsto l)$ for the $\eta'\in V$ that is the same as $\eta$ except
that $\eta'(k)=l.$ 
Now take any 
$$\lambda\in \{\eta\in I_n: \sigma^{-1}\{\eta\}= \tau^{-1}\{\eta\}=\{\eta\}\}\smallsetminus \Delta x.$$
This $\lambda$ exists, because $\sigma$ and $\tau$ are finite transformations and $\A_n$ is a dilation with enough spare dimensions.
We have by cylindric axioms (a)
$${\sf s}_{\sigma}x={\sf s}_{\sigma k}^{\lambda}{\sf s}_{\sigma (k\mapsto \lambda)}p.$$
We also have (b)
$${\sf s}_{\tau k}^{\lambda}({\sf d}_{\lambda, \sigma k}\land {\sf s}_{\sigma} p)
={\sf d}_{\tau k, \sigma k} {\sf s}_{\sigma} p,$$
and (c)
$${\sf s}_{\tau k}^{\lambda}({\sf d}_{\lambda, \sigma k}\land {\sf s}_{\sigma(k\mapsto \lambda)}p)$$
$$= {\sf d}_{\tau k,  \sigma k}\land {\sf s}_{\sigma(k\mapsto \tau k)}p.$$
and (d)
$${\sf d}_{\lambda, \sigma k}\land {\sf s}_{\sigma k}^{\lambda}{\sf s}_{{\sigma}(k\mapsto \lambda)}p=
{\sf d}_{\lambda, \sigma k}\land {\sf s}_{{\sigma}(k\mapsto \lambda)}p$$
Then by (b), (a), (d) and (c), we get,
$${\sf d}_{\tau k, \sigma k}\land {\sf s}_{\sigma} p= 
{\sf s}_{\tau k}^{\lambda}({\sf d}_{\lambda,\sigma k}\cdot {\sf s}_{\sigma}p)$$
$$={\sf s}_{\tau k}^{\lambda}({\sf d}_{\lambda, \sigma k}\land {\sf s}_{\sigma k}^{\lambda}
{\sf s}_{{\sigma}(k\mapsto \lambda)}p)$$
$$={\sf s}_{\tau k}^{\lambda}({\sf d}_{\lambda, \sigma k}\land {\sf s}_{{\sigma}(k\mapsto \lambda)}p)$$
$$= {\sf d}_{\tau k,  \sigma k}\land {\sf s}_{\sigma(k\mapsto \tau k)}p.$$
The conclusion follows from the induction hypothesis.
Now $\psi_1$ respects all quasipolyadic equality operations, that is finite substitutions (with the proof as before; 
recall that we only have finite substitutions since we are considering 
$\Sg^{\Rd\A}X_1$)  except possibly for diagonal elements. 
We check those:

Recall that for a concrete Kripke frame $\F_{\bold W}$ based on ${\bold W}=(W,\leq ,V_k, W_k),$ we have
the concrete diagonal element ${\sf d}_{ij}$ is given by the tuple $(g_k: k\in K)$ such that for $y\in V_k$, $g_k(y)=1$ iff $y(i)=y(j)$.

Now for the abstract diagonal element in $\A$, we have $\psi_1({\sf d}_{ij})=(f_k:k\in K)$, such that if $k=((\Delta, \Gamma), (T,F))$ 
is a matched pair of saturated theories in $\Sg^{\Rd\A_n}X_1$, $\Sg^{\Rd\A_n}X_2$, with $n$ unique to $i$, 
we have $f_k([x])=1$ iff ${\sf s}_{x}{\sf d}_{ij}\in \Delta \cup T$ (this is well defined $\Delta\cup T\subseteq \A_n).$
 
But the latter is equivalent to ${\sf d}_{x(i), x(j)}\in \Delta\cup T$, which in turn is equivalent to $x(i)\sim x(j)$, that is 
$[x](i)=[x](j),$ and so  $(f_k)\in {\sf d}_{ij}^{\F_{\mathfrak K}}$.  
The reverse implication is the same.

We can safely assume that $X_1\cup X_2=X$ generates $\A$.
Let $\psi=\psi_1\cup \psi_2\upharpoonright X$. Then $\psi$ is a function since, by definition, $\psi_1$ and $\psi_2$ 
agree on $X_1\cap X_2$. Now by freeness $\psi$ extends to a homomorphism, 
which we denote also by $\psi$ from $\A$ into $\F_{\mathfrak K}$.
And we are done, as usual, by $\psi$, $k_0$ and $Id\in V_{k_0}$.
\end{demo}

Theorem \ref{main2}, generalizes as is, to the expanded structures by diagonal elements. That is to say, we have:

\begin{theorem}\label{main4}
Let $G$ be the semigroup of finite transformations on an infinite set 
$\alpha$ and let $\delta$ be a cardinal $>0$. Let $\rho\in {}^{\delta}\wp(\alpha)$ be such that
$\alpha\sim \rho(i)$ is infinite for every 
$i\in \delta$. Let $\A$ be the free  $G$ algebra with equality generated by $X$ restristed by $\rho$;
 that is $\A=\Fr_{\delta}^{\rho}GPHAE_{\alpha},$
and suppose that $X=X_1\cup X_2$.
Let $(\Delta_0, \Gamma_0)$, $(\Theta_0, \Gamma_0^*)$ be two consistent theories in $\Sg^{\A}X_1$ and $\Sg^{\A}X_2,$ respectively.
Assume that $\Gamma_0\subseteq \Sg^{\A}(X_1\cap X_2)$ and $\Gamma_0\subseteq \Gamma_0^*$.
Assume, further, that  
$(\Delta_0\cap \Theta_0\cap \Sg^{\A}X_1\cap \Sg^{\A}X_2, \Gamma_0)$ is complete in $\Sg^{\A}X_1\cap \Sg^{\A}X_2$. 
Then there exist a Kripke system $\mathfrak{K}=(K,\leq \{X_k\}_{k\in K}\{V_k\}_{k\in K}),$ a homomorphism $\psi:\A\to \mathfrak{F}_K,$
$k_0\in K$, and $x\in V_{k_0}$,  such that for all $p\in \Delta_0\cup \Theta_0$ if $\psi(p)=(f_k)$, then $f_{k_0}(x)=1$
and for all $p\in \Gamma_0^*$ if $\psi(p)=(f_k)$, then $f_{k_0}(x)=0$.
\end{theorem}
\begin{demo}{Proof} $\Rd\A$ is just $\A$.
\end{demo}  
\section{Results in logical form}

We start by describing our necessary syntactical and semantical notions to follow.
Informally a language is a triple $(V, P, G)$ where $V$ is a set providing an infinite supply of variables, 
$P$ is a another set of predicates disjoint from $V,$
and $G$ is a semigroup of transformations on $V$. 
There is no restriction on the arity of $p\in P$; sometimes referred to as the rank of $p$,  that is the arity may be infinite. 
Formulas are defined recursively the usual way. Atomic formulas are of the form $p\bar{v}$, the length of $\bar{v}$ is equal to the arity of $p$.
If $\phi, \psi$ are formulas and $v\in V,$
then $\phi\lor\psi$, $\phi\land \psi$, $\phi\to \psi$, $\exists v\phi,$ $\forall v\phi$ are formulas. For each $\tau\in G$, 
${\sf S}({\tau})$ is a unary  operation on formulas, that is, for any formula $\phi$, ${\sf S}{(\tau)}\phi$ is another formula,  
reflecting the metalogical operation of simultaneous substitution
of variables (determined by $\tau$) for variables, such that the substitution is free. 
Notice that although we allow infinitary predicates, quantifications are defined only on finitely many variables, 
that is the scope of quantifiers is finite.  

We will also deal with the case when we have equality; for this purpose we add a newlogical symbol $=$ and we view it, as usual, as a binary relation.

We recall some basic semantical notions for intuitionistic logic but adpated to the presence of atomic formulas possibly having infinite length. 
An intuitionistic or Kripke frame is a triple $\bold W=(W, R, \{D_w\}_{w\in W})$ where $W$ is a non-empty set called worlds, 
preordered by $R$ and $D_w$ is a non-empty subset of $D$
called the domain of $w$ for any $w$, and the monotoncity condition of domains is satisfied:
$$(\forall w,w'\in W)[wRw'\implies  D_w\subseteq D_{w'}.]$$
On the other hand, an  intuitionistic or Kripke model is a 
quadruple $\bold M=(W, R, \{D_{w}\}_{w\in W} \models),$ where $(W, R, \{D_w\}_{w\in W})$ 
is an intuitionistic frame, $\models$ is a tenary relation between worlds, formulas, and assignments (maps from $V$ to $D$).
We write $x\models \phi[s]$ if  $(x, \phi ,s)\in \models$. This tenary relation $\models$ satisfies for
any predicate $p$, any $s\in {}^{V}D$, any formulas $\phi$, $\psi$ and any $x\in W$ the following:
$$ \text { It is not the case that }x\models \bot,$$
$$(\forall y\in W)(x\models p[s]\land xRy\implies y\models p[s]),$$
$$x\models (\phi\land \psi)[s]\Longleftrightarrow x\models \phi[s] \text { and } x\models \psi[s],$$
$$x\models (\phi\lor \psi)[s]\Longleftrightarrow  x\models \phi[s]\text { or }\phi\models \psi[s],$$
$$x\models (\phi\to \psi)[s]\Longleftrightarrow \forall y(xRy\implies(y\models \phi[s]\implies y\models \psi[s])).$$
For $s$ a function $s^k_a$ is the function defined by $s^k_a(i)=s(i)$ when $i\neq k$ and $s^k_a(k)=a$. Continuing the definition:
$$x\models \forall v\phi[s]\Longleftrightarrow( (\forall y)(xRy\implies (\forall a\in D_y)y\models \phi[s^v_a]))),$$
$$x\models \exists v\phi[s]\Longleftrightarrow (\exists a\in D_x)(x\models \phi[s^v_a])),$$
$$x\models {\sf S}({\tau})\phi[s]\Longleftrightarrow x\models \phi[\tau\circ s].$$
Evidently the model is completely determined by the frame $(W, R, \{D_{w}\}_{w\in W})$ and by $\models$ on atomic formulas.
That is for each for each $p\in P$ and each world $x$ and $s\in{}^VD_x$, $p$ determines a possibly infinitary relation $p_x\subseteq {}^VD_x$,
and we stipulate that $x\models p[s]$ if $s\in p_x$. 
If we have equality, then for the world $x$ and $s\in {}^VD_x$, we add the clause $x\models v_1=v_2$ if $s(v_1)=s(v_2)$.

We now define a calculas (in a usual underlying set theory $ZFC$, say) 
that we prove to be complete with respect to Kripke semantics; 
this will follow from our stronger proven result that such logics 
enjoy the interpolation property. We first deal with the equality free case. In such a case, our calculas is inspired by that of Keisler \cite{K}.

Let $V$ and $P$ be disjoint sets of symbols, such that $V$ is infinite, $\rho$ a function with domain is $P$ and whose range is $Set$ (the class of all sets)
or $Ord$
(the class of all ordinals).
\footnote{Strictly speaking, in $ZFC$ we cannot talk about classes, but classes can be stimulated rigorously with formulas; 
in our context we chose
not to be pedantic about it. Alternatively,  we could have  replaced $Set$ ($Ord$) by a set of sets (ordinals), 
but the notation $Set$ ($Ord$) is more succint and ecconomic.} 

Let $G\subseteq {}^VV$. We define a logic $\mathfrak{L}_G$ in the following way.
The symbols of $\mathfrak{L}_G$ consists of: 
\begin{enumarab}
\item The falsity symbol $\bot$.
\item the disjunction $\lor$, conjunction $\land$, and the implication symbol $\to$.
\item universal quantification symbol $\forall$.
\item existential quantification symbol $\exists$.
\item the individual variables $v\in V$ and predicates $p\in P.$
\end{enumarab}
We assume that $\bot, \lor, \land,\to, \forall, \exists$ are not members of $V$ nor $P$.
An atomic formula is an ordered pair $(p,x)$ 
where $p\in P$ and $x\in {}^{\rho(p)}V.$ We call $\rho(p)$ the rank of $p$. Formulas are defined the usual way by recursion; in this respect 
we regard $(\phi\to \phi)$ as an ordered triple and so are formulas involving other connectives including $\exists v\phi$ and $\forall v \phi.$
(In the former formula, the brackets are not syntactic brackets because we do no have brackets in our language.)
The set $V_f(\phi)$ of free variables and the set $V_b(\phi)$ of bound variables in a formula $\phi$ are defined recursively the usual way.
That is:
 \begin{enumarab}
\item If $\phi$ is an atomic formula $(p,x)$, then $V_f(\phi)$ is the range of $x$.
\item if $\phi=\bot$, then $V_f(\bot)=\emptyset.$
\item If $\phi$ is $(\psi\lor \theta)$ or $(\psi\land \theta)$ or $(\psi\to \theta)$, then $V_f(\phi)=V_f(\psi)\cup V_f(\theta).$
\item If $\phi=(\forall v \psi)$ or $(\exists v\psi)$,  then $V_f(\psi)=V_f(\phi)\sim \{v\}$.
Now for the bound variables $V_b(\phi)$:
\item If $\phi$ is an atomic formula $(p,x)$, then $V_b(\phi)=\emptyset.$ 
\item if $\phi=\bot$, then $V_b(\bot)=\emptyset.$
\item If $\phi$ is $(\psi\lor \theta)$ or $(\psi\land \theta)$ or $(\psi\to \theta),$ then $V_b(\phi)=V_b(\psi)\cup V_b(\theta).$
\item If $\phi=(\forall v \psi)$, then $V_b(\psi)=V_b(\phi)\cup \{v\}.$
\item If $\phi=(\exists v \psi)$, then $V_b(\psi)=V_b(\phi)\cup \{v\}.$
\end{enumarab}
Note that the variables occurring in a formula $\phi$, denoted by $V(\phi)$  is equal to $V_f(\phi)\cup V_b(\phi)$ which could well be infinite.
For $\tau\in G$ and $\phi$ a formula, ${\sf S}(\tau)\phi$ (the result of substituting each variable $v$ in $\phi$ by $\tau(v)$)
is defined recursively and so is ${\sf S}_f(\tau)\phi$ (the result of substituting each free variable $v$ by $\tau(v)$).
\begin{enumarab}
\item If $\phi$ is atomic formula $(p,x),$ then ${\sf S}(\tau)\phi=(p,\tau\circ x).$
\item if $\phi=\bot,$ then ${\sf S}(\tau)\bot=\bot$ 
\item If $\phi$ is $(\psi\lor \theta),$ then  ${\sf S}(\tau)\phi=({\sf S}(\tau)\psi\lor {\sf S}(\tau)\theta).$ The same for other propositional connectives.
\item If $\phi=(\forall v\phi),$ then ${\sf S}(\tau)\phi=(\forall\tau(v){\sf S}(\tau)\phi).$
\item If $\phi=(\exists  v\phi),$ then ${\sf S}(\tau)\phi=(\exists\tau(v){\sf S}(\tau)\phi).$

\end{enumarab}
To deal with free substitutions, that is when the resulted substituted variables remain free,
 we introduce a piece of notation that proves helpful. For any function $f\in {}^XY$ and any set $Z$, we let 
$$f|Z=\{(x, f(x)): x\in X\cap Z\}\cup \{(z,z)|z\in Z\sim X\}.$$
Then $f|Z$ always has domain $Z$ and $0|Z$ is the identity function on $Z$.

If $\tau\in \bigcup\{^WV: W\subseteq V\}$, and $\phi$ is a formula, let ${\sf S}(\tau)\phi={\sf S}(\tau|V)\phi$ and ${\sf S}_f(\tau)\phi=S_f(\tau|V)\phi$.

For free subtitution  the first three  clauses are the same, but if
$\phi=(\forall v \psi)$, then ${\sf S}_f({\tau})\phi=(\forall v {\sf S}_f(\sigma)\psi)$ 
and if $\phi=(\exists v \psi),$ then ${\sf S}_f({\tau})\phi=(\exists v {\sf S}_f(\sigma)\psi)$ 
where $\sigma=\tau\upharpoonright (V\sim \{v\})\upharpoonright V$.
Now we specify the axioms and the rules of inference.

The axioms are:
\begin{enumarab}
\item Axioms for propositional intuitionistic logic (formulated in our syntax).
\item $((\forall v(\phi\to \psi)\to (\phi\to \forall v \psi)))$ where $v\in (V\sim V_{f} \phi)).$
\item $((\forall v(\phi\to \psi)\to (\exists v\phi\to \psi)))$ where $v\in (V\sim V_{f} \phi)).$

\item $(\forall v\phi\to {\sf S}_f(\tau)\phi)$, when  $\tau(v) \notin (V\sim V_b(\phi)).$
\item $({\sf S}_f(\tau)\phi\to (\exists v \phi))$, when  $\tau(v) \notin (V\sim V_b(\phi)).$

\end{enumarab}
The rules of inference are:
\begin{enumarab}
\item Form $\phi$, $(\phi\to \psi)$ infer $\psi.$ (Modus ponens.)
\item From $\phi$ infer $(\forall v\phi).$ (Rule of generalization.)
\item From ${\sf S}_f(\tau)\phi$ infer $\phi$ whenever $\tau\in {}^{V_f(\phi)} (V\sim V_b(\phi))$ and $\tau$ is one to one. (Free substitution.)
\item From $\phi$ infer ${\sf S}(\tau)\phi$ whenever $\tau\in {}^{V(\phi)}V$ is one to one (Substitution).
\end{enumarab}
Now if we have $=$ as a primitive symbol, we add the following axioms (in this case no more rules of inference are needed):
\begin{enumarab}
\item $v=v$
\item $v=w\to w=v$
\item If $\phi$ is a formula and $\tau, \sigma$ are substitutions that agree on the indices of the free variables occuring in $\phi$, then
${\sf S}_f(\tau)\phi={\sf S}_f(\sigma)\phi.$
\end{enumarab}

We write ${\mathfrak L}_G$ for logics without equality, and we write ${\mathfrak L}_G^{=}$ for those with equality, when $G$ is specified
in advance.

Proofs are defined the usual way. For a set of formulas $\Gamma\cup \{\phi\}$, we write $\Gamma\vdash \phi$, 
if there is a proof
of $\phi$ from $\Gamma$.

To formulate the main results of this paper, we need some more basic definitions.
Let $\bold M=(W,R, \{D_{w}\}_{w\in W}, \models)$ be a Kripke model over $D$ and let $s\in {}^VD,$ where $V$ is the set of all variables.
A formula $\phi$ is satsifiable at $w$ under $s$ if $w\models \phi[s]$. The formula $\phi$ is satisfiable in $\bold M$ if there a $w\in W$ and $s\in {}^VD$ 
such that $w\models \phi[s].$ For a set of formulas $\Gamma$, we write $w\models \Gamma[s]$ if
$w\models \phi[s]$ for every $\phi\in \Gamma$. The set of formulas 
$\Gamma$ is satisfiable in $\bold M$ if there is a $w\in W$ and $s\in {}^VD$ 
such that $w\models \Gamma[s]$. The formula $\phi$ is valid in $\bold M$ under $s$ if $w\models \phi[s]$ 
for all $w\in W$ and $s\in {}^{V}D_w$; $\phi$ is valid in 
$\bold M$ if it is valid for any $s\in {}^VD$. A formula $\phi$ is valid in a frame $(W, R, \{D_w\}_{w\in W})$
if it is valid in every model based on $W$ after specifying the semantical consequence relation $\models$. 
A set of formulas $\Gamma$ is consistent if no 
contradiction is derivable from $\Gamma$ relative to the proof system defined above, that is, it is not the case that $T\vdash \bot.$

The custom in intuitionistic logic is to deal with pairs of theories, 
the first component dealing with a set of formulas that are 'true', and the second deals with a set formulas that are 
`false', in the intended interpretaton.
This is natural, since we do not have negation. So in fact, our algebraic counterpart proved in section 3,   
is in fact more general than the completeness theorem stated below; the latter 
follows from the special case when the second component of pairs is the theory 
$\{\bot\}$.

The following theorems hold for logics without equality. 
In the presence of infinitary substitutions, we obtain a weaker result for logics with equality. The set $V$ denoting the set of variables
in the next theorems is always infinite (which means that we will deal only with infinite dimensional algebras), 
however, $P$ (specifying the number of atomic formulas) could well be finite.

\begin{theorem}\label{com}\begin{enumroman}
\item  Let $V$ and $P$ be countable disjoint sets with $|V|\geq \omega$.  
When $G$ is a rich semigroup, then $\mathfrak{L}_G$ is strongly complete,
that is if $\Gamma$ is a consistent set of formulas, 
then it is satisfiable at a world of some model based on a Kripke frame. 
\item For arbitrary (disjoint) sets $V$ and $P$ with $|V|\geq \omega$, 
when $G$ is the semigroup of finite transformations, and $\rho\in {}^{V}Ord$ is
such that $V\sim \rho(p)$ is infinite for every $p\in P$, or $G={}^VV$ without any restrictions, 
then $\mathfrak{L}_G$ is strongly complete.
\end{enumroman}
\end{theorem}
\begin{demo}{Proof} cf. Theorem \ref{complete}, item  (1).
\end{demo}
We say that a logic $\mathfrak{L}$ has the Craig interpolation property if whenever 
$\models \phi\to \psi$ then there is a formula containing only symbols occurring in both $\phi$ and $\psi,$ $\theta$ say, such that
$\models \phi\to \theta$ and $\models \theta\to \psi.$ (By the above completeness theorem, we can replace $\models$ by $\vdash$.)
\begin{theorem}\label{interpolation} Let $\mathfrak{L}_G$ be as in the previous theorem, except that $G$ is assumed to be strongly rich.
Then $\mathfrak{L}_G$ has the interpolation property
\end{theorem} 
\begin{demo}{Proof} cf. Theorem \ref{complete}, item (2).
\end{demo}

In the case we have equality then we can prove a slightly weaker result when we have infinite substitutions. We say that the substitution 
operation ${\sf S}_{\tau}$ is finitary, if $\tau$ moves only finitely many points, otherwise,  it is called infinitary. Now we have: 

\begin{theorem}\label{interpolationeq}
\begin{enumarab}
\item  For arbitrary (disjoint) sets $V$ and $P$, with $|V|\geq \omega$, when $G$ is the semigroup of finite transformations, and $\rho\in {}^{V}Ord$ 
is such that $V\sim \rho(p)$ is infinite for every $p\in P$, then $\mathfrak{L}_G^{=}$ is strongly complete, 
and has the interpolation property.

\item When $G$ is  rich or $G={}^VV,$
then ${\mathfrak L}_G^{=}$ is weakly complete, that is, if a formula is valid in all Kripke models, then it is provable.  

\item When $G$ is strongly rich or $G={}^VV$, 
the logic $\mathfrak{L}_G^{=}$ has the following weak interpolation property. 
If $\phi$ and $\psi$ are formulas such that only finitary substitutions
were involved in their built up from atomic formulas, and   
$\models \phi\to \psi,$ then there is a formula containing only (atomic formulas, and possibly equality) 
occurring in both $\phi$ and $\psi,$ $\theta$ say, such that
$\models \phi\to \theta$ and $\models \theta\to \psi$, and $\theta$ may involve infinitary 
substitutions during its formation from atomic formulas. (By weak completenes $\models$ can be replaced by $\vdash$.)
\end{enumarab}
\end{theorem} 
\begin{demo}{Proof} From theorems \ref{main3}, \ref{main4}.
\end{demo}
\begin{theorem}\label{negative}
For arbitrary (disjoint) sets $V$ and $P$, with $|V|\geq \omega$, when $G$ is the semigroup of finite transformations, and $\rho\in {}^{V}Ord$ 
is such that $\rho(p)=V$ for every $p\in P$, then both $\mathfrak{L}_G$ and $\mathfrak{L}_G^{=}$ are essentially  incomplete, 
and fail to enjoy  the interpolation property.
\end{theorem}
\begin{demo}{Proof} This is proved in [Sayed].
\end{demo}

In intuitionistic ordinary predicate logic, interpolation theorems  proved to hold for logic without equality, 
remain to hold when we add equality. 
This is reflected by item (1) in theorem \ref{interpolationeq}. 
Indeed, in this case our logics are very close to ordinary ones. 
The sole difference is that atomic formulas could have infinite arity, 
but like the ordinary case, infinitely many variables lie outside (atomic) formulas.
The  next item in theorem \ref{interpolationeq} 
shows that the situation is not as smooth nor as evident as the ordinary classical case.

The presence of infinitary substitutions seems to make a drastic two-fold change. 
In the absence of diagonal elements, it turns negative results to positive ones,
but it in the presence of diagonal elements  the  positive results obtained are weaker. 

Indeed, we do not know whether strong completeness or usual interpolation holds for such logics, but it seems unlikely that they do.
We know that there will always be cases when infinitary substitutions are needed in the interpolant.

Our logics manifest themselves as essentially infinitary in at least two facets. 
One is that the atomic formulas can have infinite arity and the other is that (infinitary) substitutions, when available, can move infinitely many 
points. But they also have a finitary flavour since quantification is taken only on finitely many variables.
The classical counterpart of such logics has been studied frequently in algebraic logic, and they occur in the literature
under the name of finitary logics of 
infinitary relations, or typless logics \cite{HMT2}, though positive interpolation theorems for such logics are only rarely investigated [IGPL], 
for this area is 
dominated by negative results [references].

It is well known that first order predicate intuitionistic logic has the following two properties:

(*) Each proof involves finitely many formulas.

(*) A set of formulas is consistent if and only if it is satisfiable. 

In most cases, such as those logics which have infinitary propositional connectives, 
it is known to be impossible to define a notion of proof in such a way that both 
(*) and (**) are satisfied. We are thus confronted with the special situation that the logic ${\mathfrak L}_G$ behave like ordinary first order 
intuitionistic logic.  In passing, we note that (infinitary) generalizations of the classical Lowenheim-Skolem Theorem and of the 
Compactness Theorem for ${\mathfrak L}_G$ without equality
follows immediately from theorem \ref{com}.

Now we are ready to prove theorem \ref{interpolation}.
\begin{corollary}\label{complete} Let $G$ be a semigroup as in \ref{complete} and \ref{interpolation} 
\begin{enumarab}
\item $\mathfrak L_G$ is strongly complete
\item $\mathfrak L_G$ has the interpolation property
\end{enumarab}\end{corollary} 
\begin{demo}{Proof}

\begin{enumarab}
\item  We prove the theorem when $G$ is a strongly rich semigroup on $\alpha$, $\alpha$ a countable ordinal specifying the 
the number of variables in ${\mathfrak L}_G$.
Let $\{R_i:i\in \omega\}$ be the number of relation symbols available in our language each of arity $\alpha$.
We show that every consistent set of formulas $T$ is satisfiable at some world in a Kripke model.
Assume that $T$ is consistent.  Let $\A=\Fm/\equiv$ 
and let $\Gamma=\{\phi/\equiv: \phi\in T\}.$  Then $\Gamma$ generates a filter $F$. 
Then $\A\in GPHA_{\alpha}$ and $(F,\{0\})$ is consistent. 
By the above proof, it is satisfiable, that is 
there exists a  Kripke system $\bold K=(K, \leq, M_k, \{V_k\}_{k\in K})$ a homomorphism $\psi:\A\to \mathfrak{F}_{\bold K}$ and an element
$k_0\in \bold K$ and $x\in V_{k_0}$ such that for every $p\in \Gamma$, if $\psi(p)=(f_k)$ then $f_{k_0}(x)=1$.
Define for $k\in K$, $R_i$ an atomic formula and $s\in {}^{\alpha}M_k$, $k\models R_i[s]$ iff $(\psi(R/\equiv))_k(s)=1.$
This defines the desired  model.
\item When $G$ is a strongly rich semigroup, or $G={}^{\alpha}\alpha$, 
we show that for any $\beta$, $\A=\Fr_{\beta}GPHA_{\alpha}$ has the interpolation property, that is if $a\in \Sg^{\A}X_1$ and $b\in \Sg^{\A}X_2,$ 
then there exists
$c\in \Sg^{\A}(X_1\cap X_2)$ such that $a\leq c\leq b$. When $G$ is the semigroup of all finite transformations
and $\rho\in {}^{\beta}\wp(\alpha)$ is dimension restricting, the algebra $\Fr_{\beta}^{\rho}(GPHA_{\alpha})$ 
can be shown to have the interpolation property in exactly 
the same manner. We use theorem \ref{main} for the former case, 
while we use its analogue for dimension restricted free algebras, namely,
theorem \ref{main2} for the latter.  
Assume that $\theta_1\in \Sg^{\A}X_1$ and $\theta_2\in \Sg^{\A}X_2$ such that $\theta_1\leq \theta_2$.
Let
$\Delta_0=\{\theta\in \Sg^{\A}(X_1\cap X_2): \theta_1\leq \theta\}.$
If for some $\theta\in \Delta_0$ we have $\theta\leq \theta_2$, then we are done.
Else $(\Delta_0, \{\theta_2\})$ is consistent. Extend this to a complete theory $(\Delta_2, \Gamma_2)$ in $\Sg^{\A}X_2$.
Consider $(\Delta, \Gamma)=(\Delta_2\cap \Sg^{\A}(X_1\cap X_2), \Gamma_2\cap \Sg^{\A}(X_1\cap X_2))$.
Then $(\Delta\cup \{\theta_1\}), \Gamma)$ is consistent. For otherwise, for some $F\in \Delta, \mu\in \Gamma,$ 
we would have $(F\land \theta_1)\to \mu$ and $\theta_1\to (F\to \mu)$, so $(F\to \mu)\in \Delta_0\subseteq \Delta_2$ which is impossible.
Now $(\Delta\cup \{\theta_1\}, \Gamma)$ $(\Delta_2,\Gamma_2)$ are consistent with $\Gamma\subseteq \Gamma_2$ and $(\Delta,\Gamma)$
complete in $\Sg^{\A}X_1\cap \Sg^{\A}X_2$. So by theorem \ref{main}, $(\Delta_2\cup \{\theta_1\}, \Gamma_2)$ 
is satisfiable at some world in some set algbra based on a Kripke 
system, hence consistent. 
But this contradicts that
$\theta_2\in \Gamma_2, $ and we are done.
\end{enumarab}
\end{demo}

\begin{theorem} The logic ${\mathfrak L}_G^{=}$ has the weak interpolation property.
\end{theorem}

\begin{demo}{Proof}  Assume that $\theta_1\in \Sg^{\Rd\A}X_1$ and $\theta_2\in \Sg^{\Rd\A}X_2$ such that $\theta_1\leq \theta_2$.
Let $\Delta_0=\{\theta\in \Sg^{\A}(X_1\cap X_2): \theta_1\leq \theta\}.$
If for some $\theta\in \Delta_0$ we have $\theta\leq \theta_2$, then we are done.
Else $(\Delta_0, \{\theta_2\})$ is consistent, hence $(\Delta_0\cap \Sg^{\Rd\A}X_2,\theta_2)$ is consistent.
 Extend this to a complete theory $(\Delta_2, \Gamma_2)$ in $\Sg^{\Rd\A}X_2$; this is possible since $\theta_2\in \Sg^{\Rd\A}X_2$.
Consider $(\Delta, \Gamma)=(\Delta_2\cap \Sg^{\A}(X_1\cap X_2), \Gamma_2\cap \Sg^{\A}(X_1\cap X_2))$. 
It is complete in the `common language',  that is, in $\Sg^{\A}(X_1\cap X_2)$.
Then $(\Delta\cup \{\theta_1\}), \Gamma)$ is consistent in $\Sg^{\Rd\A}X_1$ and  $(\Delta_2, \Gamma_2)$ 
is consistent in $\Sg^{\Rd\A}X_2$, and 
$\Gamma\subseteq \Gamma_2.$ 
Applying the previous theorem, we get $(\Delta_2\cup \{\theta_1\}, \Gamma_2)$ is satisfiable.
Let $\psi_1, \psi_2$ and $\psi$ and $k_0$ be as in the previous proof. Then 
$\psi\upharpoonright \Sg^{\Rd\A}X_1=\psi_1$ and $\psi\upharpoonright \Sg^{\Rd\A}X_2=\psi_2$.
But $\theta_1\in \Sg^{\Rd\A}X_1$, then $\psi_1(\theta_1)=\psi(\theta_2)$. 
Similarly, $\psi_2(\theta_2)=\psi(\theta_2).$  
So, it readily follows that  $(\psi(\theta_1))_{k_0}(Id)=1$ and $(\psi(\theta_2))_{k_0}(Id)=0$.
This contradicts that $\psi(\theta_1)\leq \psi(\theta_2),$
and we are done.
\end{demo}

\begin{theorem} When $G$ consists only of finite transformations and $v\sim \rho(p)$ is infinite, then ${\mathfrak L}_G^{=}$ 
has the interpolation property.
\end{theorem}

In the next example, we show that the condition $\Gamma_0\subseteq \Gamma_0^*$ cannot be omitted. The example is an 
an algebraic version of theorem 4.31, p.121 in \cite{b}, but modified appropriately to  deal with infinitary languages.
\begin{example}\label{counter}

Let $G$ be a strongly rich semigroup on $\omega$. Let $\Lambda_{\omega}$ 
be a language with three predicate symbols each of arity $\omega$; this is a typless logic abstracting away from rank of atomic formulas, 
so that we might as well forget about the variables, since we allow them only in their natural order.
The real rank of such relation symbols will be recovered from the semantics.
Let $\bold M=(\N, \leq, \D_i)_{i\in \omega}$ be the Kripke frame with $D_i=\N$ for every $i$, 
and let $\N=\bigcup_{n\in \omega} B_n$, where $B_n$ is a sequence of pairwise disjoint infinite sets.
We define the relation $\models$ on atomic formulas. Let $m\in \N$. 
If $m=2n+1$, and $s\in {}^{\omega}\N,$ then $m\models p_0[s]$ if 
$s_0\in \bigcup_{i\leq 2n+1}B_i$, $m\models p_1[s]$ if $s_0\in \bigcup_{i\leq 2n+1}B_i$
and $m\models p_3[s]$. 

If $m=2n$, and $s\in {}^{\omega}\N,$ then $m\models p_0[s]$ if $s_0\in \bigcup_{i\leq n}B_i$ 
and $m\models p_2[s]$ if $s_0\in \bigcup_{i\leq 2n+1}B_i$ and $m\models p_3[s]$.
Let $\F_{\bold M}$ be the set algebra based on the defined above Kripke model $\bold M$.
Let $\A=\Fr_3GPHA_{\omega}$ and let $x_1, x_2, x_3$ be its generators.
Let $f$ be the unique map from $\A$ to $\F_{\bold M}$ 
such that for $i\in \{0,1,2\}$, $f(x_i)=p_i^{\bold M}$.
We have $\A\cong \Fm/\equiv$. We can assume that the isomorphism is the identity map.
Let $\Delta'=\{a\in A: f(a)=1\}$ and $\Theta'=\{a\in A: f(a)=0\}$.
Let $\Delta=\{\phi:\phi/\equiv\in \Delta'\},$ 
and $\Theta=\{\phi:\phi/\equiv\in \Gamma'\}.$
Let  
$$\Delta_1=\Delta\cup \{{\sf q}_0(x_1\lor x_2), {\sf c}_0(x_2\land x_3)\}$$
$$\Theta_1=\Theta\cup \{{\sf c}_0(x_1\land x_3)\}$$
$$\Delta_2=\Delta\cup \{{\sf q}_0(x_1\lor x_3)\}$$
$$\Theta_2=\Theta\cup \{{\sf c}_0(x_1\land x_3), c_0(x_2\land x_3)\}.$$
Then by analogy to 4.30 in \cite{b}, $(\Delta_1, \Theta_1)$, $(\Delta_2, \Theta_2)$ are consistent, but their union is not.
\end{example}

\begin{theorem}\label{mak} 
If $G$ is strongly rich or $G={}^{\alpha}\alpha$, then $Var(\mathfrak{L}_G)$ has $SUPAP.$
In particular, $GPHA_{\alpha}$ has $SUPAP$.
\end{theorem}

\begin{demo}{Proof} Cf. \cite{b} p.174.  Suppose that $\A_0, \A_1, \A_2\in Var(\mathfrak{L}_G)$.
Let $i_1:\A_0\to \A_2$ and $i_2: \A_0\to \A_2$ be embeddings. We need to find an amalgam.
We assume that $A_0\subseteq A_1\cap A_2$. For any $a\in A_i$, let $x_a^i$ be a variable such that $x_a^0=x_a^1=x_a^2$ for all $a\in A_0$ 
and the rest of the variables are distinct. Let $V_i$ be the set of variables corresponding to $\A_i$; then $|V_i|=|A_i|$. 
Let $V$ be the set of all  variables, endowed with countably infinitely many if the algebras are finite.
Then $|V|=\beta\geq \omega.$ We assume that the set of  variables $V$ of ${\mathfrak L}_G$ 
is the same as the set variables of 
the equational theory of $Var(\mathfrak{L}_G).$

We fix an assignment $s_i$ for each $i\in \{0,1,2\}$ such that 
$s_i: V_i\to A_i$ and $s_i(x_a^i)=a$ and so $s_1\upharpoonright V_0=s_2\upharpoonright V_0=s_0$. 
In view of the correspondence established in \ref{terms}, 
we identify terms of the equational theory of $Var({\mathfrak L}_G)$ with formulas of $\mathfrak{L}_G$;
which one we intend will be clear from context. 
Accordingly, we write $\A_i\models \psi\leftrightarrow \phi$ if $\bar{s}_i(\psi)=\bar{s}_i(\phi)$, where $\bar{s_i}$ is the unique extension of $s_i$ to 
the set all terms.
Let $\Fm_i$ be the set of formulas of $\mathfrak{L}_{G}$ in the variables $x_a^i$, $a\in A_i$, 
and let $\Fm$ be the set of all formulas built up from the set of all variables.
(Note that $Fm_i$ can be viewed as the set of terms built up from the variables $x_a^i$, 
and $Fm$ is the set of all terms built up from the set of all variables, defining operations corresponding to 
connectives turn them to absolutely free algebras.)

For $i=1,2$, let $T_i=\{\psi\in \Fm_i: \A_i\models \psi=1\}$, and let $T=\{\psi\in \Fm: T_1\cup T_2\vdash \psi\}$.

We will first prove (*):

For $\{i,j\}=\{1,2\},$ $\psi\in \Fm_i$ and  $\phi\in \Fm_j,$ we have
$T\vdash \psi\leftrightarrow \phi$ iff $(\exists c\in \Fm_0)(\A_i\models \psi\leq c\land \A_j\models c\leq \phi.)$

Only one direction is non trivial.
Assume that $T\vdash \psi\leftrightarrow \phi.$
Then there exist finite subsets $\Gamma_i\subseteq T_i$ and $\Gamma_j\subseteq T_j$ such that 
$\Gamma_i\cup  \Gamma_j\vdash \psi\leftrightarrow \phi.$
Then, by the deduction theorem for propositional intuitinistic logics, we get
$${\mathfrak L}_G\vdash \bigwedge \Gamma_i\to (\bigwedge \Gamma_j\to (\psi\to \phi)),$$
and so 
$${\mathfrak L}_G\vdash (\bigwedge \Gamma_i\land \psi)\to (\bigwedge \Gamma_j\to \phi).$$
Notice that atomic formulas and variables occuring in the last deduction are finite. So the interpolation theorem formulated for $G$ countable algebras
apply also, and indeed by this  interpolation theorem \ref{interpolation} for ${\mathfrak L}_G$, there is a formula $c\in \Fm_0$ such that 
such that $\vdash \bigwedge \Gamma_i\land \psi\to c$ and $\vdash c\to (\bigwedge 
\Gamma_j\to \phi.)$ Thus $\A_i\models \psi\leq c$ and $\A_j\models c\leq \phi$.
We have proved (*).

Putting $\psi=1,$ we get
$T\vdash \phi$ iff ($\exists c\in \Fm_0)(\A_i\models 1\leq c\land \A_j\models c\leq \phi)$ iff $\A_j\models \phi=1.$
Define on $\Fm$ the relation $\psi\sim \phi$ iff $T\vdash \alpha\leftrightarrow \beta$. Then $\sim$ is a congruence on $\Fm$.
Also for $i=1,2$ and $\psi,\phi\in Fm_i$, we have
$T\vdash \psi\sim \phi$ iff $\A_i\models \psi=\phi$.
Let $\A=\Fm/\sim$, and $e_i=\A_i\to \A$ be defined by
$e_i(a)=x_a^i/\sim$.
Then clearly $e_i$ is one to one. If $a\in \A_0$, then $x_a^0=x_a^1=x_a^2$ hence $e_1(a)=e_2(a)$.
Thus $\A$ is an amalgam via $e_1$ and $e_2.$
We now show that the superamalgamation property holds. Suppose $\{j,k\}=\{1,2\}$, $a\in \A_j$, $b\in \A_k$ and $e_j(a)\leq e_k(b)$. 
Then  $(e_j(a)\to e_k(b))=1$, so $(x_a^j\to x_b^k)=1$, that is $T\vdash (x_a^j\to x_b^k)$. Hence there exists 
$c\in \Fm_0$ such that $(\A_j\models x_a\leq c\land \A_k\models c\leq x_b)$. Then $a\leq c$ and $c\leq b.$

By taking ${\mathfrak L}_G$ to be the logic based on $\alpha$ many variables, and ${\mathfrak L}_G$ has 
countably many atomic formulas each containing $\alpha$ many variables in 
their natural order, we get that
$V=Var({\mathfrak L}_G)$, hence $V$ has $SUPAP$. 
\end{demo}

\end{document}